\documentclass[a4paper,12pt,twoside]{article}
\usepackage{amsmath,amssymb,ifthen,url}
\usepackage[amsmath,thmmarks]{ntheorem}
\usepackage{paper}
\usepackage[all]{xy}
\usepackage{mathrsfs}
\SelectTips{cm}{11}

\newcommand{\rig}{\mathrm{rig}}
\newcommand{\ad}{\mathrm{ad}}
\newcommand{\pg}{\mathrm{pg}}

\DeclareMathOperator{\Ind}{Ind}
\DeclareMathOperator{\Spf}{Spf}
\DeclareMathOperator{\Spa}{Spa}
\DeclareMathOperator{\spp}{sp}
\DeclareMathOperator{\Sh}{Sh}
\DeclareMathOperator{\Res}{Res}
\DeclareMathOperator{\Lie}{Lie}
\DeclareMathOperator{\Spl}{Spl}
\DeclareMathOperator{\GSp}{GSp}

\renewcommand{\theenumi}{(\roman{enumi})}

\renewcommand{\Box}{S}

\title{Potentially good reduction loci\\ of Shimura varieties}
\author{Naoki Imai and Yoichi Mieda}
\def\headertitle{Potentially good reduction loci of Shimura varieties}

\begin{document}

\maketitle

\begin{firstfootnote}
\end{firstfootnote}

\begin{abstract}
In this paper, we give a notion of the potentially good reduction locus of a Shimura variety.
It consists of the points which should be related with motives having potentially good reductions in some sense.
We show the existence of such locus for a Shimura variety of preabelian type. 
Further, we construct a partition of the adic space associated to a Shimura variety of preabelian type, 
which is expected to describe degenerations of motives. Using this partition, we prove that the cohomology of 
the potentially good reduction locus is isomorphic to the cohomology of a Shimura variety 
up to non-supercuspidal parts. 
\end{abstract}

\section{Introduction}
Let $(G,X)$ be a Shimura datum, and $\Sh_K(G,X)$ the Shimura variety attached to $(G,X)$
and a compact open subgroup $K$ of $G(\A^\infty)$.
It is known to be defined over a number field $E$, called the reflex field, 
which is canonically determined by $(G,X)$.
We fix a prime number $p$ and a place $v$ of $E$ above $p$, and write $E_v$ for the completion of $E$ at $v$.
The main theme of this article is the ``potentially good reduction locus'' of $\Sh_K(G,X)_{E_v}=\Sh_K(G,X)\otimes_EE_v$.

To explain what this locus is, let us first assume that $(G,X)$ is of PEL type, in which case
$\Sh_K(G,X)$ parametrizes abelian varieties with additional PEL structures.
We denote by $\mathcal{A}$ the universal abelian scheme over $\Sh_K(G,X)$.
If moreover the PEL datum is unramified at $p$ and $K=K_{p,0}K^p$ where $K_{p,0}$ is hyperspecial,
by extending the moduli problem to $\mathcal{O}_{E_v}$, we can obtain a good integral model $\mathscr{S}_{K^p}$
of $\Sh_K(G,X)$ over $\mathcal{O}_{E_v}$ (see \cite{MR1124982}).
This model is quite important in the study of the $\ell$-adic cohomology of Shimura varieties;
see \cite{MR1163241} for instance.
Let us denote by $\mathscr{S}^\wedge_{K^p}$ the formal completion of $\mathscr{S}_{K^p}$ along its special fiber,
and by $\mathscr{S}^{\wedge\rig}_{K^p}$ the rigid generic fiber of it.
Then, $\mathscr{S}^{\wedge\rig}_{K^p}$ is naturally identified with a quasi-compact rigid-analytic
open subset of $\Sh_K(G,X)_{E_v}$.
For a finite extension $F$ of $E_v$, an $F$-valued point $x$ of $\Sh_K(G,X)_{E_v}$ lies in $\mathscr{S}^{\wedge\rig}_{K^p}$
if and only if the abelian variety $\mathcal{A}_x$ over $F$ has (potentially) good reduction.
In this sense, $\mathscr{S}^{\wedge\rig}_{K^p}$ can be considered as the locus over which $\mathcal{A}$ has
(potentially) good reduction.
By this reason, we will write $\Sh_K(G,X)_{E_v}^{\mathrm{pg}}=\mathscr{S}^{\wedge\rig}_{K^p}$, and call it
the potentially good reduction locus.
We also have a rigid-analytic open subspace $\Sh_K(G,X)_{E_v}^{\mathrm{pg}}$ of $\Sh_K(G,X)_{E_v}$
for a compact open subgroup $K$ whose $p$-part is smaller than $K_{p,0}$, by taking the inverse image.
The $\ell$-adic cohomology of $\Sh_K(G,X)_{E_v}^{\mathrm{pg}}$ can be computed by using the cohomology of
the nearby cycle complex, provided that $\Sh_K(G,X)_{E_v}$ has a suitable integral model over $\mathcal{O}_{E_v}$.

In this paper, we will introduce the notion of the potentially good reduction locus for a general Shimura variety.
The rough idea is as follows. Let us fix a prime number $\ell$.
In the PEL type case, let $\mathcal{L}$ be the $\ell$-adic automorphic \'etale sheaf on $\Sh_K(G,X)$
attached to the standard representation of $G$. Then, for a finite extension $F$ of $E_v$ and an $F$-valued point $x$
of $\Sh_K(G,X)$, the stalk $\mathcal{L}_{\overline{x}}$ can be identified with the rational $\ell$-adic Tate module
$V_\ell\mathcal{A}_{\overline{x}}$. By the Neron-Ogg-Shafarevich criterion, we conclude that $\mathcal{A}_x$ has
potentially good reduction if and only if $\mathcal{L}_{\overline{x}}$ is potentially unramified
(resp.\ potentially crystalline) when $\ell\neq p$ (resp.\ $\ell=p$).
This observation urges us to define in the general case that an $F$-valued point $x$ of $\Sh_K(G,X)_{E_v}$
is of potentially good reduction if $\mathcal{L}_{\overline{x}}$ is potentially unramified/crystalline
for every automorphic \'etale sheaf $\mathcal{L}$.
Actually in the paper, we look at the torsor over $x$ obtained as the pull-back of
$\varprojlim_{K'\subset K}\Sh_{K'}(G,X)_{E_v}\to \Sh_K(G,X)_{E_v}$, which is more concise
but essentially equivalent to the above way by the Tannakian duality.
Our potentially good reduction locus is defined as
a quasi-compact open subset of the adic space $\Sh_K(G,X)_{E_v}^\ad$ attached to $\Sh_K(G,X)_{E_v}$,
whose $F$-valued points consist of those of potentially good reduction for every $F$.
It is unique, if exists.
We will show the existence of the potentially good reduction locus $\Sh_K(G,X)_{E_v}^{\mathrm{pg}}$ when the Shimura datum
$(G,X)$ is of preabelian type. Recall that $(G,X)$ is said to be of preabelian type if there exists a Shimura datum $(G',X')$
of Hodge type such that $(G^\ad,X^\ad)\cong (G'^\ad,X'^\ad)$.
This class contains almost all Shimura data in practice.
As in \cite[Introduction]{MR546620} and \cite[\S 9]{MR2192012},
a Shimura variety is believed to have a moduli interpretation by motives,
if the weight homomorphism for $(G,X)$ is defined over the rational number field.
The subset $\Sh_K(G,X)^{\mathrm{pg}}_{E_v}$ is expected to parametrize motives
with potentially good reduction at $v$.

We are also interested in what happens outside the locus $\Sh_K(G,X)^{\mathrm{pg}}_{E_v}$.
In the PEL type case, degenerations of abelian varieties occur; if a Shimura variety parametrizes motives,
then degenerations of motives should occur.
Based on this observation, we will construct a partition of $\Sh_K(G,X)^\ad_{E_v}$ into
finitely many locally closed constructible subsets
labeled by conjugacy classes of certain kind of adelic parabolic subgroups of $G$,
so that the piece corresponding to $G$ equals $\Sh_K(G,X)^{\mathrm{pg}}_{E_v}$.
It is closely related to the theory of integral toroidal compactifications.
Actually, in the PEL type case, we may also use the integral toroidal compactification 
developed in \cite{Kai-Wen} to construct our partition 
(see \cite[\S 7]{MR3221715}); there should be some more cases to which the method in \cite[\S 7]{MR3221715}
can be applied (for example, \cite{Mad-torHod}).
However, our argument here is almost totally rigid-geometric, and requires only the existence of
the integral toroidal compactifications of the Siegel modular varieties with hyperspecial level at $p$
(\cite{MR1083353}) as an input from the integral theory.
Note also that our partition is independent of any choice,
unlike the toroidal compactification that depends on the choice of a cone decomposition.

By using the partition above, we can compare the $\ell$-adic cohomology of the tower $\{\Sh_K(G,X)^{\mathrm{pg}}_{E_v}\}_K$
and that of $\{\Sh_K(G,X)\}_K$.
We assume that $(G,V)$ is of preabelian type and satisfies the condition SV6 in \cite[p.~311]{MR2192012}.
Let $\ell$ be a prime number different from $p$, and $\mathcal{L}$ an $\ell$-adic
automorphic \'etale sheaf on $\Sh_K(G,X)$ corresponding to an algebraic representation of $G^c$
over $\overline{\Q}_\ell$, where $G^c$ is the quotient of $G$ defined in \cite[p.~347]{MR1044823}.
The statement is as follows:

\begin{thm}[Theorem \ref{thm:main-thm}]\label{thm:main-thm-intro}
 In the kernel and the cokernel of the natural map
 \[
 \varinjlim_{K}H^i_c(\Sh_{K}(G,X)_{\overline{E}_v}^{\pg},\mathcal{L})\to
 \varinjlim_{K}H^i_c(\Sh_{K}(G,X)_{\overline{E}_v},\mathcal{L}),
 \]
 no irreducible supercuspidal representation of $G(\Q_{p'})$ appears as a subquotient for any prime number $p'$.
\end{thm}
Recall that an irreducible smooth representation of $G(\Q_{p'})$ is said to be supercuspidal
if it does not appear as a subquotient of the parabolically induced representations
from any proper parabolic subgroup.
Loosely speaking, this theorem is a consequence of the observation that the partition of
the complement $\Sh_{K}(G,X)^\ad_{E_v}\setminus\Sh_{K}(G,X)_{E_v}^{\pg}$ is ``geometrically induced''
from proper parabolic subgroups of $G(\Q_{p'})$.
It will be worth noting that our method is totally geometric, so that it is also valid
in the torsion coefficient case.
See Theorem \ref{thm:main-thm-mod-l} for an analogue of Theorem \ref{thm:main-thm-intro} with
the $\overline{\F}_\ell$-coefficients. 

We have already mentioned that in the PEL type case the $\ell$-adic cohomology $H^i_c(\Sh_{K}(G,X)_{\overline{E}_v}^{\pg},\mathcal{L})$
can be computed as the cohomology of a nearby cycle complex.
Hence, in this case the theorem above says that the nearby cycle cohomology is isomorphic to
the compactly supported cohomology up to non-supercuspidal representations.
This result is useful, since it connects the cohomology of Shimura varieties and that of Rapoport-Zink spaces;
see Section \ref{subsec:example} for a simple example in this direction.
Recently, during the preparation of this article, Lan and Stroh obtained a stronger result
that the nearby cycle cohomology is isomorphic to the compactly supported cohomology
in the cases where reasonable integral toroidal compactifications exist (see \cite{LSNcycle}). 
However, we have decided to include our weaker result in this paper, since the argument is totally different. 

We sketch the outline of this paper. 
In Section \ref{sec:Gal-rep}, we consider Galois representations of a $p$-adic field with values in
a general connected reductive group $G$. Under some condition, we attach a parabolic subgroup of $G$ to such a representation.
In Section \ref{sec:rig-semi-ab}, we give some preliminary results on adic spaces and semi-abelian schemes. 
In Section \ref{sec:Shvar}, we recall some notation and results on Shimura varieties. 
In Section \ref{sec:partition}, we construct a partition of the adic space associated to 
a Shimura variety of preabelian type by using results obtained in Section \ref{sec:Gal-rep}.
The potentially good reduction locus is introduced here, as a piece of the constructed partition.
In Section \ref{sec:cohomology}, we prove the theorem comparing the cohomology of potentially good reduction loci 
with that of Shimura varieties. 
In Section \ref{sec:PELcase}, we specialize our results to Shimura varieties of PEL type, 
and discuss a simple application. 

\bigbreak

\noindent{\bfseries Acknowledgment}\quad
The authors are grateful to Kai-Wen Lan, who kindly informed them the result in \cite{LSNcycle}.
This work was supported by JSPS KAKENHI Grant Numbers 24740019, 26610003, 15H03605.

\bigbreak

\noindent{\bfseries Notation}\quad
Put $\widehat{\Z}=\prod_{\text{prime $p$}}\Z_p$ and $\A^\infty=\widehat{\Z}\otimes_\Z\Q$.
For a prime $p$, put $\widehat{\Z}^p=\prod_{\text{prime $p'\neq p$}}\Z_{p'}$ and 
$\A^{\infty,p}=\widehat{\Z}^p\otimes_\Z\Q$.
More generally, for a finite set of primes $S$, we put $\A_S=\prod_{\ell\in S}\Q_\ell$, 
$\widehat{\Z}^S=\prod_{\text{prime $p\notin S$}}\Z_p$ and
${\A}^{\infty,S}=\widehat{\Z}^S\otimes_\Z\Q$.

For a scheme $X$ over a field $F$ and an extension field $L$ of $F$, we write $X_L$ for the base change of $X$
to $L$. Similar notation will be used for adic spaces.

For an algebraic group $G$, let $Z(G)$ denote the center of $G$, and $G^{\ad}=G/Z(G)$ the adjoint group of $G$. 
For a field $L$ over which $G$ is defined, we write $\mathbf{Rep}_L(G)$ for the Tannakian category of
finite-dimensional algebraic representations of $G$ over $L$.

Every sheaf and cohomology are considered in the \'etale topology.

\section{Preliminaries on Galois representations}\label{sec:Gal-rep}
In this section, fix a $p$-adic field $F$ and its algebraic closure $\overline{F}$.
Let $\ell$ be a prime number and $G$ a connected reductive group over $\overline{\Q}_\ell$.
Consider a continuous homomorphism $\phi\colon \Gal(\overline{F}/F)\to G(\overline{\Q}_\ell)$.

\begin{defn}\label{defn:phipro}
\begin{enumerate}
 \item Assume that $\ell\neq p$. 
       We say that $\phi$ is potentially unramified if $\xi\circ\phi$ is potentially unramified 
       for any $\xi\in \mathbf{Rep}_{\overline{\Q}_\ell}(G)$. 
 \item Assume that $\ell=p$.
       We say that $\phi$ is potentially crystalline if $\xi\circ\phi$ is potentially crystalline
       for any $\xi\in \mathbf{Rep}_{\overline{\Q}_p}(G)$. 
 \item Assume that $\ell=p$. 
       We say that $\phi$ is de Rham if $\xi \circ \phi$ is de Rham (or equivalently, potentially semistable)
       for any $\xi\in \mathbf{Rep}_{\overline{\Q}_p}(G)$.
\end{enumerate}
\end{defn}

To measure how far $\phi$ is from potentially unramified or potentially crystalline,
we consider the monodromy filtration on $\xi\circ\phi$ for each $\xi\in\mathbf{Rep}_{\overline{\Q}_\ell}(G)$.
First, we assume that $\ell\neq p$. 
Then, for each $(\xi,V_\xi)\in\mathbf{Rep}_{\overline{\Q}_\ell}(G)$,
we obtain the $\ell$-adic representation $(\xi\circ\phi,V_\xi)$ of $\Gal(\overline{F}/F)$ and
its monodromy filtration $M_\bullet V_\xi$.

\begin{lem}\label{lem:monodromy-filt-stab-l}
 Assume that $\ell\neq p$.
 \begin{enumerate}
  \item\label{item:P-parab-l}
       The stabilizer $P_\xi$ of the filtration $M_\bullet V_\xi\subset V_\xi$ is a parabolic subgroup of $G$.
  \item\label{item:P-faithful-l}
       If $\xi\in\mathbf{Rep}_{\overline{\Q}_\ell}(G)$ is faithful, then $P_{\xi}$ stabilizes $M_\bullet V_{\xi'}$ for every $\xi'\in\mathbf{Rep}_{\overline{\Q}_\ell}(G)$.
       In particular $P_{\xi}$ for faithful $\xi$ is independent of $\xi$.
       We write $P_\phi$ for this $P_\xi$.
  \item\label{item:P=G-l}
       The homomorphism $\phi$ is potentially unramified if and only if $P_\phi=G$.
  \item\label{item:P-ad-l}
       For the composite $\phi^\ad\colon \Gal(\overline{F}/F)\xrightarrow{\phi}G(\overline{\Q}_\ell)\to G^\ad(\overline{\Q}_\ell)$, we have $P_{\phi^\ad}=P^\ad_\phi$, where $P^\ad_\phi$ denotes the image of $P_\phi$ in $G^\ad$.
  \item\label{item:P-fin-ext-l}
       For a finite extension $F'$ of $F$ contained in $\overline{F}$, we put
       $\phi'=\phi\vert_{\Gal(\overline{F}/F')}$. Then we have $P_{\phi'}=P_\phi$.
 \end{enumerate}
\end{lem}

\begin{prf}
 The assertion \ref{item:P-parab-l} follows from \cite[Lemma 1.1.1, Lemma 1.1.3]{MR2669706},
 since $V_\xi\mapsto M_\bullet V_\xi$ gives a filtration on the Tannakian category $\mathbf{Rep}_{\overline{\Q}_\ell}(G)$.

 Let us prove \ref{item:P-faithful-l}. 
 For integers $m,m'\ge 0$, the monodromy filtration on $V_\xi^{\otimes m}\otimes V_\xi^{\vee\otimes m'}$
 can be written by using $M_\bullet V_\xi$ (see \cite[Proposition 1.6.9]{MR601520}).
 Therefore it is stable under $P_\xi$.
 As every representation $(\xi',V_{\xi'})$ of $G$ appears
 as a direct summand of $V_\xi^{\otimes m}\otimes V_\xi^{\vee\otimes m'}$ for some integers $m,m'\ge 0$,
 the filtration $M_\bullet V_{\xi'}$ is also preserved by $P_\xi$. This concludes the proof.

 For \ref{item:P=G-l}, note that $\phi$ is potentially unramified if and only if
 the monodromy operator $N$ on $V_\xi$ is zero for every $\xi\in\mathbf{Rep}_{\overline{\Q}_\ell}(G)$.
 If this condition is satisfied, we have $M_iV_\xi=0$ for $i<0$ and $M_iV_\xi=V_\xi$ for $i\ge 0$,
 thus $P_\phi=G$.
 Conversely assume that $P_\phi=G$, in other words, $M_iV_\xi$ is $G$-stable for every $\xi$ and $i$.
 It suffices to show that $N=0$ on $V_\xi$ for each irreducible representation $\xi$ of $G$.
 Since $M_iV_\xi$ is $G$-stable, there exists a unique integer $i_0$ such that $M_iV_\xi=0$ $(i<i_0)$ 
 and $M_iV_\xi=V_\xi$ $(i\ge i_0)$.
 Hence we have $N(V_\xi)=N(M_{i_0}V_\xi)\subset M_{i_0-2}V_\xi=0$, as desired.
 
 We prove \ref{item:P-ad-l}. Clearly we have $P_{\phi}^\ad\subset P_{\phi^\ad}$.
 For the reverse inclusion, take $g\in G(\overline{\Q}_\ell)$ which is mapped into $P_{\phi^\ad}(\overline{\Q}_\ell)$
 under $G(\overline{\Q}_\ell)\to G^\ad(\overline{\Q}_\ell)$. It suffices to show that
 $g$ stabilizes $M_\bullet V_\xi$ for each irreducible representation $(\xi,V_\xi)$ of $G$.
 Put $W=V_\xi\otimes V_\xi^\vee$. Since $\xi$ is irreducible, the center of $G$ acts trivially on $W$.
 Therefore $W$ can be regarded as a representation of $G^\ad$. In particular, the monodromy filtration
 $M_\bullet W$ on $W$ is stable under the action of $g$.

 Let $j_0$ be the minimal integer such that $M_{j_0}V_\xi^\vee=V_\xi^\vee$.
 Fix an integer $i_0$. Then we have 
 $M_{i_0+j_0}W=\sum_{i+j=i_0+j_0}M_iV_\xi\otimes M_jV_\xi^\vee$. Note that $M_{i_0}V_\xi$ can be recovered
 from $M_{i_0+j_0}W$ by
 \begin{equation*}
  M_{i_0}V_\xi=\bigcap_{f\in \Hom_{\overline{\Q}_\ell}(V_\xi^\vee,\overline{\Q}_\ell)\setminus \{0\}}(\id\otimes f)(M_{i_0+j_0}W).\tag{$*$}
 \end{equation*}
 Since $M_{i_0+j_0}W=g(M_{i_0+j_0}W)=\sum_{i+j=i_0+j_0}g(M_iV_\xi)\otimes g(M_jV_\xi^\vee)$, the right hand side
 of $(*)$ is also equal to $g(M_{i_0}V_\xi)$. Hence we conclude that $M_\bullet V_\xi$ is stable under $g$.

 The claim \ref{item:P-fin-ext-l} is clear, since the monodromy filtration $M_\bullet V_\xi$ does not change
 after restricting $\phi$ to $\Gal(\overline{F}/F')$.
\end{prf}

Next we consider the case $\ell=p$.
We shall introduce the notion of the monodromy filtration on a $p$-adic Galois representation.
Let $L$ be a finite extension of $\Q_p$ and $V$ a finite-dimensional de Rham $L$-representation
of $\Gal(\overline{F}/F)$. We regard $V$ as a $\Q_p$-representation of $\Gal(\overline{F}/F)$
and consider $D_{\mathrm{pst}}(V)$,
where $D_{\mathrm{pst}}$ is the functor introduced in \cite[\S 5.6]{MR1293972}.
If we write $\Q_p^{\mathrm{ur}}$ for the maximal unramified extension of $\Q_p$ contained in $\overline{F}$,
$D_{\mathrm{pst}}(V)$ is an $L\otimes_{\Q_p}\Q_p^{\mathrm{ur}}$-module equipped with several structures.
Among them, we have the monodromy operator on $D_{\mathrm{pst}}(V)$,
from which the monodromy filtration $M_\bullet D_{\mathrm{pst}}(V)$ on $D_{\mathrm{pst}}(V)$ is
naturally induced.

\begin{defn}\label{defn:has-monodromy-filt-rep-fin}
 Let $V$ be a de Rham $L$-representation of $\Gal(\overline{F}/F)$.
 We say that $V$ has the monodromy filtration if there exists a $\Gal(\overline{F}/F)$-stable $\Q_p$-subspace
 $M_iV$ of $V$ such that $D_{\mathrm{pst}}(M_iV)=M_iD_{\mathrm{pst}}(V)$ for each $i$.
 Such a subspace $M_iV$ is unique and stable under the action of $L$ if it exists,
 thanks to the fact that $D_{\mathrm{pst}}$ is fully faithful.
\end{defn}

Now, let $V$ be a finite-dimensional de Rham $\overline{\Q}_p$-representation of $\Gal(\overline{F}/F)$.
We can find a subfield $L$ of $\overline{\Q}_p$ which is finite over $\Q_p$
and a $\Gal(\overline{F}/F)$-stable $L$-subspace $V_L$ of $V$ such that $V_L\otimes_L\overline{\Q}_p=V$.
The $L$-representation $V_L$ is de Rham.

\begin{defn}\label{defn:has-monodromy-filt-rep}
 The condition that $V_L$ has the monodromy filtration is independent of the choice of $L$ and $V_L$.
 If it is the case, we say that $V$ has the monodromy filtration, and put
 $M_iV=M_iV_L\otimes_L\overline{\Q}_p$, which is easily seen to be independent of $L$ and $V_L$.
 We call $M_\bullet V$ the monodromy filtration of $V$.
\end{defn}

\begin{lem}\label{lem:monodromy-filt-properties}
 Let $V$, $W$ be finite-dimensional $\overline{\Q}_p$-representations of $\Gal(\overline{F}/F)$
 which are de Rham and have the monodromy filtrations.
 \begin{enumerate}
  \item\label{item:monod-summand} Let $V'$ be a direct summand of $V$ as a $\Gal(\overline{F}/F)$-representation.
	Then $V'$ is de Rham and has the monodromy filtration.
  \item\label{item:monod-tensor} The representations $V\otimes W$ and $V^\vee$ are de Rham and have the monodromy filtrations.
       The monodromy filtration $M_\bullet (V\otimes W)$ (resp.\ $M_\bullet (V^\vee)$) is given by
       $M_n(V\otimes W)=\sum_{i+j=n}M_iV\otimes M_jW$ (resp.\ $M_n(V^\vee)=(V/M_{-n-1}V)^\vee$).
 \end{enumerate}
\end{lem}

\begin{prf}
 Let $L$ be a finite extension of $\Q_p$.
 We have only to consider $L$-representations in place of $\overline{\Q}_p$-representations.
 In the following, let $V$, $W$ be finite-dimensional $L$-representations
 of $\Gal(\overline{F}/F)$ which are de Rham and have the monodromy filtrations.

 We prove \ref{item:monod-summand}. We write $V=V'\oplus V''$. 
 Then we have $D_{\mathrm{pst}}(V)=D_{\mathrm{pst}}(V')\oplus D_{\mathrm{pst}}(V'')$ and
 $M_iD_{\mathrm{pst}}(V)=M_iD_{\mathrm{pst}}(V')\oplus M_iD_{\mathrm{pst}}(V'')$.
 Since the essential image of the functor $D_{\mathrm{pst}}$ is stable under direct factors
 (see \cite[Th\'eor\`eme 5.6.7]{MR1293972}),
 there exists a $\Gal(\overline{F}/F)$-stable subspace $M_iV'$ of $V'$ such that 
 $D_{\mathrm{pst}}(M_iV')=M_iD_{\mathrm{pst}}(V')$. Hence $V'$ has the monodromy filtration.

 Next consider \ref{item:monod-tensor}. It is known that $V\otimes_{\Q_p}W$ is a de Rham representation, and 
 $D_{\mathrm{pst}}(V\otimes_{\Q_p}W)=D_{\mathrm{pst}}(V)\otimes_{\Q_p^{\mathrm{ur}}} D_{\mathrm{pst}}(W)$.
 The monodromy filtration on $D_{\mathrm{pst}}(V\otimes_{\Q_p}W)$ is given by
 \[
  M_nD_{\mathrm{pst}}(V\otimes_{\Q_p}W)=\sum_{i+j=n}M_iD_{\mathrm{pst}}(V)\otimes_{\Q_p^{\mathrm{ur}}} M_jD_{\mathrm{pst}}(W).
 \]
 Therefore, if we put $M_n(V\otimes_{\Q_p}W)=\sum_{i+j=n}M_iV\otimes_{\Q_p}M_jW$, we have 
 $D_{\mathrm{pst}}(M_n(V\otimes_{\Q_p}W))=M_nD_{\mathrm{pst}}(V\otimes_{\Q_p}W)$
 by the exactness of $D_{\mathrm{pst}}$.
 Hence $V\otimes_{\Q_p}W$ has the monodromy filtration.

 Let $e\in L\otimes_{\Q_p}L$ denote the idempotent corresponding to the diagonal component
 $\Spec L\hookrightarrow \Spec (L\otimes_{\Q_p}L)$. Then, we have $e(V\otimes_{\Q_p}W)=V\otimes_LW$.
 In particular $V\otimes_LW$ is a direct summand of $V\otimes_{\Q_p}W$.
 Therefore, by (i), $V\otimes_LW$ is de Rham and has the monodromy filtration.
 Clearly, the monodromy filtration on $V\otimes_LW=e(V\otimes_{\Q_p}W)$ is given by
 \[
  M_n(V\otimes_LW)=eM_n(V\otimes_{\Q_p}W)=\sum_{i+j=n}e(M_iV\otimes_{\Q_p}M_jW)=\sum_{i+j=n}M_iV\otimes_L M_jW.
 \]
 The dual $V^\vee$ can be treated similarly.
\end{prf}

\begin{lem}\label{lem:monod-restr-p-rep}
 Let $V$ be a finite-dimensional $\overline{\Q}_p$-representation of $\Gal(\overline{F}/F)$.
 For a finite extension $F'$ of $F$ contained in $\overline{F}$,
 we put $V'=V\vert_{\Gal(\overline{F}/F')}$.
 Then, $V$ is de Rham and has the monodromy filtration if and only if
 so is $V'$.
 Moreover we have $M_\bullet V'=(M_\bullet V)\vert_{\Gal(\overline{F}/F')}$.
\end{lem}

\begin{prf}
 As in the proof of Lemma \ref{lem:monodromy-filt-properties}, we may replace $V$ by an $L$-representation
 of $\Gal(\overline{F}/F)$, where $L$ is a finite extension of $\Q_p$.
 By definition, $V$ is de Rham if and only if $V'$ is de Rham.
 Suppose that $V$ and $V'$ are de Rham.
 We have $D_{\mathrm{pst}}(V')=D_{\mathrm{pst}}(V)\vert_{\Gal(\overline{F}/F')}$.
 Therefore, if $V$ has the monodromy filtration $M_\bullet V$,
 then $V'$ has the monodromy filtration $(M_\bullet V)\vert_{\Gal(\overline{F}/F')}$.
 Conversely, assume that $V'$ has the monodromy filtration $M_\bullet V'$.
 Since $D_{\mathrm{pst}}(M_iV')=M_iD_{\mathrm{pst}}(V')=M_iD_{\mathrm{pst}}(V)\subset D_{\mathrm{pst}}(V)$
 is stable under $\Gal(\overline{F}/F)$, so is $M_iV'\subset V'=V$.
 Therefore, $M_\bullet V'$ gives the monodromy filtration of $V$. This concludes the proof.
\end{prf}

Now, let $\phi\colon \Gal(\overline{F}/F)\to G(\overline{\Q}_p)$ be a continuous homomorphism, 
as in the beginning of this section.

\begin{defn}
 Assume that $\phi$ is de Rham. We say that $\phi$ has the monodromy filtration if
 $V_\xi$ has the monodromy filtration for every $\xi\in \mathbf{Rep}_{\overline{\Q}_p}(G)$.
\end{defn} 

\begin{lem}\label{lem:monodromy-filt-p}
 Assume that there exists a faithful algebraic representation $(\xi,V_\xi)$ of $G$ such that
 $V_\xi$ is de Rham and has the monodromy filtration. Then $\phi$ is de Rham and has the monodromy filtration.
\end{lem}

\begin{prf}
 Thanks to Lemma \ref{lem:monodromy-filt-properties}, we can use the same argument
 as in the proof of Lemma \ref{lem:monodromy-filt-stab-l} \ref{item:P-faithful-l}.
\end{prf}

\begin{lem}\label{lem:monodromy-filt-stab-p}
 Assume that $\phi$ is de Rham and has the monodromy filtration.
 \begin{enumerate}
  \item\label{item:P-parab-p}
       The stabilizer $P_\xi$ of the filtration $M_\bullet V_\xi\subset V_\xi$ is a parabolic subgroup of $G$.
  \item\label{item:P-faithful-p}
       If $\xi$ is faithful, then $P_{\xi}$ stabilizes $M_\bullet V_{\xi'}$ for every representation $\xi'$.
       In particular $P_{\xi}$ for faithful $\xi$ is independent of $\xi$.
       We write $P_\phi$ for this $P_\xi$.
  \item\label{item:P=G-p}
       The homomorphism $\phi$ is potentially crystalline if and only if $P_\phi=G$.
  \item\label{item:P-ad-p} 
       The composite $\phi^\ad\colon \Gal(\overline{F}/F)\xrightarrow{\phi}G(\overline{\Q}_p)\to G^\ad(\overline{\Q}_p)$ is de Rham and has the monodromy filtration. Moreover we have $P_{\phi^\ad}=P^\ad_\phi$.
 \end{enumerate}
\end{lem}

\begin{prf}
 This can be proved in the same way as Lemma \ref{lem:monodromy-filt-stab-l}.
\end{prf}

\begin{lem}\label{lem:monod-restr-p}
 For a finite extension $F'$ of $F$ contained in $\overline{F}$,
 $\phi'=\phi\vert_{\Gal(\overline{F}/F')}$ is de Rham and has the monodromy filtration
 if and only if so is $\phi$.
 Moreover, if the above conditions are satisfied, we have $P_{\phi'}=P_\phi$.
\end{lem}

\begin{prf}
 This is an immediate consequence of Lemma \ref{lem:monod-restr-p-rep}.
\end{prf}

\begin{cor}\label{cor:adjoint-equivalence}
 Let $\ell$ be a prime number.
 \begin{enumerate}
  \item Assume that $\ell \neq p$. 
	Then $\phi$ is potentially unramified if and only if $\phi^{\ad}$ is potentially unramified. 
  \item Assume that $\ell=p$, $\phi$ is de Rham and has the monodromy filtration. 
	Then $\phi$ is potentially crystalline if and only if $\phi^{\ad}$ is potentially crystalline. 
 \end{enumerate}
\end{cor}

\begin{prf}
 The first assertion follows from Lemma \ref{lem:monodromy-filt-stab-l} \ref{item:P=G-l}, \ref{item:P-ad-l},
 and the second from Lemma \ref{lem:monodromy-filt-stab-p} \ref{item:P=G-p}, \ref{item:P-ad-p}.
\end{prf}

\begin{rem}
 In fact, Corollary \ref{cor:adjoint-equivalence} (ii) holds without assuming that $\phi$ has
 the monodromy filtration 
 (we have only to consider the monodromy filtration on the image of $D_{\mathrm{pst}}$).
 However we do not need this fact later.
\end{rem}

\begin{rem}\label{rem:parab-rationality-l}
 Let $\ell$ be a prime number, and assume that $\phi$ is de Rham and has the monodromy filtration if $\ell=p$.
 If $G$ is defined over $\Q_\ell$ and the image of $\phi$ is contained in $G(\Q_\ell)$,
 the parabolic subgroup $P_\phi$ is defined over $\Q_\ell$.
 Indeed, we can take a faithful representation $\xi$ which is defined over $\Q_\ell$,
 and then the monodromy filtration $M_\bullet V_\xi$ is also defined over $\Q_\ell$.
\end{rem}

\section{Rigid geometry and semi-abelian schemes}\label{sec:rig-semi-ab}
\subsection{Notation for adic spaces}\label{subsec:notation-adic}
Throughout this paper, we will use the framework of adic spaces introduced by Huber
(\cf \cite{MR1207303}, \cite{MR1306024}, \cite{MR1734903}). Here we recall some notation briefly.

Let $S$ be a noetherian scheme and $S_0$ a closed subscheme of $S$. We denote the formal completion
of $S$ along $S_0$ by $\mathcal{S}$. Put $\mathcal{S}^\rig=t(\mathcal{S})_a$, where
$t(\mathcal{S})$ is the adic space associated to $\mathcal{S}$ (\cf \cite[\S 4]{MR1306024})
and $t(\mathcal{S})_a$ denotes the open adic subspace of $t(\mathcal{S})$
consisting of analytic points.
It is a quasi-compact analytic adic space.

Let $X$ be a scheme of finite type over $S$. Put $X_0=X\times_SS_0$ and denote
the formal completion of $X$ along $X_0$ by $\widehat{X}$. Then we can construct an adic space
$\widehat{X}^\rig$ in the same way as $\mathcal{S}^\rig$. The induced morphism
$\widehat{X}^\rig\to \mathcal{S}^\rig$ is of finite type.
On the other hand, we can construct another adic space $X\times_S\mathcal{S}^\rig$.
 Indeed, since we have morphisms of locally ringed spaces
$(\mathcal{S}^\rig,\mathcal{O}_{\mathcal{S}^\rig})\to (t(\mathcal{S}),\mathcal{O}_{t(\mathcal{S})})\to (S,\mathcal{O}_S)$ (for the second one, see \cite[Remark 4.6 (iv)]{MR1306024}), we can make the fiber product $X\times_S\mathcal{S}^\rig$ in the sense of
\cite[Proposition 3.8]{MR1306024}. 
For simplicity, we write $X^\ad$ for $X\times_S\mathcal{S}^\rig$, though it depends on
$(S,S_0)$.
Since the morphism $\mathcal{S}^\rig\to S$ factors through $S^0=S\setminus S_0$,
we have $X\times_S\mathcal{S}^\rig=(X\times_SS^0)\times_{S^0}\mathcal{S}^\rig$.
In particular, $X^\ad$ depends only on $X\times_SS^0$.
The natural morphism $X^\ad\to \mathcal{S}^\rig$ is locally of finite type,
but not necessarily quasi-compact; see the following example.

\begin{exa}
 Let $R$ be a complete discrete valuation ring and $F$ its fraction field.
 Consider the case where $S=\Spec R$ and $S_0$ is the closed point of $S$.
 Then, for an $S$-scheme $X$ of finite type, $X^\ad$ can be regarded as the rigid space over $F$
 associated to a scheme $X\times_S\Spec F$ over $F$.
 For example, $(\A^1_S)^\ad=(\A^1_F)^\ad$ is the rigid-analytic affine line over $F$ and thus is not
 quasi-compact. On the other hand, $(\widehat{\A}^1_S)^\rig$ is the unit disc
 ``$\lvert z\rvert\le 1$'' in $(\A^1_F)^\ad$, which is quasi-compact.
\end{exa}

\begin{lem}\label{lem:fiber-product}
 The functors $X\mapsto \widehat{X}^\rig$ and $X\mapsto X^\ad$ commute with
 fiber products.
\end{lem}

\begin{prf}
 For the functor $X\mapsto \widehat{X}^\rig$, it can be checked easily
 (\cf \cite[Lemma 3.4]{MR2211156} and \cite[Proof of Lemma 4.4 (v)]{formalnearby}).
 Consider the functor $X\mapsto X^\ad$. 
 Let $Y\to X\leftarrow Z$ be a diagram of $S$-schemes of finite type.
 What we should prove is
 \[
  (Y\times_XZ)\times_S\mathcal{S}^\rig\cong (Y\times_S\mathcal{S}^\rig)\times_{X\times_S\mathcal{S}^\rig}(Z\times_S\mathcal{S}^\rig).
 \]
 It is not totally automatic, since $Y\times_XZ$ is not a fiber product in the category of locally ringed spaces.
 It follows from the fact that morphisms of locally ringed spaces
 $\Spa (A,A^+)\to \Spec B$ for a complete affinoid ring $(A,A^+)$ and a ring $B$
 correspond bijectively to ring homomorphisms $B\to A$
 (this fact is used implicitly in \cite[Remark 4.6 (iv)]{MR1306024} to define $t(\mathcal{S})\to S$).
\end{prf}

Let us compare $\widehat{X}^\rig$ and $X^\ad$; by the commutative diagram
\[
 \xymatrix{%
 \widehat{X}^\rig\ar[r]\ar[d]& X\ar[d]\\
 \mathcal{S}^\rig\ar[r]& S
 }
\]
and the universality of the fiber product $X\times_S\mathcal{S}^\rig$, we have a natural
morphism $\widehat{X}^\rig\to X^\ad$.

\begin{lem}\label{lem:rig-ad}
 \begin{enumerate}
  \item If $X$ is separated over $S$, $\widehat{X}^\rig\to X^\ad$ is an open immersion.
  \item If $X$ is proper over $S$, $\widehat{X}^\rig\to X^\ad$ is an isomorphism.
 \end{enumerate}
\end{lem}

\begin{prf}
 See \cite[Remark 4.6 (iv)]{MR1306024}.
\end{prf}

\begin{rem}\label{rem:rig-ad-base-change}
 Let $f\colon S'\to S$ be a morphism of finite type and $S'_0=S'\times_SS_0$.
 We denote by $\mathcal{S}'^\rig$ the formal completion of $S'$ along $S'_0$.
 Then, all constructions above are compatible with the base change
 by $f$. More precisely, for a scheme $X$ of finite type over $S$, 
 we have $(X\times_SS')^{\wedge\rig}\cong \widehat{X}^\rig\times_{\mathcal{S}^\rig}\mathcal{S}'^\rig$
 and $(X\times_SS')^{\text{$S'$-ad}}\cong X^\ad\times_{\mathcal{S}^\rig}\mathcal{S}'^\rig$.
 Here $(-)^{\text{$S'$-ad}}$ denotes the functor $(-)^\ad$ for the base $(S',S'_0)$,
 namely, $(-)^\text{$S'$-ad}=(-)\times_{S'}\mathcal{S}'^\rig$.
\end{rem}

In the remaining part of this subsection,
assume that $S$ is the spectrum of a complete discrete valuation ring $R$ and
$S_0$ is the closed point of $S$. For a scheme of finite type $X$ over $S$, 
we have a natural morphism of locally and topologically ringed spaces 
$(t(\widehat{X}),\mathcal{O}_{\!t(\widehat{X})}^+)\to (\widehat{X},\mathcal{O}_{\widehat{X}})$
(\cf \cite[Proposition 4.1]{MR1306024}).
Note that the underlying continuous map $t(\widehat{X})\to X_0$ is different from the map
$t(\widehat{X})\to X$ considered above.
We denote the composite $\widehat{X}^\rig\hookrightarrow t(\widehat{X})\to X_0$
by $\spp_{\widehat{X}}$, or simply by $\spp$.

Let $Y$ be a closed subscheme of $X_0$ and $\mathcal{X}$ the formal completion of $X$ along $Y$.
Then we can consider the generic fiber $t(\mathcal{X})_\eta=S^0\times_St(\mathcal{X})$ of the adic space $t(\mathcal{X})$.
This is so-called the rigid generic fiber of $\mathcal{X}$ due to Raynaud and Berthelot, in the context of
adic spaces.
If $Y=X_0$, then $t(\mathcal{X})_\eta=\widehat{X}^\rig$.

\begin{lem}\label{lem:rig-gen-fiber}
 The natural morphism $t(\mathcal{X})_\eta\to \widehat{X}^\rig$ induced from
 $\mathcal{X}\to \widehat{X}$ is an open immersion.
 Its image coincides with $\spp^{-1}(Y)^\circ$, where $(-)^\circ$ denotes the interior in $\widehat{X}^\rig$.
\end{lem}

\begin{prf}
 See \cite[Lemma 3.13 i)]{MR1620118}.
\end{prf}

Let $X$ be an adic space locally of finite type over $\mathcal{S}^\rig=\Spa(F,R)$,
where $F$ denotes the fraction field of $R$.
For a point $x$ of $X$, we write $\kappa_x$ and $\kappa^+_x$ for the residue field and the valuation ring at $x$,
respectively.
We say that $x\in X$ is classical if $\kappa_x$ is a finite extension of $F$.
We denote the set of classical points of $X$ by $X(\cl)$.
Further, for a subset $Y$ of $X$, we put $Y(\cl)=X(\cl)\cap Y$.

\begin{lem}\label{lem:cl-pt-constructible}
 Let $X$ be an adic space locally of finite type over $\Spa(F,R)$.
 \begin{enumerate}
  \item For constructible subsets $L_1$, $L_2$ of $X$ (see \cite[1.1.13]{MR1734903}),
	we have $L_1\subset L_2$ if and only if $L_1(\cl)\subset L_2(\cl)$.
	In particular, $L_1=L_2$ if and only if $L_1(\cl)=L_2(\cl)$.
  \item For a constructible subset $L$, we write $L^-$ (resp.\ $L^\circ$) for the closure (resp.\ interior)
	of $L$ in $X$. Then we have $L(\cl)=L^-(\cl)=L^\circ(\cl)$.
 \end{enumerate}
\end{lem}

\begin{prf}
 For (i), it suffices to show that $L_1(\cl)\subset L_2(\cl)$ implies $L_1\subset L_2$.
 Put $L=L_1\setminus L_2=L_1\cap (X\setminus L_2)$. It is a constructible subset of $X$ satisfying
 $L(\cl)=\varnothing$.
 Let $U$ be an arbitrary affinoid open subset of $X$ .
 Then, we have $(U\cap L)(\cl)=\varnothing$.
 Therefore, \cite[Corollary 4.3]{MR1207303} tells us that $U\cap L=\varnothing$.
 Now we conclude that $L=\varnothing$, that is, $L_1\subset L_2$.

 For (ii), it suffices to prove that $L(\cl)=L^-(\cl)$.
 Take $x\in L^-\setminus L$ and an affinoid open neighborhood $U$ of $x$.
 Then $x$ lies in the closure of $U\cap L$ in $U$.
 Since $U\cap L$ is a constructible subset of the spectral space $U$, 
 by \cite[Corollary of Theorem 1]{MR0251026}, there exists $y\in U\cap L$ such that $x\in \{y\}^-$.
 Therefore, by \cite[Lemma 1.1.10 ii)]{MR1734903}, the valuation $v_x$ attached to $x$ is
 not rank $1$. In particular $x$ is not classical.
 Hence we have $L(\cl)=L^-(\cl)$, as desired.
\end{prf}

The following basic lemma is also used in Section \ref{sec:partition}.

\begin{lem}\label{lem:etale-constructible-image}
 Let $f\colon X\to Y$ be a quasi-compact quasi-separated \'etale morphism between adic spaces.
 \begin{enumerate}
  \item\label{item:constr-image} For a constructible subset $L$ of $X$, the image $f(L)$ is
       a constructible subset of $Y$.
  \item\label{item:loc-cl-image} For a locally closed subset $L$ of $X$
       satisfying $f^{-1}(f(L))=L$, the image $f(L)$ is a locally closed subset of $Y$.
 \end{enumerate}
\end{lem}

\begin{prf}
 The assertion \ref{item:constr-image} can be proved in the same way as
 \cite[(1) in the proof of Lemma 2.7.4]{MR1734903}.
 We recall the argument for reader's convenience.
 We may assume that $X$ and $Y$ are quasi-compact and quasi-separated.
 Fix $y\in f(L)$. Let $\Lambda$ denote the set of
 constructible subsets of $Y$ containing $y$. We have $\bigcap_{W\in \Lambda}W=\{y\}$, 
 as $Y$ is a spectral space. Since $f^{-1}(y)$ is a finite discrete subset of $X$,
 there exists a quasi-compact open subset $U$ of $X$ such that $U\cap f^{-1}(y)=L\cap f^{-1}(y)$.
 Then we have $U\cap \bigcap_{W\in\Lambda}f^{-1}(W)=L\cap \bigcap_{W\in\Lambda}f^{-1}(W)$.
 By the quasi-compactness of $X$ with respect to the constructible topology, 
 there exists $W\in \Lambda$ such that $U\cap f^{-1}(W)=L\cap f^{-1}(W)$.
 Put $V_y=U\cap f^{-1}(W)=L\cap f^{-1}(W)$, which is a constructible subset of $X$.
 Since $f$ is \'etale, $f(U)$ is a quasi-compact open subset of $Y$. 
 Therefore $f(V_y)=f(U)\cap W$ is a constructible subset of $Y$.
 
 Since $L\cap f^{-1}(y)\subset V_y\subset L$, we have $L=\bigcup_{y\in f(L)}V_y$.
 On the other hand, $L$ is quasi-compact under the constructible topology of $X$.
 Therefore, there exist finitely many points $y_1,\ldots,y_m\in f(L)$ such that
 $L=\bigcup_{i=1}^m V_{y_i}$. 
 Now we conclude that $f(L)=\bigcup_{i=1}^mf(V_{y_i})$ is a constructible subset of $Y$, as desired. 

 Next we consider \ref{item:loc-cl-image}. Since $L$ is locally closed, it can be written in the form $U\cap W$,
 where $U$ is an open subset of $X$ and $W$ is a closed subset of $X$.
 Note that $L^-\subset W$, thus $U\cap L^-=L$.
 For simplicity we write $L'=f(L)$. Since $f$ is an open map, we can check that $f^{-1}(L')^-=f^{-1}(L'^-)$.
 Therefore we obtain
 \[
 L=U\cap L^-=U\cap f^{-1}(L')^-=U\cap f^{-1}(L'^-)
 \]
 and $f(L)=f(U)\cap L'^-$. As $f$ is \'etale, $f(U)$ is open, hence $f(L)$ is locally closed.
\end{prf}

\subsection{Etale sheaves associated to semi-abelian schemes}\label{subsec:etale-sheaves-semi-abelian}
We continue to use the notation introduced in the beginning of the previous subsection.
Let $U$ be an open subscheme of $S^0=S\setminus S_0$ and $\ell$ a prime number invertible on $U$.
Fix an integer $m>0$.

Let $G$ be a semi-abelian scheme over $S$. Namely, $G$ is a separated smooth commutative
group scheme over $S$ such that each fiber $G_s$ of $G$ at $s\in S$ is an extension of
an abelian variety $A_s$ by a torus $T_s$. 
We denote the relative dimension of $G$ over $S$ by $d$.
Assume the following:

\begin{itemize}
 \item The rank of $T_s$ (called the toric rank of $G_s$) with $s\in S_0$ is a constant $r$.
 \item $G_U=G\times_SU$ is an abelian scheme.
\end{itemize}

Under the first condition, it is known that $G_0=G\times_SS_0$ is globally an extension
\[
 0\to T_0\to G_0\to A_0\to 0,
\]
where $T_0$ is a torus of rank $r$ over $S_0$ and $A_0$ is an abelian scheme over $S_0$
(\cite[Chapter I, Corollary 2.11]{MR1083353}).

Let us consider two group spaces $\widehat{G}^\rig[\ell^m]_{U^\ad}$ and $G^\ad[\ell^m]_{U^\ad}$ over $U^\ad$, where $(-)_{U^\ad}$ denotes the restriction to $U^\ad$. 

\begin{lem}\label{lem:ad-fin-etale}
 The adic space $G^\ad[\ell^m]_{U^\ad}$ is finite \'etale of degree $\ell^{2dm}$ over $U^\ad$.
\end{lem}

\begin{prf}
 By Lemma \ref{lem:fiber-product}, we have $G^\ad[\ell^m]_{U^\ad}=(G_U[\ell^m])\times_UU^\ad$.
 Since $G_U[\ell^m]$ is finite \'etale of degree $\ell^{2dm}$ over $U$, 
 $G^\ad[\ell^m]_{U^\ad}$ is finite \'etale of degree $\ell^{2dm}$ over $U^\ad$
 (see \cite[Corollary 1.7.3 i)]{MR1734903}).
\end{prf}

\begin{lem}\label{lem:rig-fin-etale}
 The adic space $\widehat{G}^\rig[\ell^m]_{U^\ad}$ is finite \'etale
 of degree $\ell^{(2d-r)m}$ over $U^\ad$.
\end{lem}

\begin{prf}
 We may assume that $S=\Spec R$ is affine. Let $I\subset R$ be the defining ideal of $S_0$.
 By replacing $R$ by its $I$-adic completion, we can reduce to the case where $R$ is $I$-adically complete.
 Put $S_i=\Spec R/I^{i+1}$ and $G_i=G\times_SS_i$.

 By \cite[Expos\'e IX, Th\'eor\`eme 3.6, Th\'eor\`eme 3.6 bis]{SGA3}, the exact sequence
 \[
  0\to T_0\to G_0\to A_0\to 0
 \]
 can be lifted canonically to an exact sequence
 \[
  0\to T_i\to G_i\to A_i\to 0
 \]
 over $S_i$, where $T_i$ is a torus over $S_i$ and $A_i$ is an abelian scheme over $S_i$
 (see \cite[\S 3.3.3]{Kai-Wen}). 
 Let $\widehat{T}=\varinjlim_i T_i$ and 
 $\widehat{A}=\varinjlim_i A_i$ be associated formal groups over $\mathcal{S}$.
 Then $\widehat{G}$ is an extension of $\widehat{A}$ by $\widehat{T}$.

 By taking $\ell^m$-torsion points, we get an exact sequence
 \[
  0\to \widehat{T}[\ell^m]\to \widehat{G}[\ell^m]\to \widehat{A}[\ell^m]\to 0
 \]
 of formal groups over $\mathcal{S}$. 
 Since $\widehat{G}^\rig[\ell^m]\cong (\widehat{G}[\ell^m])^\rig$, it suffices to see that 
 $(\widehat{T}[\ell^m])^\rig_{U^\ad}$ (resp.\ $(\widehat{A}[\ell^m])^\rig_{U^\ad}$) is finite \'etale of
 degree $\ell^{rm}$ (resp.\ $\ell^{2(d-r)m}$) over $U^\ad$.

 First we consider $(\widehat{T}[\ell^m])^\rig_{U^\ad}$. 
 Since $\widehat{T}[\ell^m]=\varinjlim_i (T_i[\ell^m])$, it is finite flat over $\mathcal{S}=\Spf R$.
 Therefore there exists a finite flat $R$-algebra $R'$ such that $\widehat{T}[\ell^m]=\Spf R'$.
 Moreover, a scheme $T'=\Spec R'$ is naturally equipped with a structure of a commutative group scheme
 over $S=\Spec R$. Since $T'$ is killed by $\ell^m$ and $p$ is invertible on $U$, $T'_U=T'\times_SU$
 is a finite \'etale group scheme over $U$. 
 By Lemma \ref{lem:rig-ad} (ii), we have
 $(\widehat{T}[\ell^m])^\rig=(T')^{\wedge\rig}=T'^\ad=T'\times_S\mathcal{S}^\rig$.
 Therefore $(\widehat{T}[\ell^m])^\rig_{U^\ad}=T'_U\times_UU^\ad$ is finite \'etale over $U^\ad$
 (see \cite[Corollary 1.7.3 i)]{MR1734903}). Its degree is clearly $\ell^{rm}$.

 The same argument also works for $(\widehat{A}[\ell^m])^\rig_{U^\ad}$.
\end{prf}

By Lemma \ref{lem:ad-fin-etale} and Lemma \ref{lem:rig-fin-etale}, we may regard
$\widehat{G}^\rig[\ell^m]_{U^\ad}$ and $G^\ad[\ell^m]_{U^\ad}$ as locally constant constructible sheaves
over $U^\ad$. 
Since we have a natural open immersion $\widehat{G}^\rig\hookrightarrow G^\ad$,
$\widehat{G}^\rig[\ell^m]_{U^\ad}$ is a subsheaf of $G^\ad[\ell^m]_{U^\ad}$.

\begin{rem}\label{rem:etale-sheaf-base-change}
 In the setting of Remark \ref{rem:rig-ad-base-change},
 the construction above is clearly compatible with the base change by $f\colon S'\to S$.
\end{rem}

In the remaining part of this subsection, we consider the case where $S=\Spec R$ is the spectrum of
a complete discrete valuation ring $R$,
$S_0$ is the closed point of $S$ and $U=S^0=S\setminus S_0$. 
Let $\overline{\eta}$ be a geometric point lying over the unique point of $U^{\ad}$.

As in the proof of Lemma \ref{lem:rig-fin-etale}, $\widehat{G}$ is an extension
\[
 0\to \widehat{T}\to \widehat{G}\to \widehat{A}\to 0
\]
of a formal group $\widehat{A}$ by $\widehat{T}$.
Therefore, we have $\Z/\ell^m\Z$-submodules
\[
 \widehat{T}^\rig[\ell^m]_{\overline{\eta}}\subset \widehat{G}^\rig[\ell^m]_{\overline{\eta}}\subset G^\ad[\ell^m]_{\overline{\eta}}.
\]
By taking inverse limit and tensoring with $\Q_\ell$, we have
\[
 T_\ell\widehat{T}^\rig_{\overline{\eta}}\subset T_\ell\widehat{G}^\rig_{\overline{\eta}}\subset T_\ell G^\ad_{\overline{\eta}},\quad V_\ell\widehat{T}^\rig_{\overline{\eta}}\subset V_\ell\widehat{G}^\rig_{\overline{\eta}}\subset V_\ell G^\ad_{\overline{\eta}},
\]
where we put $V_{\ell}(-)=T_{\ell}(-)\otimes_{\Z_{\ell}}\Q_{\ell}$.
By Lemma \ref{lem:rig-fin-etale} and its proof, we can deduce that 
$\dim_{\Q_\ell}V_\ell\widehat{T}^\rig_{\overline{\eta}}=r$ and
$\dim_{\Q_\ell}V_\ell\widehat{G}^\rig_{\overline{\eta}}=2d-r$.

\begin{prop}\label{prop:Afilt}
 Assume that the fraction field $F$ of $R$ is a finite extension of $\Q_p$. 
 For a filtration
 \[
 0\subset V_{\ell}\widehat{T}^\rig_{\overline{\eta}}\subset V_{\ell}\widehat{G}^\rig_{\overline{\eta}}
 \subset V_{\ell} G^\ad_{\overline{\eta}}=V_{\ell} G_{\overline{\eta}}, 
\]
 we have the following:
\begin{enumerate}
 \item\label{item:Afilt1} If $\ell \neq p$, the above filtration is the weight filtration 
     of $V_{\ell} G_{\overline{\eta}}$. 
 \item\label{item:Afilt2} If $\ell=p$, then the above filtration is a filtration 
     as semistable representations of $\Gal(\overline{F}/F)$. 
     Further, this filtration induces the weight filtration on $D_{\mathrm{st}}(V_pG_{\overline{\eta}})$. 
\end{enumerate}
\end{prop}

\begin{prf}
We give a proof of \ref{item:Afilt2}. 
We can show \ref{item:Afilt1} similarly. 
 Let $\lambda$ be a polarization of $G_U$. 
Then an alternating bilinear pairing
 \[
  \langle\ ,\ \rangle_\lambda\colon V_pG_{\overline{\eta}}\times V_pG_{\overline{\eta}}\to \Q_p(1)
 \]
 is induced. 
First, we will prove that 
 $(V_p\widehat{G}^\rig_{\overline{\eta}})^\perp=V_p\widehat{T}^{\rig}_{\overline{\eta}}$. 
 Since 
 \[
  \dim_{\Q_p}V_p\widehat{T}^\rig_{\overline{\eta}}+\dim_{\Q_p} V_p\widehat{G}^\rig_{\overline{\eta}}
 =r+(2d-r)=2d=\dim_{\Q_p} V_pG_{\overline{\eta}},
 \]
 it is sufficient to prove that
 $V_p\widehat{T}^\rig_{\overline{\eta}}\subset (V_p\widehat{G}^\rig_{\overline{\eta}})^\perp$.
 Namely, we should prove that the homomorphism 
 $V_p\widehat{T}^\rig_{\overline{\eta}}\otimes_{\Q_p}V_p\widehat{G}^\rig_{\overline{\eta}}\to \Q_p(1)$ induced by $\langle\ ,\ \rangle_\lambda$ is zero.
 
 Since $V_pG_{\overline{\eta}}$ is a semistable representation of
 $\Gal(\overline{F}/F)$, so are $V_p\widehat{T}^\rig_{\overline{\eta}}$ and $V_p\widehat{G}^\rig_{\overline{\eta}}$.
 We denote the residue field of $F$ by $\kappa_F$ and put $q=\#\kappa_F$.
 Consider the action of $\varphi^{[\kappa_F:\F_p]}$ on
 $D_{\mathrm{st}}(V_p\widehat{T}^\rig_{\overline{\eta}})$
 and $D_{\mathrm{st}}(V_p\widehat{A}^\rig_{\overline{\eta}})$.
 By \cite[Expos\'e X, Th\'eor\`eme 3.2]{SGA3}, $\widehat{T}$ can be algebraized into a torus $T$ over $S$.
 Then we have $V_p\widehat{T}^\rig_{\overline{\eta}}\cong V_pT_{\overline{\eta}}$.
 Therefore every eigenvalue of $\varphi^{[\kappa_F:\F_p]}$ on $D_{\mathrm{st}}(V_p\widehat{T}^\rig_{\overline{\eta}})$ is a Weil $q^{-2}$-number (for the definition of Weil numbers, see \cite[p.~471]{MR2276777}).
 Similarly, by \cite[Proposition 3.3.3.6, Remark 3.3.3.9]{Kai-Wen},
 $\widehat{A}$ can be algebraized into an abelian scheme $A$ over $S$, and
 we have $V_p\widehat{A}^\rig_{\overline{\eta}}\cong V_pA_{\overline{\eta}}$.
 By the Weil conjecture for the crystalline cohomology of abelian varieties,
 every eigenvalue of $\varphi^{[\kappa_F:\F_p]}$ on $D_{\mathrm{st}}(V_p\widehat{A}^\rig_{\overline{\eta}})$
 is a Weil $q^{-1}$-number.
 Therefore, every eigenvalue of $\varphi^{[\kappa_F:\F_p]}$ on 
 $D_{\mathrm{st}}(V_p\widehat{T}^\rig_{\overline{\eta}}\otimes_{\Q_p}V_p\widehat{G}^\rig_{\overline{\eta}})$
 is either a Weil $q^{-4}$-number or
 a Weil $q^{-3}$-number.
 On the other hand, every eigenvalue of $\varphi^{[\kappa_F:\F_p]}$ on $D_{\mathrm{st}}(\Q_p(1))$ is equal to
 $q^{-1}$, which is a Weil $q^{-2}$-number.
 Hence any $\varphi$-homomorphism 
 $D_{\mathrm{st}}(V_p\widehat{T}^\rig_{\overline{\eta}}\otimes_{\Q_p}V_p\widehat{G}^\rig_{\overline{\eta}})\to D_{\mathrm{st}}(\Q_p(1))$ is zero.
 Since the functor $D_{\mathrm{st}}$ is fully faithful, 
 any $\Gal(\overline{F}/F)$-equivariant homomorphism
 \[
  V_p\widehat{T}^\rig_{\overline{\eta}}\otimes_{\Q_p}V_p\widehat{G}^\rig_{\overline{\eta}}\to \Q_p(1)
 \]
 is zero. Hence, we have 
 $(V_p\widehat{G}^\rig_{\overline{\eta}})^\perp=V_p\widehat{T}^{\rig}_{\overline{\eta}}$. 
Then we have a perfect pairing 
\[
 V_p\widehat{T}^\rig_{\overline{\eta}} \times 
 (V_pG_{\overline{\eta}} / 
 V_p\widehat{G}^\rig_{\overline{\eta}}) 
 \to \Q_p(1). 
\]
The claim follows from the above arguments and this perfect pairing. 
\end{prf}

\begin{cor}\label{cor:semiab-monodromy-filt}
 The semistable representation $V_pG_{\overline{\eta}}$ of $\Gal(\overline{F}/F)$ has the monodromy filtration
 in the sense of Definition \ref{defn:has-monodromy-filt-rep-fin}.
\end{cor}

\begin{prf}
 It is well-known that in this case the weight filtration and the monodromy filtration on
 $D_{\mathrm{pst}}(V_pG_\eta)=D_{\mathrm{st}}(V_pG_\eta)$ coincide up to shift.
 Therefore the claim follows from Proposition \ref{prop:Afilt} \ref{item:Afilt2}.
\end{prf}

\begin{rem}\label{rem:rig-Raynaud}
 Actually, the extension $0\to \widehat{T}\to \widehat{G}\to \widehat{A}\to 0$ considered above
 can be algebraized; namely, there exists an exact sequence
 \[
 0\to T\to G^\natural\to A\to 0
\] 
 of commutative group schemes over $S$, where $T$ and $A$ are as in the proof of Proposition \ref{prop:Afilt},
 such that its formal completion along the special fiber is isomorphic to the extension above
 (see \cite[Proposition 3.3.3.6, Remark 3.3.3.9]{Kai-Wen}). Such an extension is called the Raynaud extension
 associated to $G$.

 Our construction above is related to the Raynaud extension in the following way.
 First, we have a natural isomorphism
 $\widehat{G}^\rig[\ell^m]_{\overline{\eta}}\xrightarrow{\cong}(G^\natural)^\ad[\ell^m]_{\overline{\eta}}$,
 which is induced from an open immersion
 $\widehat{G}^\rig\cong (\widehat{G^\natural})^\rig\hookrightarrow (G^\natural)^\ad$ 
 (see Lemma \ref{lem:rig-ad} (i)). Moreover, the image of $\widehat{G}^\rig[\ell^m]_{\overline{\eta}}\hookrightarrow G^\ad[\ell^m]_{\overline{\eta}}$ coincides with the image of the map
 $G^\natural[\ell^m]_{\overline{\eta}}\to G[\ell^m]_{\overline{\eta}}$ in \cite[Corollary 4.5.3.12]{Kai-Wen}.
\end{rem}

\section{Shimura varieties}\label{sec:Shvar}
\subsection{Notation on Shimura varieties}\label{subsec:NotaSh}
Let $(G,X)$ be a Shimura datum, and $E(G,X)$ the reflex field of $(G,X)$. 
We simply write $E$ for $E(G,X)$ if there is no risk of confusion. 
There is the canonical model over $E$ of the Shimura variety for $(G,X)$, 
which we denote by $\{\Sh_K(G,X)\}_{K \subset G(\A^{\infty})}$. 
Let $K\subset G(\A^{\infty})$ be a compact open subgroup, which is always supposed to be small enough
so that $\Sh_K(G,X)$ becomes a scheme.

\subsection{Siegel modular varieties}\label{sec:Siegel}
Let $(V,\langle\ ,\ \rangle)$ be a symplectic space of dimension $2n$ over $\Q$, 
and $L$ a self-dual $\Z$-lattice of $V$.
Let $(\GSp_{2n},X_{2n})$ be the Shimura datum associated to $(V,\langle\ ,\ \rangle)$. 
Then the Shimura variety for $(\GSp_{2n},X_{2n})$ is called the Siegel modular variety. 
In this case the reflex field $E(\GSp_{2n},X_{2n})$ equals $\Q$. 
We put 
\[
 K (N) =\Ker (\GSp_{2n} (\widehat{\Z}) 
 \to \GSp_{2n} (\widehat{\Z}/N \widehat{\Z})) 
\]
for $N\ge 1$, and 
\[
 K_{\ell,m}=\Ker(\GSp_{2n}(\mathbb{Z}_{\ell})\to\GSp_{2n}(\Z/{\ell}^m\Z))
\]
for a prime number $\ell$ and $m\ge 0$.

We recall a moduli interpretation of $\Sh_K (\GSp_{2n},X_{2n})$ using integral level structures. 
For simplicity, we assume that $K=K(N)$ with $N\ge 3$. 
We consider the functor from the category of $\mathbb{Q}$-schemes to the category of sets, 
that associates $S$ to the set of isomorphism classes of triples $(A,\lambda,\eta)$, where 
\begin{itemize}
 \item $A$ is an abelian scheme over $S$,
 \item $\lambda\colon A\to A^\vee$ is a principal polarization, and
 \item $\eta \colon L/NL \xrightarrow{\cong} A[N]$ is a symplectic similitude. 
\end{itemize}
This functor is represented by $\Sh_K (\GSp_{2n},X_{2n})$ (see \cite[4.16]{MR0498581}). 

There is another moduli interpretation using rational level structures. 
Let $S$ be a connected Noetherian scheme over $\Q$, and fix a geometric point $\overline{s}$ of $S$. 
We put 
\[
 T^{\infty}(-)=\prod_{\ell} T_{\ell}(-), \quad V^{\infty}(-)=T^{\infty}(-)\otimes_{\Z}\Q, 
\]
where $\ell$ in the product ranges over all prime numbers. 
Then, $S$-valued points of $\Sh_K (\GSp_{2n},X_{2n})$ 
correspond to the isogeny classes of triples $(A,\lambda,\eta K)$, where 
\begin{itemize}
 \item $A$ is an abelian scheme over $S$,
 \item $\lambda\colon A\to A^\vee$ is a $\mathbb{Q}$-polarization, and
 \item $\eta K$ is a $\pi_1(S,\overline{s})$-invariant $K$-orbit of symplectic similitudes 
       $V_{\mathbb{A}^{\infty}}\xrightarrow{\cong}V^{\infty}A$. 
\end{itemize}
Using this description, the Hecke action of $g \in \GSp_{2n} (\A^{\infty})$ 
can be described as 
\[
 \Sh_{K} \to \Sh_{g^{-1}Kg};\ 
 [(A,\lambda,\eta K)] \mapsto [(A,\lambda, (\eta \circ g)g^{-1}Kg)] . 
\]
See \cite[4.12]{MR0498581} for the relation between two moduli interpretations. 

Assume that $K=K_{p,0}K^p$ with a compact open subgroup $K^p$ of $\GSp_{2n}(\widehat{\Z}^p)$.
Then $\Sh_{K}(\GSp_{2n},X_{2n})$ has a natural integral model $\mathscr{S}_{K^p}$ over $\mathbb{Z}_p$ 
constructed as a moduli space of principally polarized abelian schemes with level structures 
(\cf \cite[Chapter 7, \S3]{MR1304906}). 
Let $\mathcal{A}$ denote the universal abelian scheme on $\mathscr{S}_{K^p}$. 

Thanks to a work of Faltings and Chai \cite{MR1083353}, 
we have a toroidal compactification $\mathscr{S}_{K^p}^{\mathrm{tor}}$ of $\mathscr{S}_{K^p}$
over $\mathbb{Z}_p$.
We have a semi-abelian scheme on $\mathscr{S}_{K^p}^{\mathrm{tor}}$ extending $\mathcal{A}$ on $\mathscr{S}_{K^p}$, 
for which we write the same symbol $\mathcal{A}$. 

\subsection{Shimura varieties of Hodge type}\label{subsec:Hodge-type}
In this subsection, we assume that $(G,X)$ is of Hodge type. 
We take an embedding $i\colon (G,X)\hookrightarrow (\GSp_{2n},X_{2n})$ of Shimura data.
For a compact open subgroup $\widetilde{K}$ of $\GSp_{2n}(\A^\infty)$ containing $K$,
we have a natural morphism $\Sh_K(G,X)\to \Sh_{\widetilde{K}}(\GSp_{2n},X_{2n})$
to the Siegel modular variety,
which is known to be a closed immersion if $\widetilde{K}$ is small enough.
We shall recall a moduli interpretation of $\C$-points of $\Sh_K(G,X)$.
Let $V$ be the standard representation of $\GSp_{2n}$. 
By \cite[I, Proposition 3.1]{MR654325}, there exists a finite collection of tensors 
$(s_\alpha)_{\alpha\in J'}$ with $s_\alpha\in V^{m_\alpha}\otimes V^{\vee m'_\alpha}$ such that
$G$ equals the pointwise stabilizer of $(s_\alpha)_{\alpha\in J'}$ in $\GSp_{2n}$.
We put $J=J'\amalg\{0\}$, $m_0=m'_0=1$ and let $s_0$ be the symplectic form
$\langle\ ,\ \rangle\in V\otimes V^\vee$ on $V$.

\begin{prop}\label{prop:Hodge-C-moduli}
 A $\C$-valued point of $\Sh_K(G,X)$ corresponds to the isogeny class of
 triples $(A,(t_\alpha)_{\alpha\in J},\eta K)$, where
\begin{itemize}
 \item $A$ is an abelian variety over $\C$,
 \item $(t_\alpha)_{\alpha\in J}$ with
       $t_\alpha\in H_1(A,\Q)^{m_\alpha}\otimes H_1(A,\Q)^{\vee m'_\alpha}$
       is a finite collection of Hodge cycles on $A$ (see \cite[V, \S 2]{MR654325})
       such that $\pm t_0$ is a polarization of the rational Hodge structure $H_1(A,\Q)$,
 \item $\eta K$ is a $K$-orbit of $\A^\infty$-linear isomorphisms $V_{\A^\infty}\xrightarrow{\cong}V^\infty A$
       which send $s_0$ to a $(\A^{\infty})^\times$-multiple of $t_0$ and $s_\alpha$ with $\alpha\in J'$ to
       $t_\alpha$,
\end{itemize}
 satisfying the following condition $(*)$:
 \begin{itemize}
  \item[$(*)$] there exists an isomorphism $\eta_\Q\colon V\xrightarrow{\cong}H_1(A,\Q)$ such that
	       $\eta_\Q^{-1}$ sends $t_0$ to a $\Q^\times$-multiple of $s_0$, $t_\alpha$ with $\alpha\in J'$
	       to $s_\alpha$, and the Hodge structure on $H_1(A,\Q)$ to a Hodge structure on $V$
	       induced by an element of $X$ and the embedding $i\colon G\hookrightarrow \GSp_{2n}$.
 \end{itemize}
\end{prop}
For a proof, see \cite[Theorem 7.4]{MR2192012}.

\begin{lem}\label{lem:Hodge-level-str}
Let $F$ be a $p$-adic field containing the reflex field $E$, and $x$ an $F$-valued point of $\Sh_K(G,X)$.
Choose an algebraic closure $\overline{F}$ of $F$ and denote by $\overline{x}$ the corresponding
geometric point over $x$.

We take an isomorphism $\iota\colon \overline{F}\xrightarrow{\cong}\C$ over $E$,
and write $\iota\overline{x}$ for the $\C$-valued point of $\Sh_K(G,X)$ determined by $\overline{x}$
and $\iota$.
Let $(A,(t_\alpha)_{\alpha\in J},\eta K)$ be a triple in the isogeny class corresponding to $\iota\overline{x}$
such that $A=\mathcal{A}_{\overline{x}}\otimes_{\overline{F},\iota}\C$.
Here $\mathcal{A}_{\overline{x}}$ is the abelian variety corresponding to
the image of $\overline{x}$ in the Siegel modular variety.
Let us choose a representative $\eta$ of $\eta K$. 
Under $\iota$, it corresponds to a trivialization of the $K$-torsor $\pi_K^{-1}(\overline{x})$
on $\overline{x}$, where $\pi_K$ denotes the natural map $\varprojlim_{K'\subset K}\Sh_{K'}(G,X)\to \Sh_K(G,X)$.

 For a prime number $\ell$, let $\mathcal{L}_{\mathrm{Std}\circ i,\ell}$ be the $\Q_\ell$-sheaf on
 $\Sh_K(G,X)$ corresponding to the representation $\mathrm{Std}\circ i$ of $G$ on $V$.
 For the stalk $\mathcal{L}_{\mathrm{Std}\circ i,\ell,\overline{x}}$, the following hold:
 \begin{enumerate}
  \item\label{item:stalk-isom} We have a canonical $\Gal(\overline{F}/F)$-equivariant isomorphism 
	$\mathcal{L}_{\mathrm{Std}\circ i,\ell,\overline{x}}\cong V_\ell \mathcal{A}_{\overline{x}}$.
  \item\label{item:stalk-trivialization}
        Each trivialization of the $K$-torsor $\pi_K^{-1}(\overline{x})$ determines an isomorphism
	$V_{\Q_\ell}\xrightarrow{\cong}\mathcal{L}_{\mathrm{Std}\circ i,\ell,\overline{x}}$.
	The isomorphism given by the trivialization corresponding to the chosen representative $\eta$
	equals the composite of
	$V_{\Q_\ell}\xrightarrow[\cong]{\eta_\ell}V_\ell A\xrightarrow[\cong]{\iota^{-1}}V_\ell\mathcal{A}_{\overline{x}}\stackrel{\mathrm{(i)}}{\cong}\mathcal{L}_{\mathrm{Std}\circ i,\ell,\overline{x}}$,
	where $\eta_\ell$ denotes the $\ell$-part of $\eta$.
 \end{enumerate}
\end{lem}

\begin{prf}
 The first assertion is essentially a statement for the Siegel case, which is well-known.
 The second can be checked directly by working over $\C$.
\end{prf}

\subsection{Shimura varieties of preabelian type}

\begin{defn}\label{defn:preabSh}
A Shimura datum $(G,X)$ is said to be of 
preabelian type if 
there exists a Shimura datum $(G',X')$ of 
Hodge type such that 
$(G^{\ad},X^{\ad}) \cong (G'^{\ad},X'^{\ad})$. 
If a Shimura data is of 
preabelian type, 
the associated Shimura variety is said to be 
of preabelian type (\cf \cite[p.~402]{MR1796512}). 
\end{defn}

\begin{lem}\label{lem:preHod}
Assume that $(G,X)$ is of preabelian type. 
We take a Shimura datum $(G',X')$ of Hodge type such that $(G^{\ad},X^{\ad})\cong (G'^{\ad},X'^{\ad})$. 
Let $K''$ be a compact open subgroup of $G^\ad(\A^\infty)$ which contains the image of $K$
under the map $G(\A^\infty)\to G^{\ad}(\A^\infty)$.
We regard it as a compact open subgroup of $G'^{\ad}(\A^\infty)$ by the isomorphism $G^{\ad}\cong G'^{\ad}$.
Then there exist a compact open subgroup
$K'\subset G'(\mathbb{A}^{\infty})$ and $g_1,\ldots,g_m \in G'^{\ad}(\mathbb{A}^{\infty})$ 
such that the following hold:
\begin{enumerate}
\item The morphism $(G',X') \to (G'^{\ad},X'^{\ad})$ and the conjugation by $g_i$ induces the morphism 
      \[
      f_i\colon \Sh_{K'}(G',X')\to \Sh_{g_i^{-1}K''g_i}(G'^{\ad},X'^{\ad})\to \Sh_{K''}(G'^{\ad},X'^{\ad})
      \] 
      for each $i$. 
\item The morphism 
      \[
      \coprod_{1\le i\le m}f_i\colon \coprod_{1\le i\le m}\Sh_{K'}(G',X')\to\Sh_{K''}(G'^{\ad},X'^{\ad})
      \]
      is surjective. 
\end{enumerate}
\end{lem}
\begin{prf}
This follows from the definition of Shimura varieties of preabelian type and the fact that 
Hecke action is transitive on the connected components of a Shimura variety. 
\end{prf}

\section{Partition of Shimura varieties}\label{sec:partition}
\subsection{Partition of classical points}\label{subsec:Partcl}
We fix a prime number $p$ and a finite place $v$ of $E$ above $p$. 
We write $\mathcal{O}_v$ for the ring of integers of $E_v$. 

Throughout the paper, we assume that a compact open subgroup $K$ of $G(\A^\infty)$ is small enough
so that the following conditions are satisfied:
\begin{itemize}
 \item The morphism $\pi_K\colon \varprojlim_{K'\subset K}\Sh_{K'}(G,X)\to \Sh_K(G,X)$ is a torsor
       under the quotient $K_{\Sh}$ of $K$ by a closed subgroup of $K\cap Z(G)(\A^\infty)$
       (\cf \cite[Theorem 5.28]{MR2192012}).
 \item If the Shimura datum $(G,X)$ satisfies the condition SV5 in \cite[p.~311]{MR2192012},
       then $K_{\Sh}$ equals $K$.
       Note that a Shimura datum of Hodge type satisfies SV5.
\end{itemize}

\begin{defn}
 Let $x$ be a classical point of $\Sh_K(G,X)_{E_v}^{\ad}$, and $\overline{\kappa}_x$ an algebraic closure
 of $\kappa_x$. We write $\overline{x}$ for the geometric point corresponding to $\overline{\kappa}_x$.
 
 By taking the pull-back of $\pi_K\colon\varprojlim_{K'\subset K}\Sh_{K'}(G,X)\to \Sh_K(G,X)$,
 we obtain a $K_{\Sh}$-torsor $\pi_K^{-1}(x)$ on $x$.
 This torsor and its trivialization $\eta$ over $\overline{x}$ give rise to a continuous homomorphism 
 $\phi_{x,\eta}\colon \Gal(\overline{\kappa}_x/\kappa_x)\to K_{\mathrm{Sh}}$.
 If we change $\eta$, the homomorphism $\phi_{x,\eta}$ changes by a $K_{\Sh}$-conjugation.
  \begin{enumerate}
  \item We write $\phi^\ad_{x,\eta}$ for the composite 
	$\Gal(\overline{\kappa}_x/\kappa_x)\xrightarrow{\phi_{x,\eta}} K_{\mathrm{Sh}}\to G^\ad(\A^\infty)$. 
	If we change $\eta$, the homomorphism $\phi^\ad_{x,\eta}$ changes by a $K^\ad$-conjugation,
	where $K^\ad$ denotes the image of $K$ in $G^\ad(\A^\infty)$.
	When we are only interested in the $K^\ad$-conjugacy class of $\phi^\ad_{x,\eta}$,
	we often drop the subscript $\eta$ and simply write $\phi^\ad_x$ for $\phi^\ad_{x,\eta}$.
	
	For a prime number $\ell$, we denote by $\phi^\ad_{x,\eta,\ell}$ the composite of $\phi^\ad_{x,\eta}$
	and the projection $G^\ad(\A^\infty)\to G^\ad(\Q_\ell)$.
  \item Assume that $(G,X)$ satisfies the condition SV5. 
	Then we write $\phi_{x,\eta}$ for the composite 
	$\Gal(\overline{\kappa}_x/\kappa_x)\xrightarrow{\phi_{x,\eta}}K\to G(\A^\infty)$.
	As in (i), we often write $\phi_x$ for $\phi_{x,\eta}$,
	which is well-defined up to $K$-conjugacy.

	For a prime number $\ell$, we define $\phi_{x,\eta,\ell}$ similarly.
 \end{enumerate}
\end{defn}

 \begin{rem}\label{rem:autom-sheaves}
  The homomorphism $\phi^\ad_{x,\eta,\ell}$ is related to
  $\ell$-adic automorphic \'etale sheaves on $\Sh_K(G,X)$	as follows.
  Let $(\xi,V_\xi)$ be a finite-dimensional algebraic representation of $G$
  over $\overline{\Q}_\ell$ such that $\Ker\xi$ contains $\Ker(K\to K_{\Sh})$.
  Then, we have an associated smooth $\overline{\Q}_\ell$-sheaf $\mathcal{L}_\xi$ on $\Sh_K(G,X)$
  (\cf \cite[Remark III.6.1]{MR1044823}). 
  As in Lemma \ref{lem:Hodge-level-str} \ref{item:stalk-trivialization},
  the trivialization $\eta$ of $\pi^{-1}(\overline{x})$
  induces an isomorphism $\mathcal{L}_{\xi,\overline{x}}\cong V_{\xi}$.
  Hence we obtain an $\ell$-adic Galois representation 
  $\Gal(\overline{\kappa}_x/\kappa_x)\to \GL(\mathcal{L}_{\xi,\overline{x}})\cong \GL(V_\xi)$.
  \begin{enumerate}
   \item\label{item:sheaves-all} If $\xi$ factors through $G^\ad$, it is equal to the composite
	$\Gal(\overline{\kappa}_x/\kappa_x)\xrightarrow{\phi^\ad_{x,\eta,\ell}}G^\ad(\Q_\ell)\xrightarrow{\xi}\GL(V_\xi)$.
   \item\label{item:sheaves-SV5} If $(G,X)$ satisfies SV5 (hence any $\xi$ is allowable), it is equal to the composite
	$\Gal(\overline{\kappa}_x/\kappa_x)\xrightarrow{\phi_{x,\eta,\ell}}G(\Q_\ell)\xrightarrow{\xi}\GL(V_\xi)$.
  \end{enumerate}
 \end{rem}

The following proposition can be checked easily.

\begin{prop}\label{prop:phi^ad-functoriality}
 Let $(G,X)\to (G',X')$ be a morphism of Shimura data such that $Z(G)$ is mapped into $Z(G')$.
 Let $K\subset G(\A^\infty)$ and $K'\subset G'(\A^\infty)$ be compact open subgroups
 such that $K$ is mapped into $K'$.
 For $x\in \Sh_K(G,X)^\ad_{E_v}(\cl)$, we write $x'$ for the image of $x$
 under the induced morphism $\Sh_K(G,X)\to \Sh_{K'}(G',X')$.
 Then the diagram
 \[
 \xymatrix{%
 \Gal(\overline{\kappa}_x/\kappa_x)\ar[r]^-{\phi^\ad_x}\ar[d]&G^\ad(\A^\infty)\ar[d]\\
 \Gal(\overline{\kappa}_{x'}/\kappa_{x'})\ar[r]^-{\phi^\ad_{x'}}&G'^\ad(\A^\infty)
 }
\]
 is commutative up to $K'^\ad$-conjugacy, where $K'^\ad$ denotes the image of $K'$ in
 $G'^{\mathrm{ad}}(\A^\infty)$.
\end{prop}

\begin{prop}\label{prop:phi_x-property}
 Assume that $(G,X)$ is of preabelian type. 
 \begin{enumerate}
  \item\label{phi_x-preab}
       For $x\in \Sh_K(G,X)^\ad_{E_v}(\cl)$, $\phi^\ad_{x,p}$ is de Rham and has the monodromy filtration.
  \item\label{phi_x-Hodge}
       Assume that $(G,X)$ is of Hodge type. 
       Then, for $x\in \Sh_K(G,X)^\ad_{E_v}(\cl)$, $\phi_{x,p}$ is de Rham and has the monodromy filtration.
 \end{enumerate}
\end{prop}

\begin{prf}
 By Lemma \ref{lem:monodromy-filt-stab-p} \ref{item:P-ad-p}, Lemma \ref{lem:monod-restr-p},
 Lemma \ref{lem:preHod} and Proposition \ref{prop:phi^ad-functoriality},
 the assertion \ref{phi_x-preab} is reduced to \ref{phi_x-Hodge}. We prove (ii).
 Take an embedding $i\colon (G,X)\hookrightarrow (\GSp_{2n},X_{2n})$ into a Siegel Shimura datum
 and a compact open subgroup $\widetilde{K}=\widetilde{K}_p\widetilde{K}^p$ of $\GSp_{2n}(\A^\infty)$
 containing $K$.
 By shrinking $K$, we may assume that $\widetilde{K}_p\subset\widetilde{K}_{p,0}=\GSp_{2n}(\Z_p)$
 and $\widetilde{K}^p$ is small enough.
 Then, the morphism 
 \[
  \Spec \kappa_x\to \Sh_K(G,X)_{E_v}\to \Sh_{\widetilde{K}_{p,0}\widetilde{K}^p}(\GSp_{2n},X_{2n})_{\Q_p}=\mathscr{S}_{\widetilde{K}^p,\Q_p}
 \]
 uniquely extends to $\Spec \kappa^+_x\to \mathscr{S}_{\widetilde{K}^p}^{\mathrm{tor}}$.
 Let $\mathcal{A}_{\kappa^+_x}$ denote the pull-back of the universal semi-abelian scheme $\mathcal{A}$
 by this morphism. It extends the abelian variety $\mathcal{A}_x$ over $\kappa_x$.
 Therefore, the representation $V_p\mathcal{A}_{\overline{x}}$ of $\Gal(\overline{\kappa}_x/\kappa_x)$ is semistable
 and has the monodromy filtration by Corollary \ref{cor:semiab-monodromy-filt}.

 Let $\mathrm{Std}\colon \GSp_{2n}\to\GL(V)$ denote the standard representation of $\GSp_{2n}$.
 By Lemma \ref{lem:Hodge-level-str} \ref{item:stalk-isom},
 we have a $\Gal(\overline{\kappa}_x/\kappa_x)$-equivariant isomorphism
 $\mathcal{L}_{\mathrm{Std}\circ i,\overline{x},p}\cong V_p\mathcal{A}_{\overline{x}}$.
 We fix a trivialization $\eta$ of the $K$-torsor $\pi_K^{-1}(x)$ over $\overline{x}$.
 By the isomorphism $V_{\Q_p}\cong\mathcal{L}_{\mathrm{Std}\circ i,\overline{x},p}$ induced from $\eta$,
 we regard $V_{\Q_p}$ as a representation of $\Gal(\overline{\kappa}_x/\kappa_x)$.
 As in Remark \ref{rem:autom-sheaves} \ref{item:sheaves-SV5},
 it is isomorphic to $\mathrm{Std}\circ i\circ \phi_{x,\eta,p}$.
 Summing up, we obtain a $\Gal(\overline{\kappa}_x/\kappa_x)$-equivariant isomorphism
 $\mathrm{Std}\circ i\circ \phi_{x,\eta,p}\cong V_p\mathcal{A}_{\overline{x}}$.
 Therefore, we conclude that $\mathrm{Std}\circ i\circ \phi_{x,\eta,p}$ is semistable (hence de Rham)
 and has the monodromy filtration.
 Since $\mathrm{Std}\circ i$ is a faithful representation of $G$, $\phi_{x,\eta,p}$ is de Rham and
 has the monodromy filtration by Lemma \ref{lem:monodromy-filt-p}.
 This completes the proof.
\end{prf}

\begin{rem}
 Assume that $(G,X)$ satisfies the condition SV6 in \cite[p.~312]{MR2192012}.
 Recently, Liu and Zhu announced a result that the $p$-adic sheaf $\mathcal{L}_{\xi,x}$ is de Rham
 for any finite-dimensional algebraic representation $\xi$ of $G^c$ over $\overline{\Q}_p$,
 where $G^c$ is the quotient of $G$ defined in \cite[p.~347]{MR1044823}
 (\cf \cite[Theorem 1.2]{2016arXiv160206282L}). 
 This implies that $\phi^{\mathrm{ad}}_{x,p}$ is de Rham.
 We do not use this remark later.
\end{rem}

In the sequel, we assume that $(G,X)$ is of preabelian type.
Take a finite non-empty set of primes $\Box$ such that $K=K_\Box K^\Box$,
where $K_\Box$ is a compact open subgroup of $G(\A_\Box)$ and $K^\Box$ is a hyperspecial compact open subgroup
of $G(\A^{\infty,\Box})$.
We write $\mathcal{P}_{G,\Box}(K_\Box)$
for the set of $K_\Box$-conjugacy classes of $\A_\Box$-parabolic subgroups of $G$.

Let $\eta$ be a trivialization of $\pi^{-1}_K(x)$ over $\overline{x}$.
By Proposition \ref{prop:phi_x-property} and the results in Section \ref{sec:Gal-rep},
we can attach to $\phi^\ad_{x,\eta,\ell}$ the $\Q_\ell$-parabolic subgroup $P_{\phi^\ad_{x,\eta,\ell}}$
of $G^\ad$ for each $\ell\in \Box$
and $x\in \Sh_K(G,X)^\ad_{E_v}(\cl)$. By taking the product with respect to $\ell$, 
we obtain an $\A_\Box$-parabolic subgroup of $G^\ad$.
It is easy to observe that the $K^\ad_\Box$-conjugacy class 
$[\prod_{\ell\in \Box}P_{\phi^\ad_{x,\eta,\ell}}]\in\mathcal{P}_{G^\ad,\Box}(K^\ad_\Box)$ is
independent of the choice of $\eta$.
Note that the natural map $\mathcal{P}_{G,\Box}(K_\Box)\to \mathcal{P}_{G^\ad,\Box}(K^\ad_\Box)$;
$[P]\mapsto [P^\ad]$ is bijective.

\begin{defn}
 Let $[P_{x,\Box}]\in \mathcal{P}_{G,\Box}(K_\Box)$ be the $K_\Box$-conjugacy class
 that is mapped to $[\prod_{\ell\in \Box}P_{\phi^\ad_{x,\eta,\ell}}]$ under the bijection
 $\mathcal{P}_{G,\Box}(K_\Box)\to \mathcal{P}_{G^\ad,\Box}(K^\ad_\Box)$.
\end{defn}

\begin{rem}\label{rem:P_x-SV5}
 If the Shimura datum $(G,X)$ satisfies the condition SV5, we can define $[P_{x,\Box}]$ directly by using
 $\phi_{x,\eta}$. These two ways give the same result by Lemma \ref{lem:monodromy-filt-stab-l} \ref{item:P-ad-l}
 and Lemma \ref{lem:monodromy-filt-stab-p} \ref{item:P-ad-p}.
\end{rem}

By the proof of Proposition \ref{prop:phi_x-property}, we obtain the following description of $[P_{x,\Box}]$
in the Hodge type case.

\begin{cor}\label{cor:Hodge-parabolic}
 Let $(G,X)$ be a Shimura datum of Hodge type with an embedding
 $(G,X)\hookrightarrow (\GSp_{2n},X_{2n})$ into a Siegel Shimura datum.
 Assume that $K_p\subset \GSp_{2n}(\Z_p)$ and $K^p$ is small enough. 
 For $x\in \Sh_K(G,X)^\ad_{E_v}(\cl)$, fix an isomorphism $\iota\colon \overline{\kappa}_x\xrightarrow{\cong}\C$
 and let $(A,(t_\alpha),\eta K)$ be a triple in the isogeny class corresponding to the $\C$-point $\iota\overline{x}$
 of $\Sh_K(G,X)$ such that $A=\mathcal{A}_{\overline{x}}\otimes_{\overline{\kappa}_x,\iota}\C$. 
 \begin{enumerate}
  \item\label{item:semiab-red}
       The abelian variety $\mathcal{A}_x$ over $\kappa_x$ extends to a semi-abelian scheme 
       $\mathcal{A}_{\kappa_x^+}$ over $\kappa_x^+$.
       For $\ell\in \Box$, the monodromy filtration $M_\bullet V_\ell\mathcal{A}_{\overline{x}}$
       on $V_\ell\mathcal{A}_{\overline{x}}$ is a shift of the filtration in Proposition \ref{prop:Afilt}.
  \item\label{item:monodromy-parabolic}
       Fix an arbitrary representative $\eta\colon V_{\A^\infty}\xrightarrow{\cong}V^\infty A$
       of the $K$-orbit $\eta K$.
       For $\ell\in \Box$, consider the filtration 
       $(\eta_\ell^{-1}\circ\iota)(M_\bullet V_\ell\mathcal{A}_{\overline{x}})$ on $V_{\Q_\ell}$
       obtained as the inverse image of $M_\bullet V_\ell \mathcal{A}_{\overline{x}}$
       under $V_{\Q_\ell}\xrightarrow[\cong]{\eta_\ell}V_\ell A\xrightarrow[\cong]{\iota^{-1}}V_\ell\mathcal{A}_{\overline{x}}$.
       Then, this filtration is $G_{\Q_\ell}$-split in the sense of \cite[(1.1.2)]{MR2669706}.
       Moreover, if we write $P_{x,\eta,\ell}$ for the stabilizer of this filtration,
       the $K_\Box$-conjugacy class $[\prod_{\ell\in \Box}P_{x,\eta,\ell}]$ equals $[P_{x,\Box}]$.
  \item\label{item:Siegel-toric}
       If $(G,X)=(\GSp_{2n},X_{2n})$, then $P_{x,\eta,\ell}$ in (ii) is the stabilizer of
       a totally isotropic subspace of $V_{\Q_\ell}$ whose dimension
       equals the toric rank of the special fiber of $\mathcal{A}_{\kappa_x^+}$.
 \end{enumerate}
\end{cor}

\begin{prf}
 The assertion \ref{item:semiab-red} follows from the proofs of Proposition \ref{prop:phi_x-property}
 and Corollary \ref{cor:semiab-monodromy-filt}.
 
 We prove \ref{item:monodromy-parabolic}.
 The choice of $\eta$ gives a trivialization of the $K$-torsor $\pi_K^{-1}(x)$ over $\overline{x}$,
 which is denoted by the same symbol $\eta$.
 By the argument in the proof of Proposition \ref{prop:phi_x-property},
 we have $\Gal(\overline{\kappa}_x/\kappa_x)$-equivariant isomorphisms
 \[
  \mathrm{Std}\circ i\circ\phi_{x,\eta,\ell}=V_{\Q_\ell}\stackrel{\eta}{\cong}\mathcal{L}_{\mathrm{Std}\circ i,\overline{x},\ell}\cong V_\ell\mathcal{A}_{\overline{x}}.
 \]
 By Lemma \ref{lem:Hodge-level-str} \ref{item:stalk-trivialization}, their composite is equal to
 \[
  \mathrm{Std}\circ i\circ\phi_{x,\eta,\ell}=V_{\Q_\ell}\xrightarrow[\cong]{\eta_\ell} V_\ell A\xrightarrow[\cong]{\iota^{-1}} V_\ell\mathcal{A}_{\overline{x}}.
 \]
 Hence the filtration $(\eta_\ell^{-1}\circ\iota)(M_\bullet V_\ell\mathcal{A}_{\overline{x}})$
 equals the monodromy filtration $M_\bullet V_{\Q_\ell}$ on $V_{\Q_\ell}$
 with respect to the action of $\Gal(\overline{\kappa}_x/\kappa_x)$
 by $\mathrm{Std}\circ i\circ\phi_{x,\eta,\ell}$.
 Since the monodromy filtration of $\mathrm{Std}\circ i\circ \phi_{x,\eta,\ell}$
 extends to a filtration on the Tannakian category $\mathbf{Rep}_{\Q_\ell}(G_{\Q_\ell})$,
 we conclude that $M_\bullet V_{\Q_\ell}$ is $G_{\Q_\ell}$-split
 by \cite[Lemma 1.1.3]{MR2669706}. Further, by Lemma \ref{lem:monodromy-filt-stab-l} \ref{item:P-faithful-l}
 and Lemma \ref{lem:monodromy-filt-stab-p} \ref{item:P-faithful-p},
 we have $P_{x,\eta,\ell}=P_{\phi_{x,\eta,\ell}}$.
 Therefore we have $[P_{x,\Box}]=[\prod_{x\in \Box}P_{x,\eta,\ell}]$ by Remark \ref{rem:P_x-SV5}.

 The claim \ref{item:Siegel-toric} follows from Corollary \ref{cor:semiab-monodromy-filt}
 and the equality $(V_p\widehat{G}^\rig_{\overline{\eta}})^\perp=V_p\widehat{T}^{\rig}_{\overline{\eta}}$
 (and its $\ell$-adic version) in the proof of Proposition \ref{prop:Afilt}.
\end{prf}

Next we will show that the $\ell$-part of $P_{x,\Box}$ is independent of $\ell\in \Box$ in some sense.
To state the result, we need some preparation.

\begin{defn}\label{defn:adm-parab}(\cf \cite[4.5 Definition]{MR1128753}) 
 Let $G^{\ad} =G_1 \times \cdots \times G_r$ be a decomposition into $\Q$-simple factors. 
 We say that a parabolic subgroup $P$ of $G$ is an admissible $\mathbb{Q}$-parabolic subgroup 
 if there exists a parabolic subgroup $P_i$ of $G_i$ for each $i$ such that 
 $P$ is the inverse image of $P_1 \times \dots \times P_r$ and $P_i$ is either equal to $G_i$ or 
 a maximal $\mathbb{Q}$-parabolic subgroup of $G_i$ for each $i$. 
 We write $\mathcal{P}_{G,\Q}$ for the set of $G(\Q)$-conjugacy classes of admissible $\Q$-parabolic
 subgroups of $G$.

 An admissible $\A^\infty$-parabolic subgroup means a parabolic subgroup of $G_{\A^\infty}$
 which is $G(\mathbb{A}^{\infty})$-conjugate
 to an admissible $\mathbb{Q}$-parabolic subgroup of $G$.
 Let $\mathcal{P}_G(K)$ denote the set of $K$-conjugacy classes of
 admissible $\mathbb{A}^{\infty}$-parabolic subgroups of $G$.
 Further, we write $\mathcal{P}_{G,\A^\infty}$ for the set of $G(\A^\infty)$-conjugacy classes of
 admissible $\A^\infty$-parabolic subgroups of $G$.
 We have a natural map $\mathcal{P}_G(K)\to \mathcal{P}_{G,\A^\infty}$.
\end{defn}

\begin{lem}\label{lem:P_G(K)}
 \begin{enumerate}
  \item\label{item:P_G(K)-local} 
       The natural map $\mathcal{P}_G(K)\to \mathcal{P}_{G,\Box}(K_\Box)$ is injective.
  \item\label{item:P_G(K)-finite}
       The set $\mathcal{P}_G(K)$ is finite.
  \item\label{item:P_G(K)-adjoint}
       We take a hyperspecial compact open subgroup $K''^S$ of $G^\ad(\A^{\infty,S})$ containing
       the image of $K^S$, and put $K''=K_S^\ad K''^S$, which is a compact open subgroup of $G^\ad(\A^\infty)$.
       Then, the natural map $\mathcal{P}_G(K)\to \mathcal{P}_{G^\ad}(K'')$ is bijective.
 \end{enumerate}
\end{lem}

\begin{prf}
 Fix a minimal parabolic subgroup $P_0$ of $G$.
 For an admissible $\Q$-parabolic subgroup $P$ containing $P_0$, we write $\mathcal{P}_G(K)_P$
 for the subset of $\mathcal{P}_G(K)$ consisting of $K$-conjugacy classes
 which are $G(\mathbb{A}^{\infty})$-conjugate to $P$.
 Then, we have a bijection $K\backslash G(\A^\infty)/P(\A^\infty)\xrightarrow{\cong} \mathcal{P}_G(K)_P$
 given by $KgP(\A^\infty)\mapsto gPg^{-1}$.

 Let us prove \ref{item:P_G(K)-local}. For an admissible $\Q$-parabolic subgroup $P$ containing $P_0$,
 we have
 \[
  K\backslash G(\A^\infty)/P(\A^\infty)=K_\Box\backslash G(\A_\Box)/P(\A_\Box)\times K^\Box\backslash G(\A^\Box)/P(\A^\Box)=K_\Box\backslash G(\A_\Box)/P(\A_\Box).
 \]
 by the Iwasawa decomposition $K^\Box P_0(\A^\Box)=G(\A^\Box)$.
 This implies that the composite $\mathcal{P}_G(K)_P\hookrightarrow \mathcal{P}_G(K)\to \mathcal{P}_{G,\Box}(K_\Box)$
 is injective.
 It suffices to show that the images of $\mathcal{P}_G(K)_P$ and $\mathcal{P}_G(K)_{P'}$ in 
 $\mathcal{P}_{G,\Box}(K_\Box)$ are disjoint,
 where $P$ and $P'$ are distinct admissible $\Q$-parabolic subgroups containing $P_0$.
 If the images of $\mathcal{P}_G(K)_P$ and $\mathcal{P}_G(K)_{P'}$ intersect, then
 $P_{\Q_\ell}$ and $P'_{\Q_\ell}$ are $G(\Q_\ell)$-conjugate for each $\ell\in \Box$.
 By \cite[Th\'eor\`eme 4.13]{MR0207712}, this means that $P$ and $P'$ are $G(\Q)$-conjugate.
 Since they contain $P_0$, they are equal.
 Note that in particular we have $\mathcal{P}_G(K)_P\cap \mathcal{P}_G(K)_{P'}=\varnothing$.
 Hence $\mathcal{P}_G(K)$ equals $\coprod_{P\supset P_0}\mathcal{P}_G(K)_P$, 
 where $P$ runs through admissible $\Q$-parabolic subgroups of $G$ containing $P_0$.
 
 Next we prove \ref{item:P_G(K)-finite}. It suffices to show that $\mathcal{P}_G(K)_P$ is a finite set
 for each admissible $\Q$-parabolic subgroup $P$ of $G$ containing $P_0$.
 Since
 \[
 \mathcal{P}_G(K)_P\cong K\backslash G(\A^\infty)/P(\A^\infty)\cong K_\Box\backslash G(\A_\Box)/P(\A_\Box),
 \]
 it suffices to show that $K_\Box\backslash G(\A_\Box)/P(\A_\Box)$ is a finite set.
 Let $K_\Box^0$ be the product of special compact open subgroups of $G(\Q_\ell)$ for $\ell\in \Box$.
 Then we have $K_\Box^0P_0(\A_\Box)=G(\A_\Box)$.
% It is known that there exists a compact open subgroup $K_\Box^0$ of $G(\A_\Box)$
% such that $K_\Box^0P_0(\A_\Box)=G(\A_\Box)$. 
 By shrinking $K_\Box$, we may assume that $K_\Box\subset K_\Box^0$.
 Then the map $K_\Box\backslash K_\Box^0\to K_\Box\backslash G(\A_\Box)/P(\A_\Box)$ is surjective, hence
 $K_\Box\backslash G(\A_\Box)/P(\A_\Box)$ is finite.

 Finally we prove \ref{item:P_G(K)-adjoint}. Note that admissible $\Q$-parabolic subgroups of $G$
 containing $P_0$ are in bijection with those of $G^\ad$ containing $P_0^\ad$.
 Therefore, we have only to show that $\mathcal{P}_G(K)_P\to \mathcal{P}_{G^\ad}(K'')_{P^\ad}$ is
 bijective for every admissible $\Q$-parabolic subgroup of $G$ containing $P_0$.
 Further, it is equivalent to the bijectivity of
 \[
 K_\Box\backslash G(\A_\Box)/P(\A_\Box)\xrightarrow{(*)} K_\Box^\ad\backslash G^\ad(\A_\Box)/P^\ad(\A_\Box).
 \]
 By \cite[15.1.4]{MR1642713}, we have
 \[
 G(\Q_\ell)/P(\Q_\ell)\cong (G_{\mathbb{Q}_{\ell}}/P_{\Q_\ell})(\Q_\ell)\cong (G^{\ad}_{\mathbb{Q}_{\ell}}/P^{\ad}_{\Q_\ell})(\Q_\ell)\cong G^\ad(\Q_\ell)/P^\ad(\Q_\ell)
 \]
 for each $\ell\in \Box$. Therefore the map $G(\A_\Box)/P(\A_\Box)\to G^\ad(\A_\Box)/P^\ad(\A_\Box)$
 is bijective. The bijectivity of $(*)$ easily follows from it.
\end{prf}

By the proof above, we also obtain the following:

\begin{cor}\label{cor:adm-parab-conj}
 The natural map $\mathcal{P}_{G,\Q}\to\mathcal{P}_{G,\A^\infty}$ is a bijection.
 In particular, for a compact open subgroup $K$ of $G(\A^\infty)$, we have a natural map
 $\mathcal{P}_G(K)\to \mathcal{P}_{G,\Q}$.
\end{cor}

\begin{prf}
 We use the notation in the proof of Lemma \ref{lem:P_G(K)}.
 By definition, the natural map $\mathcal{P}_{G,\Q}\to\mathcal{P}_{G,\A^\infty}$ is surjective.
 We shall show that it is injective. Take two admissible $\Q$-parabolic subgroups
 $P_1$, $P_2$ of $G$ containing $P_0$. If $P_1$ and $P_2$ are $G(\A^\infty)$-conjugate,
 then $[P_1]\in \mathcal{P}_G(K)_{P_1}\cap \mathcal{P}_G(K)_{P_2}$ for every compact open subgroup $K$
 of $G(\A^\infty)$. By the proof of Lemma \ref{lem:P_G(K)}, it implies that $P_1=P_2$.
 This completes the proof.
\end{prf}

\begin{prop}\label{prop:P_x-adm}
 For $x\in \Sh_K(G,X)_{E_v}^\ad(\cl)$, there uniquely exists an element $[P_x]\in \mathcal{P}_G(K)$
 which is mapped to $[P_{x,\Box}]$ under the injection $\mathcal{P}_G(K)\hookrightarrow \mathcal{P}_{G,\Box}(K_\Box)$
 in Lemma \ref{lem:P_G(K)} \ref{item:P_G(K)-local}. It is independent of $\Box$.
\end{prop}

To prove this proposition, we use the following lemma.

\begin{lem}\label{lem:mixed-Hodge}
 Let $(G,X)$ be a Shimura datum of Hodge type with an embedding $i\colon (G,X)\hookrightarrow (\GSp_{2n},X_{2n})$ 
 into a Siegel Shimura datum. Recall that $V$ denotes the standard representation of $\GSp_{2n}$.

 Let $W$ be a totally isotropic subspace of $V$, and define a filtration $W_\bullet V$ on $V$ as follows:
 \[
  W_0V=V,\quad W_{-1}V=W^\perp,\quad W_{-2}V=W,\quad W_{-3}V=0.
 \]
 We write $P$ for the stabilizer of $W_\bullet V$ in $G$.

 Let $L$ be a field of characteristic $0$. Assume that 
 the filtration $W_\bullet V\otimes_\Q L$ on $V_L=V\otimes_\Q L$ is $G_L$-split
 in the sense of \cite[(1.1.2)]{MR2669706}.
 Then, $P$ is an admissible $\Q$-parabolic subgroup of $G$.
 Further, if we write $P'$ for the stabilizer of $W_\bullet V$ in $\GSp_{2n}$,
 we have $P'=i_*P$ in the notation of \cite[2.1.28]{Mad-torHod}.
 Namely, the cocharacter of $\GSp_{2n}$ associated to $P'$ as in \cite[4.1]{MR1128753}
 is equal to $i\circ \lambda$, where $\lambda$ is the cocharacter of $G$ associated to $P$.
\end{lem}

\begin{prf}
 We write $U$ for the subgroup of $P$ consisting of elements acting on $\mathrm{gr}^W_\bullet V$ trivially,
 and $\nu$ for the cocharacter $\mathbb{G}_m\to \GL(\mathrm{gr}^W_\bullet V)$ determined from
 the grading on $\mathrm{gr}^W_\bullet V$.

 By \cite[Lemma 1.1.1]{MR2669706}, $P_L$ is a parabolic subgroup of $G_L$,
 $U_L$ is the unipotent radical of $P_L$, and the cocharacter
 $\nu_L\colon \mathbb{G}_m\to \GL(\mathrm{gr}^W_\bullet V\otimes_\Q L)$ over $L$
 factors through the closed subgroup $P_L/U_L$.
 Therefore, we conclude that $P$ is a parabolic subgroup of $G$,
 $U$ is the unipotent radical of $P$, and the cocharacter $\nu\colon \mathbb{G}_m\to \GL(\mathrm{gr}^W_\bullet V)$
 factors through $P/U$.
 This means that the filtration $W_\bullet V$ is $G$-split by \cite[Lemma 1.1.1]{MR2669706}.
 Take a cocharacter $w\colon \mathbb{G}_m\to G$ over $\Q$ which induces the filtration $W_\bullet V$ on $V$.
 It induces a filtration on the Tannakian category $\mathbf{Rep}_{\Q}(G)$.
 Let us prove that this filtration is Cayley in the sense of \cite[V, Definition 2.3]{MR1044823}.
 Take an arbitrary element $h\in X^+$. 
 By \cite[IV, Example 1.1 (c)]{MR1044823} (\cf \cite[4.2.1]{MR723182}), 
 $W_\bullet V$ and $h$ give a mixed Hodge structure on $V$.
 Therefore, \cite[IV, Proposition 1.3]{MR1044823} tells us that $w$ and $h$ define a mixed Hodge structure
 on $V_\xi$ for all objects $(\xi,V_\xi)$ of $\mathbf{Rep}_{\Q}(G)$.
 Hence the filtration induced from $w$ is Cayley, as desired.

 Now, by applying \cite[V, Proposition 2.4]{MR1044823} to each simple factor of 
 $(G^{\ad},X^{\ad})=(G_1,X_1) \times \cdots \times (G_r,X_r)$, we conclude that $P$ is admissible
 (we use the same argument as in the proof of Lemma \ref{lem:monodromy-filt-stab-l} \ref{item:P-ad-l} to
 pass to the adjoint group).
 The equality $P'=i_*P$ is proved in \cite[4.16]{MR1128753}.
\end{prf}

\begin{prf}[Proposition \ref{prop:P_x-adm}]
 Only the existence of $[P_x]$ requires a proof.
 By Lemma \ref{lem:monodromy-filt-stab-l} \ref{item:P-fin-ext-l}, 
 Lemma \ref{lem:monod-restr-p}, Lemma \ref{lem:preHod},
 Proposition \ref{prop:phi^ad-functoriality} and Lemma \ref{lem:P_G(K)} \ref{item:P_G(K)-adjoint},
 we may assume that $(G,X)$ is of Hodge type.
 We use the notation in Section \ref{subsec:Hodge-type}.
 By shrinking $K$, we may assume that $K_p\subset \widetilde{K}_{p,0}=\GSp_{2n}(\Z_p)$ and $K^p$ is small enough.
 We use the notation in Corollary \ref{cor:Hodge-parabolic}.
 Let $M_x=(G^{\natural},\underline{Y}\to G^{\natural}_{\kappa_x})$ be the degeneration datum
 corresponding to $\mathcal{A}_{\kappa^+_x}$ under the functor $M$ in \cite[Chapter III, Corollary 7.2]{MR1083353}.
 It gives a $1$-motive $M_{\overline{x}}$ over $\overline{\kappa}_x$ (\cf \cite[\S 10.1]{MR0498552}).
 For each prime $\ell\in \Box$, the $\ell$-adic realization $H_1(M_{\overline{x}},\Q_\ell)$ of $M_{\overline{x}}$
 is identified with $V_\ell\mathcal{A}_{\overline{x}}$, and equipped with the weight filtration
 $W_{\bullet,\overline{x},\ell}$, which coincides with the monodromy filtration
 $M_\bullet V_\ell\mathcal{A}_{\overline{x}}$ on $V_\ell\mathcal{A}_{\overline{x}}$ up to a shift.

 The $1$-motive $M_{\overline{x}}$ and the fixed isomorphism $\iota\colon \overline{\kappa}_x\xrightarrow{\cong}\C$
 gives rise to a $1$-motive $M_{\iota\overline{x}}$ over $\C$.
 Its Betti realization $H_1(M_{\iota\overline{x}},\Q)$ is naturally isomorphic to $H_1(A,\Q)$
 (recall that $(A,(t_\alpha)_{\alpha\in J},\eta K)$ denotes the triple corresponding to the $\C$-point $\iota\overline{x}$).
 We denote the weight filtration on it by $W_{\bullet,\iota\overline{x},\Q}$.
 It is known that $W_{-2,\iota\overline{x},\Q}$ is totally isotropic and $W_{-1,\iota\overline{x},\Q}=W_{-2,\iota\overline{x},\Q}^\perp$
 with respect to the polarization $\pm t_0$ on $A$.
 For $\ell\in \Box$, we write $\varepsilon_\ell$ for the comparison isomorphism
 $H_1(A,\Q)\otimes_\Q\Q_\ell\xrightarrow{\cong} V_\ell A$.
 The composite 
 \[
  H_1(A,\Q)\otimes_\Q\Q_\ell\xrightarrow[\cong]{\varepsilon_\ell} V_\ell A\xrightarrow[\cong]{\iota^{-1}} V_\ell\mathcal{A}_{\overline{x}}=H_1(M_{\overline{x}},\Q_\ell)
 \]
 carries the filtration $W_{\bullet,\iota\overline{x},\Q}\otimes_\Q\Q_\ell$ onto $W_{\bullet,\overline{x},\ell}$
 defined above.

 We take a representative $\eta\colon V_{\A^\infty}\xrightarrow{\cong}V^\infty A$ of the $K$-orbit $\eta K$ and
 an isomorphism $\eta_\Q\colon V\xrightarrow{\cong}H_1(A,\Q)$ as in the condition $(*)$
 of Proposition \ref{prop:Hodge-C-moduli}.
 For each prime $\ell$, we can easily observe that 
 $g_\ell=\eta_\ell^{-1}\circ \varepsilon_\ell\circ (\eta_{\Q}\otimes\Q_\ell)$
 preserves the tensors $(s_\alpha)_{\alpha\in J'}$. Therefore $g_\ell$ lies in $G(\Q_\ell)$.

 We put $W_\bullet V=\eta_\Q^{-1}(W_{\bullet,\iota\overline{x},\Q})$
 (it depends on the choice of $\overline{x}$, $\iota$ and $\eta_\Q$),
 and denote by $P'$ the stabilizer of $W_\bullet V$ in $G$.
 For $\ell\in \Box$, we have
 \[
  (\eta_\ell^{-1}\circ\iota)(W_{\bullet,\overline{x},\ell})=(\eta_\ell^{-1}\circ \varepsilon_\ell)(W_{\bullet,\iota\overline{x},\Q}\otimes_\Q\Q_\ell)=g_\ell(W_\bullet V\otimes_\Q\Q_\ell).
 \]
 By Corollary \ref{cor:Hodge-parabolic} \ref{item:monodromy-parabolic}, $g_\ell(W_\bullet V\otimes_\Q\Q_\ell)$
 is $G_{\Q_\ell}$-split. Hence $W_\bullet V\otimes_\Q\Q_\ell$ is also $G_{\Q_\ell}$-split.
 Therefore, Lemma \ref{lem:mixed-Hodge} tells us that $P'$ is an admissible $\Q$-parabolic subgroup of $G$.
 Further, by the above equality, $P_{x,\eta,\ell}$ 
 in Corollary \ref{cor:Hodge-parabolic} \ref{item:monodromy-parabolic} is equal to $g_\ell P'_{\Q_\ell}g_\ell^{-1}$.

 Set $g=(g_\ell)_{\ell\in \Box}\times 1\in G(\A^\infty)=G(\A_\Box)\times G(\A^{\infty,\Box})$.
 Then $P_x=gP'_{\A^\infty}g^{-1}$ is an admissible $\A^\infty$-parabolic subgroup of $G$,
 and its image $[P_x]$ in $\mathcal{P}_K(G)$ is mapped to $[\prod_{\ell\in \Box}P_{x,\eta,\ell}]=[P_{x,\Box}]$
 under the injection $\mathcal{P}_G(K)\hookrightarrow \mathcal{P}_{G,\Box}(K_\Box)$
 (for the last equality, see Corollary \ref{cor:Hodge-parabolic} \ref{item:monodromy-parabolic}).
 This completes the proof.
\end{prf}

\begin{cor}\label{cor:pg-point}
 For $x\in \Sh_K(G,X)_{E_v}^\ad(\cl)$ and a prime number $\ell\neq p$, the following are equivalent:
 \begin{enumerate}
  \item $\phi_{x,\ell}^\ad$ is potentially unramified.
  \item $\phi_{x,p}^\ad$ is potentially crystalline.
  \item $[P_x]=[G]$.
 \end{enumerate}
 If the above conditions are satisfied, we say that $x$ is of potentially good reduction.
\end{cor}

\begin{prf}
 By Lemma \ref{lem:monodromy-filt-stab-l} \ref{item:P=G-l}, $\phi_{x,\ell}^\ad$ is potentially unramified
 if and only if $P_{\phi_{x,\ell}^\ad}=G^\ad_{\Q_\ell}$. By definition this is clearly equivalent to
 $[P_x]=[G]$. Similarly we can prove the equivalence of (ii) and (iii).
\end{prf}

\begin{defn}
 We denote the map $\Sh_K(G,X)_{E_v}^\ad(\cl)\to \mathcal{P}_G(K)$; $x\mapsto [P_x]$ by $\Phi_K$.
\end{defn}

\begin{prop}\label{prop:Phi-functoriality}
 Let $(G,X)$ be a Shimura datum of preabelian type.
 \begin{enumerate}
  \item\label{item:Phi-Hecke}
       The map $\Phi_K$ is Hecke-equivariant in the following sense.
       Let $K$, $K'$ be compact open subgroups of $G(\A^\infty)$, and $g$ an element of $G(\A^\infty)$
       such that $g^{-1}Kg\subset K'$. Then the diagram
       \[
       \xymatrix{%
       \Sh_K(G,X)_{E_v}^\ad(\cl)\ar[r]^-{\Phi_K}\ar[d]^-{g}&\mathcal{P}_G(K)\ar[d]^-{[P]\mapsto [g^{-1}Pg]}\\
       \Sh_{K'}(G,X)_{E_v}^\ad(\cl)\ar[r]^-{\Phi_{K'}}&\mathcal{P}_{G}(K')
       }
       \]
       is commutative.
  \item\label{item:Phi-ad}
       Let $K$ be a compact open subgroup of $G(\A^\infty)$, and $K''$ a compact open subgroup of
       $G^\ad(\A^\infty)$ containing the image of $K$.
       We write $E^\ad$ for the reflex field of $(G^\ad,X^\ad)$ and $v^\ad$ the place of $E^\ad$ below $v$.
       Then the diagram
       \[
       \xymatrix{%
       \Sh_K(G,X)_{E_v}^\ad(\cl)\ar[r]^-{\Phi_K}\ar[d]&\mathcal{P}_G(K)\ar[d]^-{[P]\mapsto [P^\ad]}\\
       \Sh_{K''}(G^\ad,X^\ad)_{E^\ad_{v^\ad}}^\ad(\cl)\ar[r]^-{\Phi_{K''}}&\mathcal{P}_{G^\ad}(K'')
       }
       \]
       is commutative.
  \item\label{item:Phi-Hodge}
       Assume that $(G,X)$ is of Hodge type, and let $i\colon (G,X)\hookrightarrow (\GSp_{2n},X_{2n})$
       be an embedding into a Siegel Shimura datum.
       Let $K$ be a compact open subgroup of $G(\mathbb{A}^{\infty})$,
       and $\widetilde{K}$ a compact open subgroup of $\GSp_{2n}(\mathbb{A}^{\infty})$ containing $K$.
       Then the diagram
       \[
       \xymatrix{%
       \Sh_K(G,X)_{E_v}^\ad(\cl)\ar[r]^-{\Phi_K}\ar[d]&\mathcal{P}_G(K)\ar[d]^-{[P]\mapsto [i_*P]}\\
       \Sh_{\widetilde{K}}(\GSp_{2n},X_{2n})_{\Q_p}^\ad(\cl)\ar[r]^-{\Phi_{\widetilde{K}}}&\mathcal{P}_{\GSp_{2n}}(\widetilde{K})
       }
       \]
       is commutative.
 \end{enumerate}
\end{prop}

\begin{prf}
 The assertions \ref{item:Phi-Hecke} and \ref{item:Phi-ad} are immediate consequences of 
 Proposition \ref{prop:phi^ad-functoriality} and Lemma \ref{lem:P_G(K)} \ref{item:P_G(K)-local}.
 The claim \ref{item:Phi-Hodge} follows from Lemma \ref{lem:mixed-Hodge} and the construction of $[P_x]$
 in the proof of Proposition \ref{prop:P_x-adm}.
\end{prf}

In the remaining part of this section, we will prove that the map $\Phi_K$ comes from a partition
of $\Sh_K(G,X)^\ad_{E_v}$ into locally closed constructible subsets.

\begin{thm}\label{thm:Sh-partition}
 For each $[P]\in \mathcal{P}_G(K)$, there uniquely exists a locally closed constructible subset 
 $C_{[P]}$ of $\Sh_K (G,X)_{E_v}^{\ad}$ such that 
 \[
 C_{[P]}(\cl)=\{x\in \Sh_K(G,X)_{E_v}^{\ad}(\cl)\mid \Phi_K(x)=[P]\}. 
 \]
 Furthermore, the subset $C_{[G]}$ is open and quasi-compact.
\end{thm}

\begin{rem}\label{rem:C_P-cover}
 The subsets $\{C_{[P]}\}_{[P]\in \mathcal{P}_G(K)}$ in Theorem \ref{thm:Sh-partition} are mutually disjoint
 and cover $\Sh_K(G,X)^\ad_{E_v}$.
 Indeed, by Lemma \ref{lem:cl-pt-constructible} (i), it can be checked at the level of classical points,
 which is obvious.
\end{rem}
Theorem \ref{thm:Sh-partition} will be proved in Section \ref{subsec:existence-partition}.
Admitting this theorem, we have the following definition.
 
\begin{defn}\label{defn:pg-locus}
 We put $\Sh_K(G,X)_{E_v}^{\pg}=C_{[G]}$, and call it the potentially good reduction locus of
 $\Sh_K(G,X)^\ad_{E_v}$.
 It is a quasi-compact open subset of $\Sh_K(G,X)^\ad_{E_v}$ characterized by the following property:
 \begin{itemize}
  \item for $x\in \Sh_K(G,X)^\ad_{E_v}(\cl)$, $x$ lies in $\Sh_K(G,X)_{E_v}^{\pg}$ if and only if
	$x$ is of potentially good reduction in the sense of Corollary \ref{cor:pg-point}.
 \end{itemize}
\end{defn}

\begin{exa}
 When $\Sh_K(G,X)$ is proper over $E$, $G$ has no proper parabolic subgroup defined over $\Q$.
 Therefore, we have $\mathcal{P}_G(K)=\{[G]\}$ and $\Sh_K(G,X)_{E_v}^{\pg}=\Sh_K(G,X)_{E_v}^{\ad}$.
\end{exa}

\subsection{Partition in the Siegel case}
In this subsection, we will give a proof of Theorem \ref{thm:Sh-partition} in the Siegel case.
We use the notation in Section \ref{sec:Siegel}.
In particular, recall that $(V,\langle\ ,\ \rangle)$ is a symplectic space of dimension $2n$ over $\Q$
and $L$ is a self-dual $\Z$-lattice of $V$. For simplicity, we write $\mathcal{S}_K$ for
$\Sh_K(\GSp_{2n},X_{2n})_{\Q_p}^\ad$.

It is well-known that conjugacy classes of maximal parabolic subgroups of $\GSp_{2n}$
are parametrized by integers $0\le r\le n$; the class corresponding to $r$ consists of
parabolic subgroups obtained as stabilizers of $r$-dimensional totally isotropic subspaces of $V$. 
Namely, $\mathcal{P}_{\GSp_{2n},\Q}\cong \{0,1,\ldots,n\}$ under the notation in Definition \ref{defn:adm-parab}.
We write $\mathcal{P}_{\GSp_{2n}}(K)_r$ for the inverse image of $r$
under the map 
\[
 \mathcal{P}_{\GSp_{2n}}(K)\to \mathcal{P}_{\GSp_{2n},\Q}\cong \{0,1,\ldots,n\}
\]
(see Corollary \ref{cor:adm-parab-conj}).
Clearly we have $\mathcal{P}_{\GSp_{2n}}(K)=\coprod_{0\le r\le n}\mathcal{P}_{\GSp_{2n}}(K)_r$.
We put $\mathcal{P}_{\GSp_{2n}}(K)_{\le r}=\coprod_{r'\le r}\mathcal{P}_{\GSp_{2n}}(K)_{r'}$.
Note that $\mathcal{P}_{\GSp_{2n}}(K)_{0}=\mathcal{P}_{\GSp_{2n}}(K)_{\le 0}=\{[G]\}$.

\begin{prop}\label{prop:Siegel-partition-r}
 There uniquely exists a constructible open subset $\mathcal{S}_{K,\mathopen{]}\le r\mathclose{[}}$
 of $\mathcal{S}_K$ such that $x\in \mathcal{S}_K(\cl)$ belongs to $\mathcal{S}_{K,\mathopen{]}\le r\mathclose{[}}$
 if and only if $\Phi_K(x)\in \mathcal{P}_{\GSp_{2n}}(K)_{\le r}$.
 Moreover, $\mathcal{S}_{K,\mathopen{]}\le 0\mathclose{[}}$ is quasi-compact.
\end{prop}

\begin{prf}
 By Proposition \ref{prop:Phi-functoriality} \ref{item:Phi-Hecke}, we may shrink $K$ freely.
 Therefore, we may assume that $K_p$ is contained in $K_{p,0}=\GSp_{2n}(\Z_p)$.

 For $r\ge 0$, let $\mathscr{S}_{K^p,\le r}^{\mathrm{tor}}$ 
 be the subset of $\mathscr{S}_{K^p}^{\mathrm{tor}}$ consisting of $x \in \mathscr{S}_{K^p}^{\mathrm{tor}}$ 
 such that the toric rank of the semi-abelian variety $\mathcal{A}_x$ is at most $r$. 
 By \cite[Lemma 3.3.1.4]{Kai-Wen}, it is an open subset of $\mathscr{S}_{K^p}^{\mathrm{tor}}$. 

 Since $\mathscr{S}_{K^p}^{\mathrm{tor}}$ is proper over $\mathbb{Z}_p$, 
 we may consider the specialization map 
 \[
 \spp \colon (\mathscr{S}^{\mathrm{tor}}_{K^p,\Q_p})^\ad
 =(\mathscr{S}_{K^p}^{\mathrm{tor}})^\ad 
 =(\mathscr{S}_{K^p}^{\mathrm{tor}})^{\wedge\rig}
 \to \mathscr{S}^{\mathrm{tor}}_{K^p,\mathbb{F}_p} 
 \]
 introduced in Section \ref{subsec:notation-adic} (for the second equality, see Lemma \ref{lem:rig-ad} (ii)). 
 For an integer $r\ge 0$, we put
 $\mathcal{S}_{K_{p,0}K^p,\mathopen{]}\le r\mathclose{[}}=\spp^{-1}(\mathscr{S}^\mathrm{tor}_{K^p,\le r,\mathbb{F}_p}) \cap \mathscr{S}_{K^p,\Q_p}^\ad$.
 It is a constructible open subset of $\mathscr{S}_{K^p,\Q_p}^\ad=\mathcal{S}_{K_{p,0}K^p}$.
 Further, let $\mathcal{S}_{K,\mathopen{]}\le r\mathclose{[}}$ be the inverse image of
 $\mathcal{S}_{K_{p,0} K^p, \mathopen{]}\le r\mathclose{[}}$ under the natural morphism
 $\mathcal{S}_K \to \mathcal{S}_{K_{p,0}K^p}$.
 By Corollary \ref{cor:Hodge-parabolic} \ref{item:Siegel-toric}, it satisfies the desired property.
 Since $\mathcal{S}_{K_{p,0}K^p,\mathopen{]}\le 0\mathclose{[}}=\spp^{-1}(\mathscr{S}_{K^p,\mathbb{F}_p})$,
 it is quasi-compact. Hence $\mathcal{S}_{K,\mathopen{]}\le 0\mathclose{[}}$ is also quasi-compact.
\end{prf}

For $0\le r\le n$, we put
$\mathcal{S}_{K,\mathopen{]}r\mathclose{[}}=\mathcal{S}_{K,\mathopen{]}\le r\mathclose{[}}\setminus \mathcal{S}_{K,\mathopen{]}\le r-1\mathclose{[}}$,
where $\mathcal{S}_{K,\mathopen{]}\le -1\mathclose{[}}$ means $\varnothing$.
It is a locally closed constructible subset of $\mathcal{S}_K$.

\begin{lem}\label{lem:conn-comp-finite}
 The set $\pi_0(\mathcal{S}_{K,\mathopen{]}r\mathclose{[}})$ of connected components of
 $\mathcal{S}_{K,\mathopen{]}r\mathclose{[}}$ is finite, and
 consists of locally closed constructible subsets of $\mathcal{S}_K$.
 Hence $\mathcal{S}_{K,\mathopen{]}r\mathclose{[}}$ is topologically the disjoint union of
 elements of $\pi_0(\mathcal{S}_{K,\mathopen{]}r\mathclose{[}})$.
\end{lem}

\begin{prf}
 By shrinking $K$, we may assume that $K_p$ is contained in $K_{p,0}=\GSp_{2n}(\Z_p)$.
 Then the claim is the case $X=\mathscr{S}_{K^p,\Q_p}^{\mathrm{tor}}$,
 $U=\mathscr{S}_{K^p,\Q_p}$, $U'=\Sh_K(\GSp_{2n},X_{2n})_{\Q_p}$,
 and $L=\spp^{-1}(\mathscr{S}^\mathrm{tor}_{K^p,\le r,\mathbb{F}_p})\setminus \spp^{-1}(\mathscr{S}^\mathrm{tor}_{K^p,\le r-1,\mathbb{F}_p})$
 of the subsequent general lemma.
\end{prf}

\begin{lem}
 Let $F$ be a $p$-adic field. Let $X$ be a purely $d$-dimensional proper smooth scheme over $F$,
 and $Y$ a closed subscheme of $X$ whose dimension is less than $d$.
 We put $U=X\setminus Y$, and consider a finite \'etale surjection $f\colon U'\to U$.

 For a locally closed constructible subset $L$ of $X^\ad$, we put $L_{U'}=(f^{\ad})^{-1}(L\cap U^\ad)$.
 Then, the set $\pi_0(L_{U'})$ is finite, and consists of locally closed constructible subsets of $U'^\ad$.
 In particular, $L_{U'}$ is topologically the disjoint union of elements of $\pi_0(L_{U'})$.
\end{lem}

\begin{prf}
 Let $X'$ be the normalization of $X$ in $U'$.
 By the resolution of singularities, there exists a purely $d$-dimensional proper smooth scheme $X''$ over $F$
 and a proper birational morphism $\phi\colon X''\to X'$ which induces an isomorphism
 $\phi^{-1}(U')\xrightarrow{\cong}U'$.
 Let us denote the composite $X''\xrightarrow{\phi}X'\to X$ by $\phi'$.
 By replacing $X$, $Y$, $L$ with $X''$, $X''\times_XY$, $(\phi'^{\ad})^{-1}(L)$ respectively,
 we may assume that $U'=U$.

 We denote by $L^\circ$ (resp.\ $L_U^\circ$) the interior of $L$ (resp.\ $L_U$) in $X^\ad$ (resp.\ $U^\ad$).
 Clearly we have 
 $L_U^\circ=L^\circ\setminus Y^\ad$. We fix a prime number $\ell\neq p$, and
 consider the following commutative diagram:
 \[
  \xymatrix{%
 H^0(L_{\overline{F}},\F_\ell)\ar[r]\ar[d]^-{(1)}& H^0(L^\circ_{\overline{F}},\F_\ell)\ar[d]^-{(2)}\\
 H^0(L_{U,\overline{F}},\F_\ell)\ar[r]& H^0(L^\circ_{U,\overline{F}},\F_\ell)\lefteqn{.}
 }
 \]
 Since $X^\ad$ is proper of finite type over $\Spa(F,\mathcal{O}_F)$,
 $H^0(L_{\overline{F}},\F_\ell)$ is a finite-dimensional
 $\F_\ell$-vector space by \cite[Proposition 3.16 i)]{MR1620118}. Therefore, to show the finiteness of
 $\pi_0(L_U)$, it suffices to prove that the map $(1)$ is
 an isomorphism. On the other hand, by \cite[Theorem 3.7]{MR1620118}, the horizontal maps are isomorphisms.
 Hence it suffices to prove that (2) is an isomorphism.
 Note that, by \cite[Lemma 1.3 iii)]{MR1620114}, $L^\circ$ and $L_U^\circ$ are taut, and then
 $L^\circ\cap Y^\ad$ is also taut.
 Therefore we can consider the compactly supported cohomology of these spaces.
 Since $\dim (L^\circ\cap Y^\ad)<d$,
 we have $H^{2d-1}_c((L^\circ\cap Y^\ad)_{\overline{F}},\F_\ell)=H^{2d}_c((L^\circ\cap Y^\ad)_{\overline{F}},\F_\ell)=0$.
 This implies that the natural map 
 $H^{2d}_c(L^\circ_{U,\overline{F}},\F_\ell)\to H^{2d}_c(L^\circ_{\overline{F}},\F_\ell)$
 is an isomorphism. By the Poincar\'e duality, we conclude that the map (2) is an isomorphism.

 By the finiteness of $\pi_0(L_U)$, every element $C$ of $\pi_0(L_U)$ is
 an open and closed subset of $L_U$.
 Since $L$ is locally closed constructible, so is $C$.
\end{prf}

\begin{lem}\label{lem:conn-comp-Psi}
 There uniquely exists a map $\Psi_K\colon \pi_0(\mathcal{S}_{K,\mathopen{]}r\mathclose{[}})\to\mathcal{P}_G(K)$
 satisfying the following: 
 for $C\in \pi_0 (\mathcal{S}_{K,\mathopen{]}r\mathclose{[}})$ and $x\in C(\cl)$,
 we have $\Psi_K(C) =\Phi_K(x)$.
\end{lem}

\begin{prf}
 Let $C\in \pi_0(\mathcal{S}_{K,\mathopen{]}r\mathclose{[}})$.
 Then, by Lemma \ref{lem:conn-comp-finite} and Lemma \ref{lem:cl-pt-constructible} (i),
 we have $C(\cl)\neq \varnothing$.
 Therefore, it suffices to show that $\Phi_K(x)$ is independent of the choice of $x\in C(\cl)$. 
 By Proposition \ref{prop:Phi-functoriality} \ref{item:Phi-Hecke}, we may assume that $K_p$ is contained
 in $K_{p,0}=\GSp_{2n}(\Z_p)$. Recall that in this case 
 $\mathcal{S}_{K,\mathopen{]}\le r\mathclose{[}}$ is obtained as
 the inverse image under $\mathcal{S}_K\to \mathcal{S}_{K_{p,0}K^p}$ of 
 $\mathcal{S}_{K_{p,0}K^p,\mathopen{]}\le r\mathclose{[}}$,
 which is equal to
 $\spp^{-1}(\mathscr{S}^\mathrm{tor}_{K^p,\le r,\mathbb{F}_p}) \cap \mathscr{S}_{K^p,\Q_p}^\ad$
 (see the proof of Proposition \ref{prop:Siegel-partition-r}).
 
 Let $\mathcal{S}_{K,\mathopen{]}r\mathclose{[}}^\circ$ be the interior of
 $\mathcal{S}_{K,\mathopen{]}r\mathclose{[}}$ in $\mathcal{S}_K$.
 Then, the inverse image of $\mathcal{S}_{K_{p,0}K^p,\mathopen{]}r\mathclose{[}}^\circ$ under
 $\mathcal{S}_K\to \mathcal{S}_{K_{p,0}K^p}$ equals $\mathcal{S}_{K,\mathopen{]}r\mathclose{[}}^\circ$.
 We write $C^\circ$ for the interior of $C$ in $\mathcal{S}_K$.
 It is connected by \cite[Theorem 3.7]{MR1620118}, and included in
 $\mathcal{S}_{K,\mathopen{]}r\mathclose{[}}^\circ$.
 Since $C$ is constructible in $\mathcal{S}_K$, we have $C(\cl)=C^\circ(\cl)$
 by Lemma \ref{lem:cl-pt-constructible} (ii).
 Hence it suffices to show that $\Phi_K(x)$ is independent of the choice of $x\in C^\circ(\cl)$. 

 We apply the construction introduced in Section \ref{subsec:etale-sheaves-semi-abelian}
 to the case where $S=\mathscr{S}^\mathrm{tor}_{K^p,\le r}$, $S_0=\mathscr{S}^\mathrm{tor}_{K^p,r,\mathbb{F}_p}$, 
 $U=\mathscr{S}_{K^p,\Q_p}$ and $G=\mathcal{A}\vert_{\mathscr{S}^\mathrm{tor}_{K^p,\le r}}$.
 We write $\widehat{\mathcal{T}}$ for the corresponding $\widehat{T}$.
 By Lemma \ref{lem:rig-gen-fiber}, 
 $U^\ad=U\times_S t(\mathcal{S})_a=\mathcal{S}_{K_{p,0}K^p,\mathopen{]}r\mathclose{[}}^\circ$ 
 in this case. 
 Therefore, for each $m\ge 0$ and a prime $\ell$, we have three locally constant constructible sheaves 
 \[
 \widehat{\mathcal{T}}^\rig[\ell^m]_{\mathcal{S}_{K_{p,0}K^p,\mathopen{]}r\mathclose{[}}^\circ}\subset
 \widehat{\mathcal{A}}^\rig[\ell^m]_{\mathcal{S}_{K_{p,0}K^p,\mathopen{]}r\mathclose{[}}^\circ}\subset
 \mathcal{A}^\ad[\ell^m]_{\mathcal{S}_{K_{p,0}K^p,\mathopen{]}r\mathclose{[}}^\circ}.
\]
 We put
 \begin{align*}
  \mathcal{V}_\ell&=\Bigl(\varprojlim_{m}\mathcal{A}^\ad[\ell^m]_{\mathcal{S}_{K_{p,0}K^p,\mathopen{]}r\mathclose{[}}^\circ}\vert_{C^\circ}\Bigr)\otimes\Q_\ell,\quad
  \mathcal{F}_\ell=\Bigl(\varprojlim_{m}\widehat{\mathcal{A}}^\rig[\ell^m]_{\mathcal{S}_{K_{p,0}K^p,\mathopen{]}r\mathclose{[}}^\circ}\vert_{C^\circ}\Bigr)\otimes\Q_\ell,\\
  \mathcal{T}_\ell&=\Bigl(\varprojlim_{m}\widehat{\mathcal{T}}^\rig[\ell^m]_{\mathcal{S}_{K_{p,0}K^p,\mathopen{]}r\mathclose{[}}^\circ}\vert_{C^\circ}\Bigr)\otimes\Q_\ell.
 \end{align*}
 They are smooth $\ell$-adic sheaves over $C^\circ$.

 Now we use the moduli interpretation with rational level structures of $\mathcal{S}_K$.
 Fix a geometric point $\overline{x}_0$ of $C^\circ$, and let $\eta K$ be the 
 $\pi_1(C^\circ,\overline{x}_0)$-invariant $K$-orbit of isomorphisms
 $V_{\A^\infty}\xrightarrow{\cong} V^\infty\mathcal{A}_{\overline{x}_0}$ corresponding to
 the universal level structure on $\mathcal{A}\vert_{C^\circ}$.
 For $x\in C^\circ(\cl)$, the rational $K$-level structure on $\mathcal{A}_x$
 corresponding to $x$ itself is obtained in the following manner.
 Fix a geometric point $\overline{x}$ lying over $x$. Since $C^\circ$ is connected,
 there exists an isomorphism $\pi_1(C^\circ,\overline{x}_0)\to \pi_1(C^\circ,\overline{x})$, which is
 canonical up to $\pi_1(C^\circ,\overline{x}_0)$-conjugacy. If we fix such an isomorphism,
 for a smooth sheaf $\mathcal{G}$ on $C^\circ$, we have a functorial isomorphism
 $\mathcal{G}_{\overline{x}_0}\xrightarrow{\cong}\mathcal{G}_{\overline{x}}$ compatible with the $\pi_1$-actions.
 In particular, the smooth sheaf $(\varprojlim_{N}\mathcal{A}^\ad[N]_{C^\circ})\otimes\Q$ determines
 an isomorphism $V^\infty\mathcal{A}_{\overline{x}_0}\xrightarrow{\cong}V^\infty\mathcal{A}_{\overline{x}}$.
 By composing it with each element of $\eta K$, we obtain a $K$-orbit of 
 isomorphisms $V_{\A^\infty}\xrightarrow{\cong}V^\infty\mathcal{A}_{\overline{x}}$, which turns out to be
 $\pi_1(C^\circ,\overline{x})$-invariant. Since the action of $\pi_1(x,\overline{x})$
 on $V^\infty\mathcal{A}_{\overline{x}}$ factors through $\pi_1(x,\overline{x})\to \pi_1(C^\circ,\overline{x})$,
 this orbit gives a rational $K$-level structure on $\mathcal{A}_x$.

 Fix a representative $\eta$ of $\eta K$ and write $\eta_x$ for the composite of $\eta$
 and the isomorphism $V^\infty\mathcal{A}_{\overline{x}_0}\xrightarrow{\cong}V^\infty\mathcal{A}_{\overline{x}}$.
 We take a prime number $\ell$ and consider the $\ell$-part 
 \[
  \eta_{x,\ell}\colon V_{\Q_\ell}\xrightarrow[\eta_\ell]{\cong}V_\ell\mathcal{A}_{\overline{x}_0}\xrightarrow[(*)]{\cong}V_\ell\mathcal{A}_{\overline{x}}
 \]
 of $\eta_x$. Note that the isomorphism $(*)$ is given by the smooth $\ell$-adic sheaf $\mathcal{V}_\ell$
 introduced above. Moreover, by Corollary \ref{cor:Hodge-parabolic} \ref{item:semiab-red},
 Proposition \ref{prop:Afilt}, Corollary \ref{cor:semiab-monodromy-filt} and its proof,
 the monodromy filtrations on $V_\ell\mathcal{A}_{\overline{x}_0}$ and $V_\ell\mathcal{A}_{\overline{x}}$
 are given by
 \[
 0\subset \mathcal{T}_{\ell,\overline{x}_0}\subset \mathcal{F}_{\ell,\overline{x}_0}\subset \mathcal{V}_{\ell,\overline{x}_0},\quad 0\subset \mathcal{T}_{\ell,\overline{x}}\subset \mathcal{F}_{\ell,\overline{x}}\subset \mathcal{V}_{\ell,\overline{x}},
 \]
 respectively. Since $\mathcal{F}_\ell$ and $\mathcal{T}_\ell$ are smooth sheaves,
 the isomorphism $(*)$ carries the first filtration to the second. Hence we have
 $\eta_{x,\ell}^{-1}(M_\bullet V_\ell\mathcal{A}_{\overline{x}})=\eta_\ell^{-1}(M_\bullet V_\ell\mathcal{A}_{\overline{x}_0})$, which is independent of $x\in C^\circ(\cl)$.
 Therefore, the parabolic subgroup $P_{x,\iota_x\circ\eta_x,\ell}$ in
 Corollary \ref{cor:Hodge-parabolic} \ref{item:monodromy-parabolic}, 
 where $\iota_x\colon \overline{\kappa}_x\xrightarrow{\cong}\C$ is a fixed isomorphism,
 is also independent of $x$. By Corollary \ref{cor:Hodge-parabolic} \ref{item:monodromy-parabolic}
 and Lemma \ref{lem:P_G(K)} \ref{item:P_G(K)-local}, we conclude that
 $\Phi_K(x)=[P_x]\in \mathcal{P}_{\GSp_{2n}}(K)$ is independent of $x$.
\end{prf}

Now we can prove Theorem \ref{thm:Sh-partition} in the Siegel case.

\begin{prf}[of Theorem \ref{thm:Sh-partition} for the Siegel case]
 For $[P]\in\mathcal{P}_{\GSp_{2n}}(K)$, take $0\le r\le n$ such that $[P]$ lies in $\mathcal{P}_{\GSp_{2n}}(K)_r$.
 We put 
 \[
  C_{[P]}=\bigcup_{\substack{C\in\pi_0(\mathcal{S}_{K,\mathopen{]}r\mathclose{[}}),\\\Psi_K(C)=[P]}}C.
 \]
 It is a constructible subset of $\mathcal{S}_K$ by Lemma \ref{lem:conn-comp-finite}.
 Since each $C\in \pi_0(\mathcal{S}_{K,\mathopen{]}r\mathclose{[}})$ is open in
 $\mathcal{S}_{K,\mathopen{]}r\mathclose{[}}$, $C_{[P]}$ is a locally closed subset of $\mathcal{S}_K$.
 We can also check that $x\in\mathcal{S}_K(\cl)$ lies in $C_{[P]}$ if and only if $\Phi_K(x)=[P]$.
 We have already checked in Proposition \ref{prop:Siegel-partition-r}
 that $C_{[G]}=\mathcal{S}_{K,\mathopen{]}\le 0\mathclose{[}}$ is quasi-compact open.
\end{prf}

\begin{cor}\label{cor:C_P-disjoint-Siegel}
 For $[P]\in \mathcal{P}_G(K)_r$ and $[P']\in \mathcal{P}_G(K)_{r'}$, assume that
 $r>r'$ or $r=r'$ and $[P]\neq [P']$.
 Then we have $C_{[P]}^-\cap C_{[P']}=\varnothing$,
 where $C_{[P]}^-$ denotes the closure of $C_{[P]}$ in $\mathcal{S}_K$.
\end{cor}

 \begin{prf}
  First assume that $r>r'$. Since the complement $\mathcal{S}^c_{K,\mathopen{]}\le r'\mathclose{[}}$
  of $\mathcal{S}_{K,\mathopen{]}\le r'\mathclose{[}}$ is a closed subset of $\mathcal{S}_K$,
  we have
  \[
  C_{[P]}^-\cap C_{[P']}\subset \mathcal{S}^c_{K,\mathopen{]}\le r'\mathclose{[}}\cap \mathcal{S}_{K,\mathopen{]}\le r'\mathclose{[}}=\varnothing.
  \]
  Next assume that $r=r'$ and $[P]\neq [P']$.
  By construction, $C_{[P]}$ is closed in $\mathcal{S}_{K,\mathopen{]}r\mathclose{[}}$.
  Therefore, we have 
 \[
  C_{[P]}^-\cap C_{[P']}=(C_{[P]}^-\cap \mathcal{S}_{K,\mathopen{]}r\mathclose{[}})\cap C_{[P']}=C_{[P]}\cap C_{[P']}=\varnothing.
 \]
 For the last equality, see Remark \ref{rem:C_P-cover}.
 \end{prf}

\subsection{Existence of partition}\label{subsec:existence-partition}
In this subsection, we complete the proof of Theorem \ref{thm:Sh-partition} by reducing to the Siegel case.
First we consider the Hodge type case.

\begin{lem}\label{lem:P_G(K)-Hodge}
Let $(G,X)$ be a Shimura datum of Hodge type with an embedding $i\colon (G,X)\hookrightarrow (\GSp_{2n},X_{2n})$ 
into a Siegel Shimura datum. 
For a compact open subgroup $K$ of $G(\A^\infty)$, 
there exists a compact open subgroup $\widetilde{K}$ of $\GSp_{2n}(\A^\infty)$ 
containing $K$ such that the map $\mathcal{P}_G(K)\to \mathcal{P}_{\GSp_{2n}}(\widetilde{K});$
$[P]\mapsto [i_*P]$ is injective.
\end{lem}

\begin{prf}
 It suffices to prove the following claim:
 \begin{quote}
  for $[P_1],[P_2]\in \mathcal{P}_G(K)$ with $[P_1]\neq [P_2]$, there exists a compact open subgroup
  $\widetilde{K}_{[P_1],[P_2]}$ of $\GSp_{2n}(\A^\infty)$ containing $K$ such that
  $[i_*P_1]\neq [i_*P_2]$ in $\mathcal{P}_{\GSp_{2n}}(\widetilde{K}_{[P_1],[P_2]})$.
 \end{quote}
 Indeed, the intersection of $\widetilde{K}_{[P_1],[P_2]}$ for all pairs $([P_1],[P_2])$
 with $[P_1]\neq [P_2]$ satisfies the desired condition.

 Fix a compact open subgroup $\widetilde{K}_0$ of $\GSp_{2n}(\A^\infty)$ containing $K$.
 Take representatives $P_1$, $P_2$ of $[P_1]$, $[P_2]$, respectively.
 We put $Z=\{g\in \widetilde{K}_0\mid g(i_*P_1)g^{-1}=i_*P_2\}$, which is
 clearly a closed subset of $\widetilde{K}_0$.
 Therefore, a subset $KZ$ of $\widetilde{K}_0$ is compact, hence closed.
 We prove that $1\notin KZ$. If $1\in KZ$, there exists $k\in K$ such that
 $k(i_*P_1)k^{-1}=i_*P_2$. Taking intersections with $G$, we obtain
 $kP_1k^{-1}=P_2$, which contradicts the assumption $[P_1]\neq [P_2]$.
 Therefore, we can find a compact open normal subgroup $\widetilde{K}_1$ of $\widetilde{K}_0$
 such that $\widetilde{K}_1\cap KZ=\varnothing$.
 Then, $\widetilde{K}_{[P_1],[P_2]}=K\widetilde{K}_1$ is a compact open subgroup of $\GSp_{2n}(\A^\infty)$
 satisfying $\widetilde{K}_{[P_1],[P_2]}\cap Z=\varnothing$.
 This concludes the proof.
\end{prf}

\begin{prf}[of Theorem \ref{thm:Sh-partition} for the Hodge type case]
 We assume that $(G,X)$ is of Hodge type, and take an embedding $i\colon (G,X)\hookrightarrow (\GSp_{2n},X_{2n})$
 into a Siegel Shimura datum.
 Further, we take a compact open subgroup $\widetilde{K}\subset \GSp_{2n} (\mathbb{A}^{\infty})$
 as in Lemma \ref{lem:P_G(K)-Hodge}.

 Let $[P]\in \mathcal{P}_G(K)$.
 Since Theorem \ref{thm:Sh-partition} is known for the Siegel case,
 we have a locally closed constructible subset $C_{[i_*P]}$ of
 $\Sh_{\widetilde{K}}(\GSp_{2n},X_{2n})_{\Q_p}^{\ad}$. 
 Let $C_{[P]}$ be the inverse image of $C_{[i_*P]}$ in 
 $\Sh_K(G,X)_{E_v}^{\ad}$. Then, $C_{[P]}$ satisfies the desired condition
 by Proposition \ref{prop:Phi-functoriality} \ref{item:Phi-Hodge}. 
 The subset $C_{[G]}$ is open and quasi-compact,
 since $C_{[\GSp_{2n}]}=\mathcal{S}_{\widetilde{K},\mathopen{]}\le 0\mathclose{[}}$ is open and quasi-compact.
\end{prf}

The following lemma is the Hodge type version of Corollary \ref{cor:C_P-disjoint-Siegel}.

\begin{lem}\label{lem:C_P-disjoint-Hodge}
 Let $(G,X)$ be a Shimura datum of Hodge type. 
 Take an embedding $i\colon (G,X)\hookrightarrow (\GSp_{2n},X_{2n})$ into a Siegel Shimura datum. 
 For a compact open subgroup $K$ of $G(\A^\infty)$ and an integer $0\le r\le n$,
 we write $\mathcal{P}_G(K)_r$ for the inverse image of $r$ under the composite
 \[
  \mathcal{P}_G(K)\to \mathcal{P}_{G,\Q}\xrightarrow{i_*}\mathcal{P}_{\GSp_{2n},\Q}\cong \{0,1,\ldots,n\}.
 \]
 For $[P]\in \mathcal{P}_G(K)_r$ and $[P']\in \mathcal{P}_G(K)_{r'}$, assume that
 $r>r'$ or $r=r'$ and $[P]\neq [P']$.
 Then we have $C_{[P]}^-\cap C_{[P']}=\varnothing$,
 where $C_{[P]}^-$ denotes the closure of $C_{[P]}$ in $\Sh_K(G,X)_{E_v}^\ad$.

 In particular, if $[P]$, $[P']$ are distinct elements of $\mathcal{P}_G(K)$ such that
 $P$ is $G(\A^\infty)$-conjugate to $P'$, then we have $C_{[P]}^-\cap C_{[P']}=\varnothing$.
\end{lem}

\begin{prf}
 We take $\widetilde{K}$ as in Lemma \ref{lem:P_G(K)-Hodge}.
 By definition, we have $[i_*P]\in \mathcal{P}_{\GSp_{2n}}(\widetilde{K})_r$ and 
 $[i_*P']\in \mathcal{P}_{\GSp_{2n}}(\widetilde{K})_{r'}$.
 Since the map $\mathcal{P}_G(K)\to \mathcal{P}_{\GSp_{2n}}(\widetilde{K})$ is injective,
 we have $[i_*P]\neq [i_*P']$ if $r=r'$.
 Hence Corollary \ref{cor:C_P-disjoint-Siegel} tells us that $C^-_{[i_*P]}\cap C_{[i_*P']}=\varnothing$.

 Since the natural morphism $\Sh_{K}(G,X)^\ad_{E_v}\to\Sh_{\widetilde{K}}(\GSp_{2n},X_{2n})^\ad_{\Q_p}$
 maps $C_{[P]}$ and $C_{[P']}$ into $C_{[i_*P]}$ and $C_{[i_*P']}$, respectively
 (see the construction of $C_{[P]}$ in the proof of Theorem \ref{thm:Sh-partition}),
 the set $C_{[P]}^-\cap C_{[P']}$ is mapped into $C_{[i_*P]}^-\cap C_{[i_*P']}=\varnothing$.
 Therefore we conclude that $C^-_{[P]}\cap C_{[P']}=\varnothing$.

 The last claim follows from the observation that if $P$ and $P'$ are $G(\A^\infty)$-conjugate, then $r=r'$.
\end{prf}

Now we can prove Theorem \ref{thm:Sh-partition} for the preabelian type case. 

\begin{prf}[of Theorem \ref{thm:Sh-partition}]
 We choose a compact open subgroup $K''$ of $G^\ad(\A^\infty)$ in such a way as in
 Proposition \ref{lem:P_G(K)} \ref{item:P_G(K)-adjoint}, and use the notation in Lemma \ref{lem:preHod}.
 We write $E'$ for the reflex field of $(G',X')$, and choose a place $v'$ of $E'$ above $v^\ad$
 (for the definition of $v^\ad$, see Proposition \ref{prop:Phi-functoriality} \ref{item:Phi-ad}).

 We have natural maps
 \[
 \mathcal{P}_{G'}(K')\to\mathcal{P}_{G'^\ad}(g_i^{-1}K''g_i)\xrightarrow{g_i}\mathcal{P}_{G'^\ad}(K'')\cong \mathcal{P}_{G^\ad}(K'')\cong \mathcal{P}_G(K)
 \]
 for each $i$ (see Lemma \ref{lem:P_G(K)} \ref{item:P_G(K)-adjoint}).
 Let $S_i$ be the inverse image of $[P]\in \mathcal{P}_G(K)$ under this map.
 We put $C^{\ad}_{[P]}=\bigcup_{1\le i\le m}f_i(\bigcup_{[Q]\in S_i}C_{[Q]})$.
 Since $C_{[Q]}$ for each $[Q]\in S_i$ is constructible,
 Lemma \ref{lem:etale-constructible-image} \ref{item:constr-image} tells us that $C^{\ad}_{[P]}$ is constructible.
 Let us prove that $C^{\ad}_{[P]}$ is locally closed.
 By Lemma \ref{lem:cl-pt-constructible} (i), 
 Proposition \ref{prop:Phi-functoriality} \ref{item:Phi-Hecke}, \ref{item:Phi-ad} and the constructibility
 of $C^{\ad}_{[P]}$, we can check that the inverse image of $C^{\ad}_{[P]}$ under $\coprod_{1\le i\le m}f_i$
 equals $\coprod_{1\le i\le m}\bigcup_{[Q]\in S_i}C_{[Q]}$. 
 Therefore, Lemma \ref{lem:etale-constructible-image} \ref{item:loc-cl-image} tells us that
 it suffices to prove that $\bigcup_{[Q]\in S_i}C_{[Q]}$ is locally closed. 
 Note that we have the commutative diagram
 \[
 \xymatrix{% 
 \mathcal{P}_{G'}(K')\ar[r]\ar[d]&\mathcal{P}_{G'^\ad}(g_i^{-1}K''g_i)\ar[r]^-{g_i}\ar[d]&\mathcal{P}_{G'^\ad}(K'')\ar@{<->}[r]^{\cong}\ar[d]&\mathcal{P}_{G^\ad}(K'')\ar[d]&\mathcal{P}_G(K)\ar[l]_-{\cong}\ar[d]\\
 \mathcal{P}_{G',\Q}\ar[r]^{\cong}&\mathcal{P}_{G'^\ad,\Q}\ar@{=}[r]&\mathcal{P}_{G'^\ad,\Q}\ar@{<->}[r]^-{\cong}&\mathcal{P}_{G^\ad,\Q}&\mathcal{P}_{G,\Q}\ar[l]_-{\cong}\lefteqn{,}
 }
 \]
 where the vertical arrows are the maps in Corollary \ref{cor:adm-parab-conj}.
 From this diagram we can observe that all elements in $S_i$ are $G'(\A^\infty)$-conjugate.
 Now Lemma \ref{lem:C_P-disjoint-Hodge} tells us that the closure $C_{[Q]}^-$ of $C_{[Q]}$ for $[Q]\in S_i$
 does not intersect  $\bigcup_{[Q']\in S_i\setminus \{[Q]\}}C_{[Q']}$. 
 Therefore $C_{[Q]}$ is closed (hence open) in $\bigcup_{[Q']\in S_i}C_{[Q']}$,
 from which we conclude that $\bigcup_{[Q]\in S_i}C_{[Q]}$
 is locally closed in $\Sh_{K'}(G',X')^\ad_{E'_{v'}}$, as desired.

 Let $C_{[P]}$ be the inverse image of $C^{\ad}_{[P]}$ under the map
 \[
  \Sh_K(G,X)^\ad_{E_v}\to \Sh_{K''}(G^\ad,X^\ad)^\ad_{E^\ad_{v^\ad}}\cong \Sh_{K''}(G'^\ad,X'^\ad)^\ad_{E^\ad_{v^\ad}}
 \]
 (here $E^\ad$ and $v^\ad$ are as in Proposition \ref{prop:Phi-functoriality} \ref{item:Phi-ad}).
 Then $C_{[P]}$ satisfies the desired condition by Proposition \ref{prop:Phi-functoriality} \ref{item:Phi-Hecke}, \ref{item:Phi-ad}.
 
 If $[P]=[G]$, we have $S_i=\{[G']\}$. Hence the openness and the quasi-compactness of $C_{[G]}$ follows from
 those of $C_{[G']}$.
\end{prf}

\begin{lem}\label{lem:C_P-disjoint}
 Let $(G,X)$ be a Shimura datum of preabelian type. Take a Shimura datum $(G',X')$ of Hodge type
 such that $(G^\ad,X^\ad)\cong(G'^\ad,X'^\ad)$, and an embedding $i\colon (G',X')\hookrightarrow (\GSp_{2n},X_{2n})$ into
 a Siegel Shimura datum. For a compact open subgroup $K$ of $G(\A^\infty)$ and an integer $0\le r\le n$,
 we write $\mathcal{P}_G(K)_r$ for the inverse image of $r$ under the composite
 \[
  \mathcal{P}_G(K)\to \mathcal{P}_{G,\Q}\cong\mathcal{P}_{G^\ad,\Q}\cong\mathcal{P}_{G'^\ad,\Q}\cong \mathcal{P}_{G',\Q}\xrightarrow{i_*}\mathcal{P}_{\GSp_{2n},\Q}\cong \{0,1,\ldots,n\}.
 \]
 For $[P]\in \mathcal{P}_G(K)_r$ and $[P']\in \mathcal{P}_G(K)_{r'}$, assume that
 $r>r'$ or $r=r'$ and $[P]\neq [P']$.
 Then we have $C_{[P]}^-\cap C_{[P']}=\varnothing$,
 where $C_{[P]}^-$ denotes the closure of $C_{[P]}$ in $\Sh_K(G,X)_{E_v}^\ad$.

 In particular, if $[P]$, $[P']$ are distinct elements of $\mathcal{P}_G(K)$ such that
 $P$ is $G(\A^\infty)$-conjugate to $P'$, then we have $C_{[P]}^-\cap C_{[P']}=\varnothing$.
\end{lem}

\begin{prf}
 We use the same notation as in the proof of Theorem \ref{thm:Sh-partition} above,
 and denote the morphism $\Sh_K(G,X)^\ad_{E_v}\to \Sh_{K''}(G^\ad,X^\ad)^\ad_{E^\ad_{v^\ad}}\cong \Sh_{K''}(G'^\ad,X'^\ad)^\ad_{E^\ad_{v^\ad}}$ by $h$.
 If a point $x$ belongs to $C^-_{[P]}\cap C_{[P']}$, 
 we can find $y\in C_{[P]}$ specializing to $x$ (see \cite[Corollary of Theorem 1]{MR0251026}).
 By the construction of $C_{[P]}$, there exist $1\le i\le m$, $[Q]\in S_i$ and $y'\in C_{[Q]}$
 such that $h(y)=f_i(y')$. Since $f_i$ is finite, there exists $x'\in \Sh_{K'}(G',X')^\ad_{E'_{v'}}$
 which is a specialization of $y'$ and mapped to $h(x)$ by $f_i$.
 By Remark \ref{rem:C_P-cover}, $x'$ belongs to $C_{[Q'']}$ for a unique $[Q'']\in\mathcal{P}_{K'}(G')$.
 Take $r''$ such that $[Q'']$ lies in $\mathcal{P}_{K'}(G')_{r''}$. 
 Since $x'\in C^-_{[Q]}\cap C_{[Q'']}$, Lemma \ref{lem:C_P-disjoint-Hodge} tells us that
 we have either $r<r''$ or $[Q]=[Q'']$.

 We write $[P'']$ for the image of $[Q'']$ under the composite
 \[
 \mathcal{P}_{G'}(K')\to\mathcal{P}_{G'^\ad}(g_i^{-1}K''g_i)\xrightarrow{g_i}\mathcal{P}_{G'^\ad}(K'')\cong \mathcal{P}_{G^\ad}(K'')\cong \mathcal{P}_G(K).
 \]
 It belongs to $\mathcal{P}_G(K)_{r''}$. Since $h(x)=f_i(x')\in C^\ad_{[P'']}$, the point $x$ lies in $C_{[P'']}$.
 Hence Remark \ref{rem:C_P-cover} tells us that $[P']=[P'']$. In particular we have $r'=r''$,
 which implies $[Q]=[Q'']$. Therefore we have $[P]=[P'']=[P']$, which contradicts the assumption on $[P']$.
\end{prf}

\begin{cor}\label{cor:Hecke-inverse-disjoint}
 Let $(G,X)$ be a Shimura datum of preabelian type.
 Let $K$ and $K'$ be compact open subgroups of $G(\A^\infty)$ and
 $g\in G(\A^\infty)$ with $g^{-1}Kg\subset K'$.
 For an element $[P']$ of $\mathcal{P}_G(K')$, the inverse image of $C_{[P']}$ under
 the Hecke action $g\colon \Sh_K(G,X)\to \Sh_{K'}(G,X)$ is equal to
 \[
  \coprod_{\substack{[P]\in\mathcal{P}_G(K),\\ [g^{-1}Pg]=[P']\text{ \upshape in $\mathcal{P}_G(K')$}}} C_{[P]}
 \]
 as topological spaces.

 In particular, for $[P]\in \mathcal{P}_K(G)$, $C_{[P]}$ is mapped to $C_{[g^{-1}Pg]}$ under
 the Hecke action by $g$.
\end{cor}

\begin{prf}
 First note that both $g^{-1}(C_{[P']})$ and $\bigcup_{[P]\in\mathcal{P}_K(G),[g^{-1}Pg]=[P']}C_{[P]}$
 are constructible subsets and have the same set of classical points by
 Proposition \ref{prop:Phi-functoriality} \ref{item:Phi-Hecke}.
 Therefore, by Lemma \ref{lem:cl-pt-constructible} (i), they are equal.
 The union $\bigcup_{[P]\in\mathcal{P}_K(G),[g^{-1}Pg]=[P']}C_{[P]}$ is set-theoretically disjoint
 by Remark \ref{rem:C_P-cover}.
 Let $[P_1]$, $[P_2]$ be distinct elements of $\mathcal{P}_G(K)$ such that $[g^{-1}P_1g]=[g^{-1}P_2g]$
 in $\mathcal{P}_{K'}(G)$. 
 Lemma \ref{lem:C_P-disjoint} tells us that $C^-_{[P_1]}\cap C_{[P_2]}=\varnothing$. 
 This implies that $C_{[P_1]}$ is closed (hence open) in $\bigcup_{[P]\in\mathcal{P}_K(G),[g^{-1}Pg]=[P']}C_{[P]}$.
 Now the proof is complete.
\end{prf}

The following lemma will be used in the next section.

\begin{lem}\label{lem:U_r}
 Let the notation be as in Lemma \ref{lem:C_P-disjoint}.
 Put $\mathcal{P}_G(K)_{\le r}=\bigcup_{r'\le r}\mathcal{P}_G(K)_{r'}$.
 Then, $U_{r,K}=\bigcup_{[P]\in \mathcal{P}_G(K)_{\le r}}C_{[P]}$ is a constructible open subset of
 $\Sh_K(G,X)^\ad_{E_v}$, and $U_{r,K}\setminus U_{r-1,K}$ equals $\coprod_{[P]\in\mathcal{P}_G(K)_r}C_{[P]}$
 as topological spaces (here we put $U_{-1,K}=\varnothing$).
\end{lem}
 
\begin{prf}
 By Lemma \ref{lem:C_P-disjoint}, the set $\bigcup_{r'>r}\bigcup_{[P]\in \mathcal{P}_G(K)_{r'}}C_{[P]}$
 is closed in $\Sh_K(G,X)^\ad_{E_v}$. Therefore, by Remark \ref{rem:C_P-cover}, $U_{r,K}$ is open
 in $\Sh_K(G,X)^\ad_{E_v}$. Since $C_{[P]}$ is constructible for every $[P]$,
 the subset $U_{r,K}$ is also constructible.

 The claim $U_{r,K}\setminus U_{r-1,K}=\coprod_{[P]\in\mathcal{P}_G(K)_r}C_{[P]}$ follows from
 Lemma \ref{lem:C_P-disjoint} by the same argument as in the proof of
 Corollary \ref{cor:Hecke-inverse-disjoint}.
\end{prf}

\begin{rem}\label{rem:U_r-Hodge}
 If $(G,X)$ is of Hodge type (namely, $(G,X)=(G',X')$), we can also construct $U_{r,K}$ in the following way.
 Take a compact open subgroup $\widetilde{K}\subset \GSp_{2n}(\mathbb{A}^{\infty})$
 containing $K$.
 Then, $U_{r,K}$ equals the inverse image of $\mathcal{S}_{\widetilde{K},\mathopen{]}\le r\mathclose{[}}$
 under 
 \[
  \Sh_{K}(G,X)^\ad_{E_v}\to \Sh_{\widetilde{K}}(\GSp_{2n},X_{2n})^\ad_{\Q_p}=\mathcal{S}_{\widetilde{K}}.
 \]
 This is an immediate consequence of Proposition \ref{prop:Phi-functoriality} \ref{item:Phi-Hodge}
 and Proposition \ref{prop:Siegel-partition-r}.
\end{rem}

\section{Cohomology of Shimura varieties}\label{sec:cohomology}
\subsection{Comparison of cohomology}\label{subsec:main-thm}
We continue to assume that $(G,X)$ is of preabelian type.
For simplicity, we further assume that $(G,X)$ satisfies SV6 in \cite[p.~311]{MR2192012}.
We simply write $\Sh_K$ for $\Sh_K(G,X)$, if there is no risk of confusion. 
Fix a prime $\ell$ which is different from $p$.
Let $G^c$ be the quotient of $G$ defined in \cite[p.~347]{MR1044823},
and $\xi$ an algebraic representation of $G^c$ on a finite-dimensional $\overline{\Q}_\ell$-vector space. 
Then we have the associated $\overline{\Q}_\ell$-sheaf $\mathcal{L}_\xi$ on $\Sh_{K}$
(see \cite[III, \S 6]{MR1044823}).
Moreover, $\mathcal{L_\xi}$ is equivariant with respect to the Hecke action.

Let $p'$ be a prime number. Let us fix a compact open subgroup $K^{p'}$ of $G(\widehat{\Z}^{p'})$. 
We consider the compactly supported cohomology
\[
 H^i_c(\Sh_{\infty,K^{p'},\overline{E}_v},\mathcal{L}_\xi)
 =\varinjlim_{K_{p'}}
 H^i_c(\Sh_{K_{p'}K^{p'},\overline{E}_v},
 \mathcal{L}_\xi).
\]
The group $G(\Q_{p'})\times \Gal(\overline{E}_v/E_v)$ naturally acts on it.
By this action, $H^i_c(\Sh_{\infty,K^{p'},\overline{E}_v},\mathcal{L}_\xi)$ becomes an admissible/continuous representation of
$G(\Q_{p'})\times \Gal(\overline{E}_v/E_v)$ 
in the sense of \cite[\S I.2]{MR1876802}. 

The group $G(\Q_{p'})\times \Gal(\overline{E}_v/E_v)$ 
naturally acts also on 
\[
 H^i_c(\Sh^{\pg}_{\infty,K^{p'},\overline{E}_v},\mathcal{L}_\xi^{\ad})
 =\varinjlim_{K_{p'}}
 H^i_c(\Sh_{K_{p'}K^{p'}, \overline{E}_v}^{\pg}, 
 \mathcal{L}_\xi^{\ad}). 
\]
See \cite[\S 1]{MR1626021} for the definition of the compactly supported $\ell$-adic cohomology
for adic spaces. 
It gives an admissible/continuous representation of $G(\Q_{p'})\times \Gal(\overline{E}_v/E_v)$
(\cf Lemma \ref{lem:coh-invariant-part} in the next subsection).

Here we use the notation in \cite[\S I.2]{MR1876802}.
Let $H$ be a locally profinite group.
For an admissible/continuous representation $V$ of $H\times \Gal(\overline{E}_v/E_v)$
over $\overline{\Q}_\ell$ and an irreducible admissible representation $\pi$ of $H$, 
put $V[\pi]=\bigoplus_{\sigma}\sigma^{\oplus m_{\pi\boxtimes\sigma}}$,
where $\sigma$ runs through finite-dimensional irreducible continuous $\overline{\Q}_\ell$-representations
of $\Gal(\overline{E}_v/E_v)$ and $m_{\pi\boxtimes\sigma}$ denotes the coefficient of $[\pi\boxtimes\sigma]$
in the image of $V$ in the Grothendieck group considered in \cite[\S I.2]{MR1876802}.
It is a semisimple continuous representation of $\Gal(\overline{E}_v/E_v)$. 

\begin{thm}\label{thm:main-thm}
 The kernel and the cokernel of the canonical homomorphism
 \[
 H^i_c(\Sh^{\pg}_{\infty,K^{p'},\overline{E}_v},\mathcal{L}_\xi^{\ad})
 \to 
 H^i_c(\Sh_{\infty,K^{p'},\overline{E}_v},\mathcal{L}_\xi)
 \]
 are non-cuspidal, namely, they have no supercuspidal subquotient of 
 $G(\Q_{p'})$.
 In particular, for an irreducible supercuspidal representation $\pi$ 
 of $G(\Q_{p'})$,
 we have an isomorphism 
 \[
  H^i_c(\Sh^{\pg}_{\infty,K^{p'},\overline{E}_v},\mathcal{L}_\xi^{\ad})[\pi]\cong 
 H^i_c(\Sh_{\infty,K^{p'},\overline{E}_v},\mathcal{L}_\xi)[\pi].
 \]
\end{thm}

This theorem will be proved in Section \ref{subsec:proof-main-thm}.

\begin{rem}\label{rem:cpt-ordinary}
 Let $H^i(\Sh_{\infty,K^{p'},\overline{E}_v},\mathcal{L}_\xi)=\varinjlim_{K_{p'}}H^i(\Sh_{K_{p'}K^{p'},\overline{E}_v},\mathcal{L}_\xi)$ 
 be the ordinary cohomology of our Shimura variety.
 This is also an admissible/continuous representation of $G(\Q_{p'})\times \Gal(\overline{E}_v/E_v)$.
 By using the minimal compactification of $\Sh_{K}$ and its natural stratification (\cf \cite[\S 3.7]{MR1149032}),
 it is easy to see that the kernel and the cokernel of
 the canonical homomorphism
 \[
 H^i_c(\Sh_{\infty,K^{p'},\overline{E}_v},\mathcal{L}_\xi)
 \to 
 H^i(\Sh_{\infty,K^{p'},\overline{E}_v},\mathcal{L}_\xi)
 \]
 are non-cuspidal as $G(\Q_{p'})$-representations (in fact, we can use the similar argument
 as in the next subsection). 
 Therefore, the kernel and the cokernel of the composite
 \[
  H^i_c(\Sh^{\pg}_{\infty,K^{p'},\overline{E}_v},\mathcal{L}_\xi^{\ad})
  \to 
  H^i_c(\Sh_{\infty,K^{p'},\overline{E}_v},\mathcal{L}_\xi)
  \to 
  H^i(\Sh_{\infty,K^{p'},\overline{E}_v},\mathcal{L}_\xi)
 \]
 are non-cuspidal by Theorem \ref{thm:main-thm}. 
\end{rem}

\begin{rem}\label{rem:cpt-IH}
 Let $\mathit{IH}^i(\Sh_{\infty,K^{p'},\overline{E}_v},\mathcal{L}_\xi)
 =\varinjlim_{K_{p'}}H^i(\Sh^{\mathrm{min}}_{K_{p'}K^{p'},\overline{E}_v},j_{!*}\mathcal{L}_\xi)$ be the intersection cohomology of our Shimura variety,
 where $j\colon \Sh_{K_{p'}K^{p'}}\hookrightarrow \Sh_{K_{p'}K^{p'}}^{\mathrm{min}}$ denotes the minimal
 compactification of $\Sh_{K_{p'}K^{p'}}$.
 Then, as in the previous remark, it is easy to see that the kernel and the cokernel of the canonical homomorphism
 \[
  H^i_c(\Sh_{\infty,K^{p'},\overline{E}_v},\mathcal{L}_\xi)
  \to
  \mathit{IH}^i(\Sh_{\infty,K^{p'},\overline{E}_v},\mathcal{L}_\xi)
 \]
 are non-cuspidal. Therefore, by Theorem \ref{thm:main-thm},
 we have an isomorphism
 \[
  H^i_c(\Sh^{\pg}_{\infty,K^{p'},\overline{E}_v},\mathcal{L}_\xi^{\ad})[\pi]\cong 
  \mathit{IH}^i(\Sh_{\infty,K^{p'},\overline{E}_v},\mathcal{L}_\xi)[\pi]
 \]
 for an irreducible supercuspidal representation $\pi$ of $G(\Q_{p'})$. 
\end{rem}

\begin{cor}\label{cor:Piupsc}
 We put
\[
 H^i_c(\Sh_{\infty,\overline{E}_v},\mathcal{L}_\xi)
 =\varinjlim_{K^{p'}}
 H^i_c(\Sh_{\infty,K^{p'},\overline{E}_v},\mathcal{L}_\xi),\
 H^i_c(\Sh^{\pg}_{\infty,\overline{E}_v},\mathcal{L}_\xi^{\ad})=
 \varinjlim_{K^{p'}} 
 H^i_c(\Sh^{\pg}_{\infty,K^{p'},\overline{E}_v},\mathcal{L}_\xi^{\ad}).
\]
 These are admissible/continuous $G(\A^\infty)\times \Gal(\overline{E}_v/E_v)$-representations.

 Let $\Pi$ be an irreducible admissible representation of $G(\A^\infty)$.
 Assume that there exists a prime $p'$ such that $\Pi_{p'}$ is a supercuspidal representation of $G(\Q_{p'})$.
 Then, $\Pi$ does not appear as a subquotient of the kernel or the cokernel of the canonical
 homomorphism
 $H^i_c(\Sh^{\pg}_{\infty,\overline{E}_v},\mathcal{L}_\xi^{\ad})\to H^i_c(\Sh_{\infty,\overline{E}_v},\mathcal{L}_\xi)$. 
 In particular, we have an isomorphism of $\Gal(\overline{E}_v/E_v)$-representations
 \[
 H^i_c(\Sh^{\pg}_{\infty,\overline{E}_v},\mathcal{L}_\xi^{\ad})[\Pi]\cong 
 H^i_c(\Sh_{\infty,\overline{E}_v},\mathcal{L}_\xi)[\Pi].
 \]
\end{cor}

\begin{prf}
 We take a compact open subgroup $K^{p'}\subset G(\widehat{\Z}^{p'})$ such that
 $\Pi^{K^{p'}}\neq 0$. If $\Pi$ appears as a subquotient of the kernel or the cokernel of 
 $H^i_c(\Sh^{\pg}_{\infty,\overline{E}_v},\mathcal{L}_\xi^{\ad})
 \to 
 H^i_c(\Sh_{\infty,\overline{E}_v},\mathcal{L}_\xi)$, 
 then $\Pi_{p'}$ appears as a subquotient of the kernel or the cokernel of
 $H^i_c(\Sh^{\pg}_{\infty,K^{p'},\overline{E}_v},\mathcal{L}_\xi^{\ad})\to H^i_c(\Sh_{\infty,K^{p'},\overline{E}_v},\mathcal{L}_\xi)$ (\cf Lemma \ref{lem:coh-invariant-part} in the next subsection).
 This contradicts Theorem \ref{thm:main-thm}. 
\end{prf}

\subsection{Proof of Theorem \ref{thm:main-thm}}\label{subsec:proof-main-thm}
Let $K$ be a compact open subgroup of $G(\A^\infty)$.
We regard $C_{[P]}$ for $[P]\in\mathcal{P}_G(K)$ 
as a pseudo-adic space (\cf \cite[\S 1.10]{MR1734903}). 
See \cite[Proposition 2.6 (i)]{MR1626021} for the definition of the compactly supported $\ell$-adic cohomology of
pseudo-adic spaces. 

\begin{prop}\label{prop:finiteness-cohomology}
 Let $[P]$ be an element of $\mathcal{P}_G(K)$.
 \begin{enumerate}
  \item\label{item:coh-fin-l-adic} 
       For a constructible $\ell$-adic sheaf $\mathcal{F}=(\mathcal{F}_n)\otimes\overline{\Q}_\ell$
       on $\Sh_{K,\overline{E}_v}$, $H^i_c(C_{[P],\overline{E}_v},\mathcal{F}^\ad)$ is
       a finite-dimensional $\overline{\Q}_\ell$-vector space.
  \item\label{item:coh-fin-torsion}
       For a constructible $\Z/\ell^n\Z$-sheaf $\mathcal{F}$ on $\Sh_{K,\overline{E}_v}$,
       $H^i_c(C_{[P],\overline{E}_v},\mathcal{F}^\ad)$ is a finitely generated $\Z/\ell^n\Z$-module.
 \end{enumerate}
\end{prop}

Before proving this proposition, we note the following general lemmas.

\begin{lem}\label{lem:coh-invariant-part}
 Let $k$ be an algebraically closed non-archimedean field,
 and $X$ an adic space locally of finite type, separated and taut over $k$.
 Let $L$ be a locally closed constructible subset of $X$, which is regarded as a pseudo-adic space. 
 Let $\pi\colon X'\to X$ be a finite \'etale Galois covering with Galois group $H$.
 We put $L'=\pi^{-1}(L)$.
 For an $\ell$-adic sheaf $\mathcal{F}=(\mathcal{F}_n)\otimes\overline{\Q}_\ell$ on $X$, the natural map
 \[
  H^i_c(L,\mathcal{F})\to H^i_c(L',\pi^*\mathcal{F})^H
 \]
 is an isomorphism.
\end{lem} 

\begin{prf}
 This lemma might be well-known, but we include its proof for reader's convenience.
 We put $\mathcal{F}'_n=\pi_*\pi^*\mathcal{F}_n$ and $\mathcal{F}'=(\mathcal{F}'_n)\otimes\overline{\Q}_\ell$.
 The group $H$ acts on $\mathcal{F}'_n$, and we have $(\mathcal{F}'_n)^H=\mathcal{F}_n$.
 Consider the map $\psi=\sum_{h\in H}h\colon \mathcal{F}'_n\to \mathcal{F}_n$.
 The composite $\mathcal{F}_n\hookrightarrow\mathcal{F}'_n\xrightarrow{\psi} \mathcal{F}_n$ equals 
 the multiplication by $\#H$.
 By taking the cohomology, we have a commutative diagram
 \[
 \xymatrix{%
 H^i_c\bigl(L,(\mathcal{F}_n)_n\bigr)\otimes\overline{\Q}_\ell\ar[r]\ar[rd]^-{\pi^*}& H^i_c\bigl(L,(\mathcal{F}'_n)_n\bigr)\otimes\overline{\Q}_\ell\ar[r]^-{H^i_c(\psi)}\ar@{=}[d]& H^i_c\bigl(L,(\mathcal{F}_n)_n\bigr)\otimes\overline{\Q}_\ell\ar[d]^-{\pi^*}\\
 & H^i_c\bigl(L',(\pi^*\mathcal{F}_n)_n\bigr)\otimes\overline{\Q}_\ell\ar[r]^-{\sum_{h\in H}h}& H^i_c\bigl(L',(\pi^*\mathcal{F}_n)_n\bigr)\otimes\overline{\Q}_\ell\lefteqn{.}
 }
 \]
 The composite of two upper horizontal arrows is the multiplication by $\#H$, which is an isomorphism.
 Therefore $\pi^*$ is injective and $H^i_c(\psi)$ is surjective. 
 The surjectivity of $H^i_c(\psi)$ implies that the image of $\pi^*$ is equal to
 that of $\sum_{h\in H}h$, that is, the $H$-invariant part of
 $H^i_c(L',(\pi^*\mathcal{F}_n)_n)\otimes\overline{\Q}_\ell$.
\end{prf}

\begin{rem}\label{rem:limit-invariant}
 By the same method, we can also prove that the natural map
 \[
 \varprojlim_{n}H^i_c(L,\mathcal{F}_n)\otimes\overline{\Q}_\ell\to \Bigl(\varprojlim_{n}H^i_c(L',\pi^*\mathcal{F}_n)\otimes\overline{\Q}_\ell\Bigr)^H
 \]
 is an isomorphism.
\end{rem}

\begin{lem}\label{lem:excision-exact}
 Let $k$ be an algebraically closed non-archimedean field.
 Let $X$ be an adic space locally of finite type, separated and taut over $k$,
 and $U$ a constructible open subset of $X$. Set $Z=X\setminus U$, 
 which is regarded as a pseudo-adic space. 
 For an $\ell$-adic sheaf $\mathcal{F}=(\mathcal{F}_n)\otimes\overline{\Q}_\ell$ over $X$,
 we have a long exact sequence
 \[
 \cdots\to H_c^i(U,\mathcal{F})\to H_c^i(X,\mathcal{F})\to H_c^i(Z,\mathcal{F})\to H_c^{i+1}(U,\mathcal{F})\to \cdots.
 \]
\end{lem}

\begin{prf}
 Since the open immersion $j\colon U\hookrightarrow X$ is quasi-compact,
 we have $H^i_c(U,\mathcal{F})=H^i_c(X,j_!j^*\mathcal{F})$ by the definition of the compactly supported cohomology.
 Therefore, the claim follows from \cite[Proposition 2.6 (i)]{MR1626021}.
\end{prf}

\begin{prf}[of Proposition \ref{prop:finiteness-cohomology}]
 We consider \ref{item:coh-fin-l-adic}. 
 Note that Corollary \ref{cor:Hecke-inverse-disjoint} and Lemma \ref{lem:coh-invariant-part}
 enable us to shrink $K$ arbitrarily. 
 
 First we consider the Hodge type case. Take an embedding $(G,X)\hookrightarrow (\GSp_{2n},X_{2n})$ into
 a Siegel Shimura datum.
 We have a constructible open subset $U_{r,K}$ for each $0\le r\le n$ by Lemma \ref{lem:U_r}.
 By the long exact sequence
 \begin{align*}
  \cdots&\to H^i_c(U_{r-1,K,\overline{E}_v},\mathcal{F})\to H^i_c(U_{r,K,\overline{E}_v},\mathcal{F})\to H^i_c(U_{r,K,\overline{E}_v}\setminus U_{r-1,K,\overline{E}_v},\mathcal{F})\\
  &\to H^{i+1}_c(U_{r-1,K,\overline{E}_v},\mathcal{F})\to \cdots
 \end{align*}
 (see Lemma \ref{lem:excision-exact}) and Lemma \ref{lem:U_r},
 it suffices to show that $H^i_c(U_{r,K,\overline{E}_v},\mathcal{F})$ is finite-dimensional
 for each $0\le r\le n$.

 Take a compact open subgroup $\widetilde{K}$ of $G(\A^\infty)$ so that
 we have a natural embedding $\Sh_K(G,X)\hookrightarrow \Sh_{\widetilde{K}}(\GSp_{2n},X_{2n})_E$.
 By shrinking $K$, we may assume that $\widetilde{K}=\widetilde{K}_p\widetilde{K}^p$ with
 $\widetilde{K}_p\subset \GSp_{2n}(\Z_p)$. 
 Let $\mathscr{S}_{\widetilde{K}}^{\mathrm{nor}}$ be the normalization of 
 $\mathscr{S}_{\widetilde{K}^p}^{\mathrm{tor}}$ in $\Sh_{\widetilde{K}}(\GSp_{2n},X_{2n})_{\mathbb{Q}_p}$, 
 and $\mathscr{S}^\mathrm{nor}_{\widetilde{K},\le r,\F_p}$ the inverse image of 
 $\mathscr{S}^\mathrm{tor}_{\widetilde{K}^p,\le r,\mathbb{F}_p}$ in $\mathscr{S}_{\widetilde{K}}^{\mathrm{nor}}$. 
 We write $Y$ for the closure of $\Sh_K(G,X)_{E_v}$ in $\mathscr{S}^{\mathrm{nor}}_{\widetilde{K},E_v}$ and
 put $Z=Y\setminus \Sh_K(G,X)_{E_v}$.
 Let $V$ denote the inverse image of $\mathscr{S}^\mathrm{nor}_{\widetilde{K},\le r,\F_p}$ under the composite
 \[
 Y^\ad\hookrightarrow (\mathscr{S}^{\mathrm{nor}}_{\widetilde{K}})_{E_v}^\ad\to (\mathscr{S}^{\mathrm{nor}}_{\widetilde{K}})_{\Q_p}^\ad=(\mathscr{S}^{\mathrm{nor}}_{\widetilde{K}})^{\wedge\rig}\xrightarrow{\spp}\mathscr{S}^{\mathrm{nor}}_{\widetilde{K},\F_p},
 \]
 which is a quasi-compact open subset of $Y^\ad$.
 Note that $V\setminus (V\cap Z^\ad)=U_{r,K}$ by Remark \ref{rem:U_r-Hodge}.
 Therefore, \cite[Theorem 3.3 (i)]{MR1626021} tells us that $H^i_c(U_{r,K,\overline{E}_v},\mathcal{F})$ is
 finite-dimensional. This completes the proof in the Hodge type case.

 Next we consider the preabelian type case. 
 We choose a compact open subgroup of $K''$ of $G^\ad(\A^\infty)$ in such a way as in
 Proposition \ref{lem:P_G(K)} \ref{item:P_G(K)-adjoint}, and use the notation in Lemma \ref{lem:preHod}.
 Then, for $[P]\in \mathcal{P}_K(G)$, the inverse image of $C_{[P^\ad]}$ under
 $\pi\colon \Sh_K(G,X)\to \Sh_{K''}(G^\ad,X^\ad)$ is equal to $C_{[P]}$.
 Therefore, by pushing forward sheaves by $\pi$, we may assume that $G=G'^{\ad}$
 (note that $(\pi_*\mathcal{F}_n)^\ad=\pi^\ad_*\mathcal{F}_n^\ad$ by \cite[Theorem 3.7.2]{MR1734903}). 
 Since the Hecke action is transitive on the connected components of $\Sh_K(G,X)_{\overline{E}_v}$, 
 we may work on a connected component $\Sh_K(G,X)^0_{\overline{E}_v}$
 of $\Sh_K(G,X)_{\overline{E}_v}$ which is a quotient of
 a connected component $\Sh_{K'}(G',X')_{\overline{E}_v}^0$ of $\Sh_{K'}(G',X')_{\overline{E}_v}$ 
 by a free action of a finite group $H$ for some $K'$. 
 By Corollary \ref{cor:Hecke-inverse-disjoint}, 
 the inverse image of $C_{[P],\overline{E}_v}\cap \Sh_K(G,X)_{\overline{E}_v}^{0,\ad}$
 for $[P]\in\mathcal{P}_K(G)$ under
 \[
  f\colon \Sh_{K'}(G',X')_{\overline{E}_v}^{0,\ad}\to \Sh_K(G,X)^{0,\ad}_{\overline{E}_v}
 \]
 equals 
 \[
  \coprod_{[P']\in \mathcal{P}_{K'}(G'),[P']\mapsto [P]}C_{[P'],\overline{E}_v}\cap \Sh_{K'}(G',X')_{\overline{E}_v}^{0,\ad}.
 \]
 Since $H^i_c(C_{[P'],\overline{E}_v},f^*\mathcal{F})$ is finite-dimensional,
 so is 
 \[
  H^i_c\bigl(f^{-1}(C_{[P],\overline{E}_v}\cap \Sh_K(G,X)_{\overline{E}_v}^{0,\ad}),f^*\mathcal{F}\bigr).
 \]
 By Lemma \ref{lem:coh-invariant-part}, $H^i_c(C_{[P],\overline{E}_v}\cap \Sh_K(G,X)_{\overline{E}_v}^{0,\ad},\mathcal{F})$
 is equal to the $H$-invariant part of the above, hence finite-dimensional.
 This concludes the proof of \ref{item:coh-fin-l-adic}.

 The assertion \ref{item:coh-fin-torsion} can be proved in the same way,
 by using the Hochschild-Serre spectral sequence in place of Lemma \ref{lem:coh-invariant-part}
 when taking a quotient.
\end{prf}

\begin{rem}
 By the same method and Remark \ref{rem:limit-invariant}, we can also prove that the natural map
 $H^i_c(C_{[P],\overline{E}_v},\mathcal{F}^\ad)\to \varprojlim_{n} H^i_c(C_{[P],\overline{E}_v},\mathcal{F}_n^\ad)\otimes \overline{\Q}_\ell$ is an isomorphism. However, we do not need this fact.
\end{rem}

Now let $K^{p'}$ be as in Section \ref{subsec:main-thm}, $K_{p'}$ a compact open subgroup of $G(\Q_{p'})$,
and $K=K_{p'}K^{p'}$.
Take a Shimura datum $(G',X')$ of Hodge type such that $(G^\ad,X^\ad)\cong(G'^\ad,X'^\ad)$,
and an embedding $i\colon (G',X')\hookrightarrow (\GSp_{2n},X_{2n})$ into a Siegel Shimura datum.
Then, as in Lemma \ref{lem:C_P-disjoint} and Lemma \ref{lem:U_r}, we obtain an increasing sequence
\[
  \{[G]\}=\mathcal{P}_G(K)_{\le 0}\subset\mathcal{P}_G(K)_{\le 1}\subset\cdots\subset\mathcal{P}_G(K)_{\le n}=\mathcal{P}_G(K)
\]
of subsets of $\mathcal{P}_G(K)$. Note that $\mathcal{P}_G(K)_r=\mathcal{P}_G(K)_{\le r}\setminus \mathcal{P}_G(K)_{\le r-1}$ is a union of fibers of the natural map $\mathcal{P}_G(K)\to \mathcal{P}_{G,\A^\infty}$.
Therefore, by refining the sequence above, we can find an increasing sequence
\[
 \{[G]\}=S_0\subsetneq S_1\subsetneq \cdots \subsetneq S_m =\mathcal{P}_G(K)
\]
of subsets of $\mathcal{P}_G(K)$ satisfying the following conditions: 
\begin{itemize}
 \item For every $0\le r\le n$, there exists $0\le j\le m$ such that $S_j=\mathcal{P}_G(K)_{\le r}$.
 \item For $[P_1],[P_2]\in\mathcal{P}_G(K)$, $P_1$ and $P_2$ are conjugate by $G(\Q_{p'})$
       if and only if $[P_1], [P_2] \in S_j \setminus S_{j-1}$ for some $0\le j\le m$
       (here we put $S_{-1}=\varnothing$).
\end{itemize}

For $0\le j\le m$, we put $T_{j,K}=\bigcup_{[P]\in S_j}C_{[P]}$. By Lemma \ref{lem:U_r},
it is a constructible open subset of $\Sh^\ad_{K,E_v}$.
Further, we put $Z_{j,K}=T_{j,K}\setminus T_{j-1,K}$ ($T_{-1,K}$ is regarded as $\varnothing$).
Lemma \ref{lem:U_r} tells us that $Z_{j,K}=\coprod_{[P]\in S_j\setminus S_{j-1}}C_{[P]}$
as topological spaces.
For a compact open subgroup $K'_{p'} \subset K_{p'}$, we put $K'=K'_{p'} K^{p'}$, and 
let $S'_j$ be the inverse image of $S_j$ under the natural map $\mathcal{P}_G(K') \to \mathcal{P}_G(K)$. 
The sequence $\{S'_j\}_{0\le j\le m}$ satisfies the same conditions as $\{S_j\}_{0\le j\le m}$ does.
Therefore, we can define $T_{j,K'}$ and $Z_{j,K'}$ in the same way as $T_{j,K}$ and $Z_{j,K}$.
Note that $T_{j,K'}$ (resp.\ $Z_{j,K'}$) is the inverse image of $T_{j,K}$ (resp.\ $Z_{j,K}$)
under $\Sh^\ad_{K',E_v}\to \Sh^\ad_{K,E_v}$.
We put 
\[
 V^i_{\le j}=\varinjlim_{K'_{p'}} H_c^i(T_{j,K',\overline{E}_v},\mathcal{L}_\xi^\ad), \quad 
 V^i_{j}=\varinjlim_{K'_{p'}} H_c^i 
(Z_{j,K',\overline{E}_v},\mathcal{L}_\xi^\ad), 
\]
where $K'_{p'}$ runs through compact open subgroups of $K_{p'}$. 

\begin{lem}\label{lem:V_i}
 \begin{enumerate}
  \item\label{item:V_i-adm} The group $G(\Q_{p'})$ naturally acts on $V^i_{\le j}$ and $V^i_{j}$, 
	and these are admissible $G(\Q_{p'})$-representations.
  \item\label{item:V_i-exact-seq}
       We have the following long exact sequence of $G(\Q_{p'})$-modules:
       \[
       \cdots\to V^i_{\le j-1}\to V^i_{\le j}\to V^i_{j}\to V^{i+1}_{\le j-1}\to \cdots.
       \]
 \end{enumerate}
\end{lem}

\begin{prf}
 First, by Corollary \ref{cor:Hecke-inverse-disjoint}, the group $G(\Q_{p'})$ acts on $V^i_{\le j}$ and $V^i_j$.
 By Lemma \ref{lem:excision-exact}, we have a long exact sequence as in \ref{item:V_i-exact-seq},
 which is obviously $G(\Q_{p'})$-equivariant. 

 Let us prove \ref{item:V_i-adm}. Clearly $V^i_{\le j}$ and $V^i_j$ are smooth $G(\Q_{p'})$-representations.
 We will show the admissibility of them.
 Take a compact open subgroup $K_{p'}$ of $G(\Q_{p'})$ and its open normal subgroup $K'_{p'}$,
 and put $K=K_{p'}K^{p'}$, $K'=K'_{p'}K^{p'}$.
 Then, we have 
 \[
 H^i_c(Z_{j,K',\overline{E}_v},\mathcal{L}_\xi^\ad)^{K_{p'}}=H^i_c(Z_{j,K,\overline{E}_v},\mathcal{L}_\xi^\ad) 
 \]
 by Lemma \ref{lem:coh-invariant-part}. 
 Taking the inductive limit with respect to $K'_{p'}$, we have
 \[
  (V^i_{j})^{K_{p'}}=H^i_c(Z_{j,K,\overline{E}_v},\mathcal{L}_\xi^\ad)
 =\bigoplus_{[P]\in S_j\setminus S_{j-1}}H^i_c(C_{[P],\overline{E}_v},\mathcal{L}_\xi^\ad).
 \]
 By Proposition \ref{prop:finiteness-cohomology} \ref{item:coh-fin-l-adic},
 it is a finite-dimensional $\overline{\Q}_\ell$-vector space.
 Therefore we conclude that $V^i_j$ is an admissible representation of $G(\Q_{p'})$.
 By the long exact sequence in \ref{item:V_i-exact-seq} and the obvious identity $V^i_0=V^i_{\le 0}$,
 we can see the admissibility of $V^i_{\le j}$ inductively.
\end{prf}

\begin{prop}\label{prop:parabolically-induced}
 We take a representative $P_j$ of an element of $S_j\setminus S_{j-1}$. 
 For a compact open subgroup $K'_{p'} \subset K_{p'}$, 
 let $[P_j]_{K'}$ denote the class of $P_j$ in $\mathcal{P}_G (K')$. 
 We put $W^i_{j}=\varinjlim_{K'_{p'}}H^i_c (C_{[P_j]_{K'},\overline{E}_v},\mathcal{L}_\xi^\ad)$, 
 where $K'_{p'}$ runs through compact open subgroups of $K_{p'}$. 
\begin{enumerate}
 \item We have a natural smooth action of $P_j(\Q_{p'})$ on $W^i_{j}$.
 \item We have a natural $G(\Q_{p'})$-equivariant isomorphism 
	$V^i_j\cong \Ind_{P_j(\Q_{p'})}^{G(\Q_{p'})}W^i_{j}$. 
 \item The $G(\Q_{p'})$-representation $V^i_j$ is non-cuspidal.
 \end{enumerate}
\end{prop}

\begin{prf}
 The claim (i) is clear from Corollary \ref{cor:Hecke-inverse-disjoint}. 
 Let us prove (ii).
 We follow the proof of \cite[Proposition 5.20]{RZ-GSp4}. By the Frobenius reciprocity, 
 we have a homomorphism of $G(\Q_{p'})$-modules $\Ind_{P_j(\Q_{p'})}^{G(\Q_{p'})}W^i_{j}\to V^i_j$. 
 We shall observe that this is bijective. 
 We take a special maximal compact subgroup $K_{p'}^0$ of $G(\Q_{p'})$. 
 For a compact open subgroup $K'_{p'} \subset K_{p'}$ which is normal in $K_{p'}^0$, we have
 \begin{align*}
  H_c^i (Z_{j,K',\overline{E}_v},\mathcal{L}_\xi^\ad)
  &\stackrel{(1)}{=}\bigoplus_{[P] \in S'_j \setminus S'_{j-1}} 
  H^i_c(C_{[P],\overline{E}_v},\mathcal{L}_\xi^\ad)\\
  &\stackrel{(2)}{=} \bigoplus_{g\in K'_{p'}\backslash G(\Q_{p'})/P_j(\Q_{p'})}  H^i_c (C_{[gP_jg^{-1}]_{K'},\overline{E}_v},\mathcal{L}_\xi^\ad)\\
  &\stackrel{(3)}{=} 
 \bigoplus_{g\in K'_{p'} \backslash K_{p'}^0/P_j(\Q_{p'})\cap K_{p'}^0} 
 H^i_c (C_{[gP_jg^{-1}]_{K'},\overline{E}_v},\mathcal{L}_\xi^\ad)\\
  &\stackrel{(4)}{\cong} 
 \Ind_{(P_j(\Q_{p'})\cap K_{p'}^0)/(P_j(\Q_{p'})\cap K'_{p'})}^{K_{p'}^0/K'_{p'}}
 H^i_c (C_{[P_j]_{K'},\overline{E}_v},\mathcal{L}_\xi^\ad).
 \end{align*}
 Here (1) follows from $Z_{j,K'}=\coprod_{[P] \in S'_j \setminus S'_{j-1}}C_{[P]}$ mentioned before, 
 (2) from the definitions of $S'_j$ and $P_j$, 
 and (3) from the Iwasawa decomposition $G(\Q_{p'})=P_j(\Q_{p'})K_{p'}^0$. 
 The isomorphism (4) is a consequence of Corollary \ref{cor:Hecke-inverse-disjoint}
 and \cite[Lemme 13.2]{MR1719811}.
 By taking the inductive limit, we obtain $K_{p'}^0$-isomorphisms
 \[
  V^i_{j}\cong \Ind_{P_j(\Q_{p'})\cap K_{p'}^0}^{K_{p'}^0}
  W^i_{j}\xleftarrow{\cong}\Ind_{P_j(\Q_{p'})}^{G(\Q_{p'})}W^i_{j} 
 \]
 (the second map is an isomorphism by the Iwasawa decomposition).
 By the proof of \cite[Lemme 13.2]{MR1719811},
 it is easy to see that the first isomorphism above is nothing but the $K_{p'}^0$-homomorphism obtained by
 the Frobenius reciprocity for 
 $P_j(\Q_{p'})\cap K_{p'}^0\subset K_{p'}^0$. 
 Therefore the composite of the two isomorphisms
 above coincides with the $G(\Q_{p'})$-homomorphism introduced at the beginning of our proof of (ii).
 Thus we conclude the proof of (ii).

 Finally consider (iii). By (ii), we have only to prove that the unipotent radical of $P_j(\Q_{p'})$
 acts trivially on $W^i_{j}$. By \cite[Lemme 13.2.3]{MR1719811}, it suffices to prove that
 $W^i_{j}$ is an admissible $P_j(\Q_{p'})$-representation. 
 For any compact open subgroup $K'_{p'}$ of $K_{p'}^0$, 
 the vector space $(W^i_{j})^{P_j(\Q_{p'})\cap K'_{p'}}$ 
 is a subspace of 
 $(\Ind_{P_j(\Q_{p'})}^{G(\Q_{p'})}W^i_{j})^{K'_{p'}}$.
 By (ii) and Lemma \ref{lem:V_i} \ref{item:V_i-adm},
 $(\Ind_{P_j(\Q_{p'})}^{G(\Q_{p'})}W^i_{j})^{K'_{p'}}\cong (V^i_j)^{K'_{p'}}$ 
 is a finite-dimensional $\overline{\Q}_\ell$-vector space.
 Hence $W^i_{j}$ is an admissible $P_j(\Q_{p'})$-representation, as desired.
\end{prf}

\begin{prf}[of Theorem \ref{thm:main-thm}]
The claim follows from Lemma \ref{lem:V_i} \ref{item:V_i-exact-seq}
and Proposition \ref{prop:parabolically-induced} (iii), because 
$V_{\le 0}^i = H^i_c(\Sh^{\pg}_{\infty,K^{p'},\overline{E}_v},\mathcal{L}_\xi^{\ad})$ and 
$V_{\le m}^i =H^i_c(\Sh_{\infty,K^{p'},\overline{E}_v},\mathcal{L}_\xi)$.
\end{prf}

\subsection{Torsion coefficients}
Since our method of proving Theorem \ref{thm:main-thm} is totally geometric,
we may also obtain an analogous result for $\ell$-torsion coefficients.
For simplicity, we will only consider a constant coefficient 
$\overline{\F}_\ell$.
We assume that $p' \neq \ell$, and put
\begin{align*}
 H^i_c(\Sh_{\infty,K^{p'},\overline{E}_v},\overline{\F}_\ell)
 &=\varinjlim_{K_{p'}}
 H^i_c(\Sh_{K_{p'}K^{p'},\overline{E}_v},\overline{\F}_\ell),\\
 H^i_c(\Sh^{\pg}_{\infty,K^{p'},\overline{E}_v},\overline{\F}_\ell)
 &=\varinjlim_{K_{p'}}
 H^i_c(\Sh^{\pg}_{K_{p'}K^{p'},\overline{E}_v},\overline{\F}_\ell).
\end{align*}
They are naturally endowed with actions of 
$G(\Q_{p'})\times \Gal(\overline{E}_v/E_v)$.
They are admissible/continuous 
$G(\Q_{p'})\times \Gal(\overline{E}_v/E_v)$-representations; 
note that we have
\begin{align*}
 H^i_c(\Sh_{\infty,K^{p'},\overline{E}_v},\overline{\F}_\ell)^{K_{p'}} 
 &=H^i_c(\Sh_{K_{p'}K^{p'},\overline{E}_v},\overline{\F}_\ell),\\
 H^i_c(\Sh^{\pg}_{\infty,K^{p'},\overline{E}_v},\overline{\F}_\ell)^{K_{p'}}
 &=H^i_c(\Sh^{\pg}_{K_{p'}K^{p'},\overline{E}_v},\overline{\F}_\ell),
\end{align*}
if $K_{p'}$ is a pro-$p'$ group (\cf \cite[Proposition 2.5]{non-cusp}). 

The following theorem can be proved in exactly the same way as Theorem \ref{thm:main-thm}
(we use Proposition \ref{prop:finiteness-cohomology} \ref{item:coh-fin-torsion} in place of
Proposition \ref{prop:finiteness-cohomology} \ref{item:coh-fin-l-adic}).

\begin{thm}\label{thm:main-thm-mod-l}
We assume that $p' \neq \ell$. 
  The kernel and the cokernel of the canonical homomorphism
 \[
  H^i_c(\Sh^{\pg}_{\infty,K^{p'},\overline{E}_v},\overline{\F}_\ell)\to 
 H^i_c(\Sh_{\infty,K^{p'},\overline{E}_v},\overline{\F}_\ell)
 \]
 have no supercuspidal subquotient of $G(\Q_{p'})$. 
 For the definition of supercuspidal representations over 
$\overline{\F}_\ell$, see \cite[II.2.5]{MR1395151}.
\end{thm}

\section{PEL type case}\label{sec:PELcase}
\subsection{Notation for Shimura varieties of PEL type}\label{subsec:notation-Sh-var}
In this section, we are interested in Shimura varieties of PEL type considered in \cite[\S 5]{MR1124982}
(see also \cite[\S 1.4]{Kai-Wen}).
We recall it briefly. Fix a prime $p$. Consider a 6-tuple $(B,\mathcal{O}_B,*,V,L,\langle\ ,\ \rangle)$, where
\begin{itemize}
 \item $B$ is a finite-dimensional simple $\Q$-algebra such that $B\otimes_\Q\Q_p$ is a product of matrix algebras over unramified extensions of $\Q_p$,
 \item $\mathcal{O}_B$ is an order of $B$ whose $p$-adic completion is a maximal order of $B\otimes_\Q\Q_p$,
 \item $*$ is a positive involution of $B$ (namely, an involution such that $\Tr(bb^*)>0$ for every non-zero 
       $b\in B$) which preserves $\mathcal{O}_B$,
 \item $V$ is a non-zero finite $B$-module,
 \item $L$ is a $\Z$-lattice of $V$ preserved by $\mathcal{O}_B$, and
 \item $\langle\ ,\ \rangle\colon V\times V\to \Q$ is a non-degenerate alternating
       $*$-Hermitian pairing with respect to the $B$-action such that $\langle x,y\rangle\in\Z$
       for every $x,y\in L$, and that $L_p=L\otimes_\Z\Z_p$ is a self-dual lattice of $V_p=V\otimes_\Q\Q_p$.
\end{itemize}
From $(B,V,\langle\ ,\ \rangle)$, we define a simple $\Q$-algebra $C=\End_B(V)$ with a unique involution $\#$
satisfying $\langle cv,w\rangle=\langle v,c^\# w\rangle$ for every $c\in C$ and $v,w\in V$.
Moreover we define an algebraic group $G$ over $\Q$ by
\[
 G(R)=\{g\in (C\otimes_\Q R)^\times\mid gg^\#\in R^\times\}
\]
for every $\Q$-algebra $R$. The condition $gg^\#\in R^\times$ is equivalent to the existence of $c(g)\in R^\times$ such that
$\langle gv,gw\rangle=c(g)\langle v,w\rangle$ for every $v,w\in V\otimes_\Q R$.
By the presence of the lattice $L$, $G$ can be naturally extended to a group scheme over $\Z$, which is also
denoted by the same symbol $G$. 

Consider an $\R$-algebra homomorphism $h\colon \C\to C\otimes_\Q\R$ preserving involutions (on $\C$, we consider the complex conjugation)
such that the symmetric real-valued bilinear form $(v,w)\mapsto \langle v,h(i)w\rangle$ on $V\otimes_\Q\R$ is positive definite.
Such a 7-tuple $(B,\mathcal{O}_B,*,V,L,\langle\ ,\ \rangle,h)$ is said to be an unramified integral PEL datum.
Note that the map $h$ induces a homomorphism $\Res_{\C/\R}\mathbb{G}_m\to G_\R$ of algebraic groups over $\R$, which is also denoted by $h$.

Let $F$ be the center of $B$ and $F^+$ the subfield of $F$ consisting of elements fixed by $*$.
The existence of $h$ tells us that $N=[F:F^+](\dim_FC)^{1/2}/2$ is an integer.
An unramified integral PEL datum falls into the following three types:
\begin{description}
 \item[type (A)] $[F:F^+]=2$.
 \item[type (C)] $[F:F^+]=1$ and $C\otimes_\Q\R$ is isomorphic to a product of $M_{2N}(\R)$.
 \item[type (D)] $[F:F^+]=1$ and $C\otimes_\Q\R$ is isomorphic to a product of $M_N(\mathbb{H})$.
\end{description}
For simplicity, we will exclude the type (D) case.

Using $h\colon \C\to C\otimes_{\Q}\R\hookrightarrow C\otimes_{\Q}\C$,
we can decompose the $B\otimes_\Q\C$-module
$V\otimes_\Q\C$ as $V\otimes_\Q\C=V_1\oplus V_2$, where $V_1$ (resp.\ $V_2$) is the subspace of $V\otimes_\Q\C$
on which $h(z)$ acts by $z$ (resp.\ $\overline{z}$) for every $z\in\C$.
We denote by $E$ the field of definition of the isomorphism class of the $B\otimes_\Q\C$-module $V_1$,
and call it the reflex field. It is a subfield of $\C$ which is finite over $\Q$.

In the sequel, we fix an unramified integral PEL datum $(B,\mathcal{O}_B,*,V,L,\langle\ ,\ \rangle,h)$.
For a compact open subgroup $K^p$ of $G(\widehat{\Z}^p)$, consider the functor $\mathscr{S}_{K^p}$ from the category of
$\mathcal{O}_{E,(p)}$-schemes to the category of sets, that associates $S$ to the set of isomorphism classes
of quadruples $(A,i,\lambda,\eta^p)$, where
\begin{itemize}
 \item $A$ is an abelian scheme over $S$,
 \item $\lambda\colon A\to A^\vee$ is a prime-to-$p$ polarization,
 \item $i\colon \mathcal{O}_B\to \End_S(A)$ is an algebra homomorphism
       such that $\lambda\circ i(b)=i(b^*)^\vee\circ \lambda$ for every $b\in\mathcal{O}_B$,
 \item $\eta^p$ is a level-$K^p$ structure of $(A,i,\lambda)$ of type 
       $(L\otimes_\Z\widehat{\Z}^p,\langle\ ,\ \rangle)$ in the sense of \cite[Definition 1.3.7.6]{Kai-Wen},
\end{itemize}
satisfying the equality of polynomials $\det_{\mathcal{O}_S}(b;\Lie A)=\det_E(b;V_1)$ in the sense of
\cite[\S 5]{MR1124982}. Recall that two quadruples $(A,i,\lambda,\eta^p)$ and $(A',i',\lambda',\eta'^p)$
are said to be isomorphic if there exists an isomorphism $f\colon A\to A'$ of abelian schemes
such that 
\begin{itemize}
 \item $\lambda=f^\vee\circ \lambda'\circ f$,
 \item $f\circ i(b)=i'(b)\circ f$ for every $b\in\mathcal{O}_B$,
 \item and $f\circ \eta^p=\eta'^p$ in the sense of \cite[Definition 1.4.1.4]{Kai-Wen}.
\end{itemize}
If $K^p$ is neat (\cf \cite[Definition 1.4.1.8]{Kai-Wen}), the functor $\mathscr{S}_{K^p}$ is represented by a quasi-projective smooth $\mathcal{O}_{E,(p)}$-scheme (see \cite[Corollary 7.2.3.9]{Kai-Wen}), 
which is also denoted by $\mathscr{S}_{K^p}$.
Here we will call it a Shimura variety of PEL type.
The group $G(\A^{\infty,p})$ naturally acts on the tower of schemes $(\mathscr{S}_{K^p})_{K^p\subset G(\widehat{\Z}^p)}$
as Hecke correspondences.

Let $\ell$ be a prime number different from $p$.
For an algebraic representation $\xi$ of $G$ on a finite-dimensional $\overline{\Q}_\ell$-vector space, 
we can define a $\overline{\Q}_\ell$-sheaf $\mathcal{L}_\xi$ on $\Sh_{K}$
(see \cite[III, \S 6]{MR1044823}).
It is equivariant with respect to the Hecke action.

\begin{rem}\label{rem:PELSh}
 \begin{enumerate}
  \item Our definition of $\mathscr{S}_{K^p}$, due to \cite{Kai-Wen}, is slightly different from that in \cite{MR1124982},
	but they give the same moduli space. See \cite[Proposition 1.4.3.4]{Kai-Wen}.
  \item \label{item:PELShSh}
	Let us recall the relation between $\mathscr{S}_{K^p}$ and Shimura varieties in	Section \ref{sec:Shvar}.
	See \cite[\S 8]{MR1124982} for detail.
	Let $X$ denote the $G(\R)$-orbit of the homomorphism $h\colon \Res_{\C/\R}\mathbb{G}_m\to G_\R$.
	Then, the pair $(G,X)$ forms a Shimura datum, and 
	$\mathscr{S}_{K^p,E}$ is isomorphic to a disjoint union of $\#\ker^1(\Q,G)$ copies of
	$\Sh_{G(\Z_p)K^p}(G,X)$. In the cases of type (C) or type (A) with $N$ even,
	it is known that $\ker^1(\Q,G)=1$.
 \end{enumerate}
\end{rem}

So far in this section, we have only considered level structures which are prime to $p$. Now we add $p^m$-level structures
on the universal abelian scheme of the generic fiber $\mathscr{S}_{K^p,E}$.
Let $\mathscr{S}_{m,K^p,E}$ be the scheme over $\mathscr{S}_{K^p,E}$ classifying principal level-$m$ structures
(\cf \cite[Definition 1.3.6.2]{Kai-Wen})
of the universal object $(\mathcal{A},i^{\mathrm{univ}},\lambda^{\mathrm{univ}})$ over $\mathscr{S}_{K^p,E}$.
We denote the structure morphism $\mathscr{S}_{m,K^p,E}\to \mathscr{S}_{K^p,E}$ by $\pr_m$,
which is finite and \'etale.
We write $\mathcal{L}_{\xi,m}$ or $\mathcal{L}_\xi$ for the inverse image of $\mathcal{L}_\xi$ by $\pr_m$.

Let $K_{p,m}$ be the compact open subgroup of $G(\Q_p)$ defined as the kernel of $G(\Z_p)\to G(\Z/p^m\Z)$. 
Then $\mathscr{S}_{m,K^p,E}$ coincides with a disjoint union of 
the Shimura variety $\Sh_{K_{p,m}K^p}(G,X)$, where we use the notation in Remark \ref{rem:PELSh} (ii).

\subsection{Compactly supported cohomology and nearby cycle cohomology}\label{subsec:cohomology-notation}
Fix a place $v$ of $E$ over $p$. We write $E_v$ for the completion of $E$ at $v$, $\mathcal{O}_v$ for
the ring of integers of $E_v$, and $\kappa_v$ the residue field of $\mathcal{O}_v$.
We put $\mathscr{S}_{K^p,\eta}=\mathscr{S}_{K^p,E_v}$,
$\mathscr{S}_{K^p,\overline{\eta}}=\mathscr{S}_{K^p,\overline{E}_v}$,
$\mathscr{S}_{K^p,v}=\mathscr{S}_{K^p,\kappa_v}$ and
$\mathscr{S}_{K^p,\overline{v}}=\mathscr{S}_{K^p,\overline{\kappa}_v}$.
Further, for $m\ge 0$ we set $\mathscr{S}_{m,K^p,\eta}=\mathscr{S}_{m,K^p,E}\otimes_EE_v$ and
$\mathscr{S}_{m,K^p,\overline{\eta}}=\mathscr{S}_{m,K^p,E}\otimes_E\overline{E}_v$.

Let $p'$ be a prime number, and $K^{p'} \subset G(\widehat{\mathbb{Z}}^{p'})$ a compact open subgroup. 
If $p'\neq p$, we assume that $K^{p'}=K_{p,m_0}K^{p',p}$ for some $m_0\ge 0$ and compact open subgroup $K^{p',p}$
of $G(\widehat{\Z}^{p',p})$.
We put
\[
 H^i_c(\mathscr{S}_{\infty,K^{p'},\overline{\eta}},\mathcal{L}_\xi)=
 \begin{cases}
 \varinjlim_{m}H^i_c(\mathscr{S}_{m,K^p,\overline{\eta}},\mathcal{L}_\xi) & \textrm{if $p'=p$},\\ 
 \varinjlim_{K_{p'}}H^i_c(\mathscr{S}_{m_0,K_{p'}K^{p',p},\overline{\eta}},\mathcal{L}_\xi) & \textrm{if $p'\neq p$}, 
\end{cases}
\]
which is an admissible/continuous $G(\Q_{p'})\times \Gal(\overline{E}_v/E_v)$-representation.
We are also interested in the nearby cycle cohomology defined as follows:
\[
 H^i_c(\mathscr{S}_{\infty,K^{p'},\overline{v}},R\psi\mathcal{L}_\xi)=
 \begin{cases}
 \varinjlim_{m}H^i_c\bigl(\mathscr{S}_{K^p,\overline{v}},R\psi (\pr_{m*}\mathcal{L}_\xi)\bigr) & \textrm{if $p'=p$},\\ 
 \varinjlim_{K_{p'}}H^i_c\bigl(\mathscr{S}_{K_{p'} K^{p',p} ,\overline{v}},R\psi (\pr_{m_0 *}\mathcal{L}_\xi)\bigr) & \textrm{if $p'\neq p$}. 
\end{cases}
\]
Obviously the group $\Gal(\overline{E}_v/E_v)$ acts on it. The following lemma gives an action of
$G(\Q_{p'})$ on $H^i_c(\mathscr{S}_{\infty,K^{p'},\overline{v}},R\psi\mathcal{L}_\xi)$.

\begin{lem}\label{lem:nearby-gp-action}
 We have a natural action of $G(\Q_{p'})$ on $H^i_c(\mathscr{S}_{\infty,K^{p'},\overline{v}},R\psi\mathcal{L}_\xi)$. 
 By this action, $H^i_c(\mathscr{S}_{\infty,K^{p'},\overline{v}},R\psi\mathcal{L}_\xi)$ becomes an admissible/continuous
 $G(\Q_{p'})\times \Gal(\overline{E}_v/E_v)$-representation.
\end{lem}

\begin{prf}
We show the claim in the case $p'=p$. 
The other cases are easier. 
 To ease notation, we omit the subscript $K^p$.

 As in \cite[\S 6]{MR2169874},
 we can construct a tower $(\mathscr{S}_m)_{m\ge 0}$ of schemes over $\mathcal{O}_v$ with finite transition maps
 such that $\mathscr{S}_m$ gives an integral model of $\mathscr{S}_{m,\eta}$ and $\mathscr{S}_0=\mathscr{S}$.
 In this situation, we have 
 \[
 H^i_c\bigl(\mathscr{S}_{\overline{v}},R\psi (\pr_{m*}\mathcal{L}_\xi)\bigr)=H^i_c(\mathscr{S}_{m,\overline{v}},R\psi\mathcal{L}_\xi),
 \]
 where $\mathscr{S}_{m,\overline{v}}=\mathscr{S}_m\otimes_{\mathcal{O}_v}\overline{\kappa}_v$.

 We put $G^+(\Q_p)=\{g\in G(\Q_p)\mid g^{-1}L_p\subset L_p\}$. For $g\in G^+(\Q_p)$, let
 $e(g)$ be the minimal non-negative integer such that $gL_p\subset p^{-e(g)}L_p$.
 Then we can construct a tower $(\mathscr{S}_{m,g})_{m\ge e(g)}$ of schemes over $\mathcal{O}_v$
 and two morphisms
 \[
 \pr\colon \mathscr{S}_{m,g}\to \mathscr{S}_m,\quad [g]\colon \mathscr{S}_{m,g}\to \mathscr{S}_{m-e(g)}
 \]
 which are compatible with the transition maps. It is known that these are proper morphisms,
 $\pr$ induces an isomorphism on the generic fibers,
 and $[g]$ induces the Hecke action of $g$ on the generic fibers
 (see \cite[Proposition 16, Proposition 17]{MR2169874}). In particular, we have a canonical cohomological
 correspondence (\cf \cite[Expos\'e III]{SGA5}, \cite{MR1431137})
 \[
  c_g\colon [g]_\eta^*\mathcal{L}_{\xi,m-e(g)}\xrightarrow{\cong} \pr_\eta^*\mathcal{L}_{\xi,m}=R\pr_\eta^!\mathcal{L}_{\xi,m}.
 \]
 Let
 \[
  R\psi(c_g)\colon [g]_{\overline{v}}^*R\psi\mathcal{L}_{\xi,m-e(g)}\to R\pr_{\overline{v}}^!R\psi\mathcal{L}_{\xi,m}
 \]
 be the specialization of $c_g$ (\cf \cite[\S 1.5]{MR1431137}, \cite[\S 6]{RZ-GSp4}).
 Since $[g]_{\overline{v}}$ is proper, this induces a homomorphism
 \[
  H^i_c(\mathscr{S}_{m-e(g),\overline{v}},R\psi\mathcal{L}_\xi)
 \xrightarrow{H_c^i(R\psi(c_g))}H^i_c(\mathscr{S}_{m,\overline{v}},R\psi\mathcal{L}_\xi).
 \]
 Taking the inductive limit, we get
 \[
  \gamma_g\colon H^i_c(\mathscr{S}_{\infty,K^p,\overline{v}},R\psi\mathcal{L}_\xi)
 \to H^i_c(\mathscr{S}_{\infty,K^p,\overline{v}},R\psi\mathcal{L}_\xi).
 \]
 From an obvious relation $c_{gg'}=c_g\circ g^*c_{g'}$ for $g,g'\in G^+(\Q_p)$,
 we deduce $\gamma_{gg'}=\gamma_{g}\circ \gamma_{g'}$ (\cf \cite[Corollary 6.3]{RZ-GSp4}).
 On the other hand, by \cite[Proposition 16 (3), Proposition 17 (3)]{MR2169874}, $\gamma_{p^{-1}}$ is
 the identity. Since $G(\Q_p)$ is generated by $G^+(\Q_p)$ and $p$ as a monoid, we can extend $\gamma_g$
 to the whole $G(\Q_p)$. By \cite[Proposition 16 (4)]{MR2169874}, the restriction of this action
 to $K_{p,0}=G(\Z_p)$ coincides with the inductive limit of the natural action of $K_{p,0}$ on
 $H^i_c(\mathscr{S}_{m,\overline{v}},R\psi \mathcal{L}_\xi)$.
 In particular, it is a smooth action. Furthermore, for integers $m'\ge m\ge 1$, we have
 \[
  H^i_c(\mathscr{S}_{m',\overline{v}},R\psi \mathcal{L}_\xi)^{K_{p,m}/K_{p,m'}}=H^i_c(\mathscr{S}_{m,\overline{v}},R\psi \mathcal{L}_\xi)
 \]
 (see \cite[Proposition 2.5]{non-cusp}). Taking inductive limit, we obtain
 \[
  H^i_c(\mathscr{S}_{\infty,K^p,\overline{v}},R\psi\mathcal{L}_\xi)^{K_{p,m}}=H^i_c(\mathscr{S}_{m,\overline{v}},R\psi \mathcal{L}_\xi).
 \]
 This implies that $H^i_c(\mathscr{S}_{\infty,K^p,\overline{v}},R\psi\mathcal{L}_\xi)$ is
 an admissible/continuous representation of $G(\Q_p)\times\Gal(\overline{E}_v/E_v)$.
\end{prf}

\begin{cor}\label{cor:PEL-comp}
 The kernel and the cokernel of the canonical homomorphism
 \[
  H^i_c(\mathscr{S}_{\infty,K^{p'},\overline{v}},R\psi\mathcal{L}_\xi)\to H^i_c(\mathscr{S}_{\infty,K^{p'},\overline{\eta}},\mathcal{L}_\xi)
 \]
 (\cf \cite[Expos\'e XIII, (2.1.7.3)]{SGA7}) are non-cuspidal.
 In particular, for an irreducible supercuspidal representation 
 $\pi$ of $G(\Q_{p'})$,
 we have an isomorphism 
 \[
  H^i_c(\mathscr{S}_{\infty,K^{p'},\overline{\eta}},\mathcal{L}_\xi)[\pi]\cong H^i_c(\mathscr{S}_{\infty,K^{p'},\overline{v}},R\psi\mathcal{L}_\xi)[\pi].
 \]
 Similar results hold for the coefficient $\overline{\F}_\ell$ for a prime number $\ell\neq p,p'$.
\end{cor}

\begin{prf}
Analogues of Theorem \ref{thm:Sh-partition} and Theorem \ref{thm:main-thm} are also valid in 
the PEL type case in this section by Remark \ref{rem:PELSh} \ref{item:PELShSh}. 
Let 
$\mathscr{S}_{K^{p},E_v}^{\pg}$ 
be the potentially good reduction locus of 
$\mathscr{S}_{K^{p},E_v}^{\ad}$. 
As in the Siegel case, 
$\mathscr{S}_{K^{p},E_v}^{\pg}$ coincides with 
the rigid generic fiber of the completion of 
$\mathscr{S}_{K^{p},\mathcal{O}_v}$ along the special fiber. 
Hence, we have an isomorphism 
\[
 H^i_c(\mathscr{S}_{K^{p},\overline{E}_v}^{\pg}, \pr_{m *}\mathcal{L}_\xi^{\ad}) 
 \cong 
 H^i_c(\mathscr{S}_{K^{p},\overline{v}},R\psi (\pr_{m *} \mathcal{L}_\xi) ) 
\]
for any non-negative integer $m$
by \cite[Theorem 3.7.2, Theorem 5.7.6]{MR1734903} and \cite[Theorem 3.1]{MR1626021}.
Taking inductive limits, we have an isomorphism 
\[
 H^i_c(\mathscr{S}_{\infty,K^{p'},\overline{E}_v}^{\pg}, \mathcal{L}_\xi^{\ad})  \cong 
 H^i_c(\mathscr{S}_{\infty,K^{p'},\overline{v}},R\psi\mathcal{L}_\xi). 
\]
Hence the claim follows from the analogue of 
Theorem \ref{thm:main-thm}. 
\end{prf}

\begin{rem}
 \begin{enumerate}
  \item The case where $p'\neq p$ in Corollary \ref{cor:PEL-comp} 
	was previously obtained by Tetsushi Ito and the second author. 
	In that case, we can use minimal compactifications over $\mathcal{O}_v$ 
	to show the claim. 
  \item In \cite{LSNcycle}, Lan and Stroh obtained a stronger result that the canonical homomorphism
	in Corollary \ref{cor:PEL-comp} is in fact an isomorphism. Their method is totally different from ours.
 \end{enumerate}
\end{rem}

\subsection{Example}\label{subsec:example}
In this subsection, we give a very simple application of 
Corollary \ref{cor:PEL-comp}.
Proofs in this subsection are rather sketchy, since the technique is more or less well-known.

Here we consider the Shimura variety for $\mathit{GU}(1,n-1)$ over $\Q$.
Let $F$ be an imaginary quadratic extension of $\Q$ and $\Spl_{F/\Q}$ the set of rational primes over which $F/\Q$ splits. We fix a field embedding $F\hookrightarrow \C$ and regard $F$ as a subfield of $\C$.
For an integer $n\ge 2$, consider the integral PEL datum
$(B,\mathcal{O}_B,*,V,L,\langle\ ,\ \rangle,h)$ as follows:
\begin{itemize}
 \item $B=F$, $\mathcal{O}_B=\mathcal{O}_F$ and $*$ is the unique non-trivial element of $\Gal(F/\Q)$.
 \item $V=F^n$ and $L=\mathcal{O}_F^n$.
 \item $\langle\ ,\ \rangle\colon V\times V\to \Q$ is an alternating pairing
       satisfying the following conditions:
       \begin{itemize}
	\item[$\bullet$] $\langle x,y\rangle\in \Z$ for every $x,y\in L$, 
	\item[$\bullet$] $\langle bx,y\rangle=\langle x,b^*y\rangle$ for every $x,y\in V$ and $b\in F$, and
	\item[$\bullet$] $G_{\R}\cong \mathit{GU}(1,n-1)$ (for the definition of $G$, see Section \ref{subsec:notation-Sh-var}).
       \end{itemize}
 \item $h\colon \C\to \End_F(V)\otimes_\Q\R\cong M_n(\C)$ is given by
       $z\mapsto \begin{pmatrix}z&0\\0&\overline{z}I_{n-1}\end{pmatrix}$,
       where the last isomorphism is induced by the fixed embedding $F\hookrightarrow \C$.
\end{itemize}
In this case, the reflex field $E$ is equal to $F$. To a neat compact open subgroup $K$ of $G(\widehat{\Z})$,
we can attach the Shimura variety $\Sh_K$ of PEL type, which is not proper over $\Spec F$. 

Put $\Sigma=\{p\in\Spl_{F/\Q}\mid L_p=L_p^\perp\}$. Then our integral Shimura datum is
unramified at every $p\in\Sigma$.
Moreover, for such $p$, $G_{\Q_p}$ is isomorphic to $\GL_n(\Q_p)\times \GL_1(\Q_p)$
(\cf \cite[\S 1.2.3]{MR2074714}).
If $K=K_{p,0}K^p$ for some compact open subgroup
$K^p$ of $G(\widehat{\Z}^p)$, we have 
$\Sh_K=\mathscr{S}_{K^p}\otimes_{\mathcal{O}_{F,(p)}}F$, where
$\mathscr{S}_{K^p}$ is the moduli space 
introduced in Section \ref{subsec:notation-Sh-var}.

Let us fix a prime number $\ell$. We put 
\[
 H_c^i(\Sh,\overline{\Q}_\ell)=\varinjlim_{K}H_c^i(\Sh_K\otimes_F\overline{F},\overline{\Q}_\ell).
\]
It is an admissible/continuous $G(\A^\infty)\times\Gal(\overline{F}/F)$-representation over $\overline{\Q}_\ell$.

\begin{thm}\label{thm:non-cusp}
 Let $\Pi$ be an irreducible admissible representation of $G(\A^\infty)$ over $\overline{\Q}_\ell$.
 Assume that there exists a prime $p\in \Sigma$ such that $\Pi_p$ is a supercuspidal 
 representation of $G(\Q_p)$.
 Then $H_c^i(\Sh,\overline{\Q}_\ell)[\Pi]=0$ unless $i=n-1$.
\end{thm}

\begin{rem}
 For proper Shimura varieties, an analogous result is known (\cite{MR1114211}, \cite[Corollary IV.2.7]{MR1876802}).
 It would be possible to give an ``automorphic'' proof of Theorem \ref{thm:non-cusp}
 by using results in \cite{MR2567740}. However, the authors think that our proof,
 consisting of purely local arguments, is simpler and has some importance.
\end{rem}

\begin{prf}
 Let $\ell'$ be another prime number and fix an isomorphism of fields $\iota\colon\overline{\Q}_\ell\cong \overline{\Q}_{\ell'}$. Then $\iota$ induces an isomorphism $H_c^i(\Sh,\overline{\Q}_\ell)[\Pi]\cong H_c^i(\Sh,\overline{\Q}_{\ell'})[\iota \Pi]$, where $\iota \Pi$ is the representation of $G(\A^\infty)$ over $\overline{\Q}_{\ell'}$
 induced by $\Pi$ and $\iota$. It is easy to observe that $\Pi_p$ is supercuspidal if and only if
 $(\iota\Pi)_p$ is supercuspidal. Therefore, we can change our $\ell$ freely, and thus we can assume that
 there exists a prime $p\in \Sigma\setminus \{\ell\}$ such that $\Pi_p$ is supercuspidal.  
 Fix such $p$ and take a place $v$ of $F$ lying over $p$. Then, for an integer $m\ge 0$ and a neat compact open
 subgroup $K^p$ of $G(\widehat{\Z}^p)$, $\Sh_{K_{p,m}K^p}\otimes_FF_v$ is isomorphic to $\mathscr{S}_{m,K^p,\eta}$
 introduced in Section \ref{subsec:cohomology-notation}.
 Therefore we have an isomorphism
 \[
  H_c^i(\Sh,\overline{\Q}_\ell)\cong \varinjlim_{m,K^p}H_c^i(\mathscr{S}_{m,K^p,\overline{\eta}},\overline{\Q}_\ell)
 =\varinjlim_{K^p}H_c^i(\mathscr{S}_{\infty,K^p,\overline{\eta}},\overline{\Q}_\ell).
 \]
 Thus it suffices to show that $H_c^i(\mathscr{S}_{\infty,K^p,\overline{\eta}},\overline{\Q}_\ell)[\pi]=0$
 for a supercuspidal representation $\pi$ of $G(\Q_p)$, a neat compact open subgroup $K^p$, 
 and an integer $i\neq n-1$.
 By Corollary \ref{cor:PEL-comp}, it is equivalent to showing that
 $H_c^i(\mathscr{S}_{\infty,K^p,\overline{v}},R\psi\overline{\Q}_\ell)[\pi]=0$.

 For an integer $h\ge 0$, let $\mathscr{S}^{[h]}_{K^p,v}$ be the reduced closed subscheme of $\mathscr{S}_{K^p,v}$
 consisting of points $x$ such that the \'etale rank of $\mathcal{A}_x[v^\infty]$ is less than
 or equal to $h$ (\cf \cite[p.~111]{MR1876802}),
 where $\mathcal{A}$ denotes the universal abelian scheme over $\mathscr{S}_{K^p}$.
 Put $\mathscr{S}^{(h)}_{K^p,v}=\mathscr{S}^{[h]}_{K^p,v}\setminus \mathscr{S}^{[h-1]}_{K^p,v}$.
 Our proof of the theorem is divided into the subsequent two lemmas.
\end{prf} 

\begin{lem}
 For every supercuspidal representation $\pi$ of $G(\Q_p)$, 
 we have 
 \[
  H_c^i(\mathscr{S}_{\infty,K^p,\overline{v}},R\psi\overline{\Q}_\ell)[\pi]
 =\Bigl(\varinjlim_{m}H_c^i\bigl(\mathscr{S}^{[0]}_{K^p,\overline{v}},(R\psi \pr_{m*}\overline{\Q}_\ell)\vert_{\mathscr{S}^{[0]}_{K^p,\overline{v}}}\bigr)\Bigr)[\pi].
 \]
\end{lem}

\begin{prf}
 First recall that $\mathscr{S}_{m,K^p,\eta}$ has a good integral model over $\mathcal{O}_v$. 
 For an integer $m\ge 0$, consider the functor from the category of $\mathcal{O}_v$-schemes to the category of
 sets, that associates $S$ to the set of isomorphism classes of 6-tuples $(A,i,\lambda,\eta^p,\eta_v,\eta_{p,0})$,
 where
 \begin{itemize}
  \item $[(A,i,\lambda,\eta^p)]\in \mathscr{S}_{K^p}(S)$,
  \item $\eta_v\colon L\otimes_\Z(v^{-m}\mathcal{O}_v/\mathcal{O}_v)\to A[v^m]$ is a Drinfeld $v^m$-level structure
	(\cf \cite[II.2]{MR1876802}), and
  \item $\eta_{p,0}\colon p^{-m}\Z/\Z\to \mu_{p^m,S}$ is a Drinfeld $p^m$-level structure.
 \end{itemize}
 Then it is easy to see that this functor is represented by a scheme $\mathscr{S}_{m,K^p}$
 which is finite over $\mathscr{S}_{K^p}$. Moreover the generic fiber of $\mathscr{S}_{m,K^p}$ can be naturally identified with 
 $\mathscr{S}_{m,K^p,\eta}$ (\cf the moduli problem $\mathfrak{X}'_U$ introduced in
 \cite[p.~92]{MR1876802}). As in \cite[III.4]{MR1876802}, we can extend the Hecke action of $G(\Q_p)$ on
 $(\mathscr{S}_{m,K^p,\eta})_{m\ge 0}$ to the tower $(\mathscr{S}_{m,K^p})_{m\ge 0}$. We have a $G(\Q_p)$-equivariant isomorphism
 \[
  H_c^i(\mathscr{S}_{\infty,K^p,\overline{v}},R\psi\overline{\Q}_\ell)\cong \varinjlim_{m}H_c^i(\mathscr{S}_{m,K^p,\overline{v}},R\psi \overline{\Q}_\ell).
 \]
 Let us denote by $\mathscr{S}^{[h]}_{m,K^p,v}$ (resp.\ $\mathscr{S}^{(h)}_{m,K^p,v}$) the inverse image of
 $\mathscr{S}^{[h]}_{K^p,v}$ (resp.\ $\mathscr{S}^{(h)}_{K^p,v}$) under $\mathscr{S}_{m,K^p}\to \mathscr{S}_{K^p}$.
 For an integer $h\ge 0$, it is easy to observe that
 \begin{align*}
  \varinjlim_{m}H_c^i\bigl(\mathscr{S}^{[h]}_{K^p,\overline{v}},(R\psi \pr_{m*}\overline{\Q}_\ell)\vert_{\mathscr{S}^{[h]}_{K^p,\overline{v}}}\bigr)&\cong \varinjlim_{m}H_c^i\bigl(\mathscr{S}^{[h]}_{m,K^p,\overline{v}},(R\psi \overline{\Q}_\ell)\vert_{\mathscr{S}^{[h]}_{m,K^p,\overline{v}}}\bigr),\\
 \varinjlim_{m}H_c^i\bigl(\mathscr{S}^{(h)}_{K^p,\overline{v}},(R\psi \pr_{m*}\overline{\Q}_\ell)\vert_{\mathscr{S}^{(h)}_{K^p,\overline{v}}}\bigr)&\cong \varinjlim_{m}H_c^i\bigl(\mathscr{S}^{(h)}_{m,K^p,\overline{v}},(R\psi \overline{\Q}_\ell)\vert_{\mathscr{S}^{(h)}_{m,K^p,\overline{v}}}\bigr)
 \end{align*}
 and that they are admissible $G(\Q_p)$-representations. Moreover, by considering the kernel of 
 the universal Drinfeld $v^m$-level structure $\eta_v^{\mathrm{univ}}$, we can decompose $\mathscr{S}_{m,K^p,v}^{(h)}$ into finitely many
 open and closed subsets indexed by the set consisting of direct summands of $L\otimes_\Z(v^{-m}\mathcal{O}_v/\mathcal{O}_v)$ with
 rank $n-h$ (\cf \cite[D\'efinition 10.4.1, Proposition 10.4.2]{MR1719811} and \cite[Definition 5.1, Lemma 5.3]{RZ-GSp4}).
 Using this partition, when $h>0$, we can prove that the $G(\Q_p)$-representation 
 \[
  \varinjlim_{m}H_c^i\bigl(\mathscr{S}^{(h)}_{m,K^p,\overline{v}},(R\psi \overline{\Q}_\ell)\vert_{\mathscr{S}^{(h)}_{m,K^p,\overline{v}}}\bigr)
 \]
 is parabolically induced from a proper parabolic subgroup of $G(\Q_p)$.
 Therefore, by the same argument as in the proof of Theorem \ref{thm:main-thm},
 we can conclude that the kernel and the cokernel of
 \[
  \varinjlim_{m}H_c^i(\mathscr{S}_{m,K^p,\overline{v}},R\psi \overline{\Q}_\ell)\to \varinjlim_{m}H_c^i\bigl(\mathscr{S}^{[0]}_{m,K^p,\overline{v}},(R\psi \overline{\Q}_\ell)\vert_{\mathscr{S}^{[0]}_{m,K^p,\overline{v}}}\bigr)
 \]
 are non-cuspidal. This completes the proof of the lemma.
\end{prf}

\begin{lem}\label{lem:basic-noncusp}
 Let $\pi$ be a supercuspidal representation of $G(\Q_p)$.
 If $i\neq n-1$, we have
 \[
  \Bigl(\varinjlim_{m}H_c^i\bigl(\mathscr{S}^{[0]}_{K^p,\overline{v}},(R\psi \pr_{m*}\overline{\Q}_\ell)\vert_{\mathscr{S}^{[0]}_{K^p,\overline{v}}}\bigr)\Bigr)[\pi]=0.
 \]
\end{lem}

\begin{prf}
 Let $\mu_h\colon \mathbb{G}_{m,\C}\to G_{\C}$ be the homomorphism of algebraic groups over $\C$
 defined as the composite of
 \[
  \mathbb{G}_{m,\C}\xrightarrow{z\mapsto (z,1)}\mathbb{G}_{m,\C}\times \mathbb{G}_{m,\C}\stackrel{(*)}{\cong} (\Res_{\C/\R}\mathbb{G}_{m,\C})\otimes_\R\C\xrightarrow{h_\C}G_{\C},
 \]
 where $(*)$ is given by $\C\otimes_\R\C\xrightarrow{\cong}\C\times\C$; $a\otimes b\mapsto (ab,a\overline{b})$.
 Fix an isomorphism of fields $\C\cong \overline{\Q}_p$ and denote by $\mu\colon \mathbb{G}_{m,\overline{\Q}_p}\to G_{\overline{\Q}_p}$ the induced cocharacter of $G_{\overline{\Q}_p}$.
 Let $b$ be a unique basic element of $B(G_{\Q_p},\mu)$
 (for the definition of $B(G,\mu)$, we refer to \cite[\S 2.1.1]{MR2074714}),
 and denote by $\mathcal{M}$ the Rapoport-Zink space associated to the local unramified PEL datum
 $(F\otimes_{\Q}\Q_p,\mathcal{O}_F\otimes_{\Z}\Z_p,*,V_p,L_p,\langle\ ,\ \rangle,b,\mu)$
 (\cf \cite[\S 2.3.5]{MR2074714}). The Rapoport-Zink space $\mathcal{M}$ is equipped with an action of
 the group $J(\Q_p)$,
 where $J$ denotes the algebraic group over $\Q_p$ associated to $b$ (\cf \cite[Proposition 1.12]{MR1393439}).
 By \cite[\S 2.3.7.1]{MR2074714}, 
 $\mathcal{M}$ is isomorphic to $\mathcal{M}_{\mathrm{LT}}\times\Q_p^\times/\Z_p^\times$, where 
 $\mathcal{M}_{\mathrm{LT}}$ is the Lubin-Tate space for $\GL_n(\Q_p)$.
 Furthermore, $J(\Q_p)$ is isomorphic to $D^\times\times \Q_p^\times$, where $D$ denotes the central
 division algebra over $\Q_p$ with invariant $1/n$. The action of $J(\Q_p)$ on $\mathcal{M}$ is identified with
 the well-known action of $D^\times\times \Q_p^\times$ on $\mathcal{M}_{\mathrm{LT}}\times\Q_p^\times/\Z_p^\times$.

 By the $p$-adic uniformization theorem of Rapoport-Zink (\cite[Theorem 6.30]{MR1393439}, \cite[Corollaire 3.1.9]{MR2074714}), we have an isomorphism
 \[
 \coprod_{\ker^1(\Q,G)}I(\Q)\backslash \mathcal{M}\times G(\A^{\infty,p})/K^p\cong \mathscr{S}_{K^p}^{\wedge},
 \]
 where $I$ is an algebraic group over $\Q$ satisfying $I(\A^\infty)\cong J(\Q_p)\times G(\A^{\infty,p})$ and
 $\mathscr{S}_{K^p}^{\wedge}$ denotes the formal completion of $\mathscr{S}_{K^p}\otimes_{\mathcal{O}_v}W(\overline{\F}_p)$
 along $\mathscr{S}_{K^p,\overline{v}}^{[0]}$, the basic locus of $\mathscr{S}_{K^p,\overline{v}}$.
 By this isomorphism, we know that $\mathscr{S}^{[0]}_{K^p,\overline{v}}$, which coincides with $\mathscr{S}_{K^p}^{\wedge}$
 as topological spaces, consists of finitely many closed points; indeed, 
 the left hand side of the isomorphism above is a finite disjoint union of formal schemes of the form
 $\Gamma\backslash \mathcal{M}$, where $\Gamma\subset J(\Q_p)$ is a discrete cocompact subgroup
 (\cf \cite[Lemme 3.1.7]{MR2074714}). Therefore,  by \cite[Theorem 3.1]{MR1395723},
 we have an isomorphism
 \begin{align*}
 H_c^i\bigl(\mathscr{S}^{[0]}_{K^p,\overline{v}},(R\psi \pr_{m*}\overline{\Q}_\ell)\vert_{\mathscr{S}^{[0]}_{K^p,\overline{v}}}\bigr)
  &=H^i\bigl(\mathscr{S}^{[0]}_{K^p,\overline{v}},(R\psi \pr_{m*}\overline{\Q}_\ell)\vert_{\mathscr{S}^{[0]}_{K^p,\overline{v}}}\bigr)\\
  &\cong H^i\bigl(\mathscr{S}_{m,K^p,\overline{\eta}}(b),\overline{\Q}_\ell\bigr),
 \end{align*}
 where $\mathscr{S}_{m,K^p,\overline{\eta}}(b)=\pr_m^{-1}(\spp^{-1}(\mathscr{S}_{K^p,v}^{[0]})^\circ)_{\overline{F}_v}$.

 Now we use the Hochschild-Serre spectral sequence (see \cite[Th\'eor\`eme 4.5.12]{MR2074714})
\begin{align*}
 &E_2^{r,s}=\varinjlim_{m}\Ext^r_{J(\Q_p)\text{-smooth}}\bigl(H^{2(n-1)-s}_c(\mathcal{M}_{K_{p,m}},\overline{\Q}_\ell)(n-1),\mathcal{A}(I)_{\mathbf{1}}^{K^p}\bigr)\\
 &\qquad\qquad\qquad\Rightarrow \varinjlim_{m}H^{r+s}\bigl(\mathscr{S}_{m,K^p,\overline{\eta}}(b),\overline{\Q}_\ell\bigr),
\end{align*}
 where $\mathcal{M}_{K_{p,m}}$ is the Rapoport-Zink space of level $K_{p,m}$,
 and $\mathcal{A}(I)_{\mathbf{1}}$ is the space of automorphic forms on $I(\A^\infty)$
 (see \cite[D\'efinition 4.5.8]{MR2074714} for detail).
 Since $J(\Q_p)=D^\times\times \Q_p^\times$ is anisotropic modulo center, 
 it is easy to see that $E_2^{r,s}=0$ unless $r=0$. If $r=0$, we have
 \begin{align*}
  E_2^{0,s}&=\varinjlim_m\Hom_{J(\Q_p)}\bigl(H^{2(n-1)-s}_c(\mathcal{M}_\infty,\overline{\Q}_\ell)(n-1),\mathcal{A}(I)_{\mathbf{1}}^{K^p}\bigr)^{K_{p,m}},
 \end{align*}
 where we put $H^i_c(\mathcal{M}_\infty,\overline{\Q}_\ell)=\varinjlim_{m}H^i_c(\mathcal{M}_{K_{p,m}},\overline{\Q}_\ell)$.

 By \cite{non-cusp},
 the $G(\Q_p)$-representation $H^{2(n-1)-s}_c(\mathcal{M}_\infty,\overline{\Q}_\ell)(n-1)$ has
 non-zero supercuspidal part only if $s=n-1$. Indeed, for an irreducible supercuspidal representation
 $\pi=\pi_1\otimes \chi$ of $G(\Q_p)=\GL_n(\Q_p)\times \GL_1(\Q_p)$,
 where $\pi_1$ is an irreducible supercuspidal representation of $\GL_n(\Q_p)$
 and $\chi$ is a character of $\GL_1(\Q_p)$, we have
 \[
 H^i_c(\mathcal{M}_\infty,\overline{\Q}_\ell)[\pi]
 =H^i_c(\mathcal{M}_{\mathrm{LT},\infty},\overline{\Q}_\ell)[\pi_1]\otimes\chi,
 \]
 as we see in \cite[p.~168]{MR2074714}. Therefore $E_2^{0,s}$ has a supercuspidal subquotient only if
 $s=n-1$.

 Hence we conclude that
 \[
 \varinjlim_m H_c^i\bigl(\mathscr{S}^{[0]}_{K^p,\overline{v}},(R\psi \pr_{m*}\overline{\Q}_\ell)\vert_{\mathscr{S}^{[0]}_{K^p,\overline{v}}}\bigr)\cong \varinjlim_{m}H^i\bigl(\mathscr{S}_{m,K^p,\overline{\eta}}(b),\overline{\Q}_\ell\bigr)
 \]
 has non-zero supercuspidal part only if $i=n-1$.
\end{prf}

We also have a similar result for the torsion coefficient case.
For a neat compact open subgroup $K^p$ of $G(\widehat{\Z}^p)$, we put
\[
 H_c^i(\Sh_{K^p},\overline{\F}_\ell)=\varinjlim_{m}H_c^i(\Sh_{K_{p,m}K^p}\otimes_F\overline{F},\overline{\F}_\ell).
\]
It is an admissible/continuous $G(\Q_p)\times\Gal(\overline{F}/F)$-representation over $\overline{\F}_\ell$.

\begin{thm}\label{thm:non-cusp-torsion}
 Let $p$ be a prime in $\Sigma\setminus \{\ell\}$ and $\pi$ an irreducible supercuspidal $\overline{\F}_\ell$-representation of $G(\Q_p)$.
 Then, for every neat compact open subgroup $K^p$ of $G(\widehat{\Z}^p)$, we have
 $H_c^i(\Sh_{K^p},\overline{\F}_\ell)[\pi]=0$ unless $i=n-1$.
\end{thm}

\begin{rem}
 \begin{enumerate}
  \item Theorem \ref{thm:non-cusp-torsion} for proper Shimura varieties is due to Shin \cite{MR3377388}. 
	His method, using Mantovan's formula, is slightly different from ours.
	The non-proper cases are also covered in his paper using results in our paper. 
  \item Using the result in \cite{Dat-ltmodl}, it is possible to describe
	the action of $W_{\Q_p}$ on $H_c^{n-1}(\Sh_{K^p},\overline{\F}_\ell)[\pi]$
	by means of the mod-$\ell$ local Langlands correspondence.
	Such study has also been carried out by Shin when the Shimura variety is proper.
 \end{enumerate}
\end{rem}

\begin{prf}
 Almost all arguments in the proof of Theorem \ref{thm:non-cusp} work well.
 The only one point which should be modified is about the vanishing of the supercuspidal part of $E_2^{r,s}$
 for $(r,s)\neq (0,n-1)$ in the proof of Lemma \ref{lem:basic-noncusp}; note that an irreducible $\overline{\F}_\ell$-representation of $J(\Q_p)$,
 being supercuspidal, is not necessarily injective in the category of smooth 
 $\overline{\F}_\ell$-representations of $J(\Q_p)$ with the fixed central character.
 For this point, we can use the same argument as that by Shin (see \cite[\S 3.2]{MR3377388}), in which he uses
 the vanishing of the supercuspidal part 
 $H^i_c(\mathcal{M}_{\mathrm{LT},\infty},\overline{\F}_\ell)_{\mathrm{sc}}$ for $i\neq n-1$ 
 (\cf \cite[proof of Proposition 3.1.1, Remarque 3.1.5]{Dat-ltmodl}) and 
 the projectivity of the $D^\times$-representation 
 $H^{n-1}_c(\mathcal{M}_{\mathrm{LT},\infty},\overline{\F}_\ell)_{\mathrm{sc}}$
 (\cf \cite[\S 3.2.2, Remarque iii)]{Dat-ltmodl}).
\end{prf}

\def\cprime{$'$} \def\cprime{$'$} \newcommand{\dummy}[1]{}
\providecommand{\bysame}{\leavevmode\hbox to3em{\hrulefill}\thinspace}
\providecommand{\MR}{\relax\ifhmode\unskip\space\fi MR }
\providecommand{\MRhref}[2]{%
  \href{http://www.ams.org/mathscinet-getitem?mr=#1}{#2}
}
\providecommand{\href}[2]{#2}

\bigbreak\bigbreak

\noindent Naoki Imai\par
\noindent Graduate School of Mathematical Sciences, 
The University of Tokyo, 3--8--1 Komaba, Meguro-ku, 
Tokyo, 153--8914, Japan\par
\noindent E-mail address: \texttt{naoki@ms.u-tokyo.ac.jp}

\bigbreak

\noindent Yoichi Mieda\par
\noindent Graduate School of Mathematical Sciences, 
The University of Tokyo, 3--8--1 Komaba, Meguro-ku, 
Tokyo, 153--8914, Japan\par
\noindent E-mail address: \texttt{mieda@ms.u-tokyo.ac.jp}

\end{document}